\documentclass[a4paper,11pt]{elsarticle}
\usepackage[utf8]{inputenc}
\usepackage{bm}
\usepackage{amsmath}
\usepackage{mathrsfs}
\usepackage{graphicx}
\usepackage{epsfig, setspace}
\usepackage{url}
\usepackage{graphicx, transparent, color}

\usepackage{etex}
\reserveinserts{28}

\DeclareSymbolFont{rsfs}{U}{rsfs}{m}{n}
\DeclareSymbolFontAlphabet{\mathcal}{rsfs}

\usepackage[margin=3cm]{geometry}
\usepackage{tikz}
\usepackage{caption}
\usetikzlibrary{shapes,calc}
\usepackage{verbatim}

\numberwithin{equation}{section}

\usepackage{epsfig}
\usepackage{amsmath}
\usepackage{pst-plot}
\usepackage{amssymb}
\usepackage{amsfonts}
\usepackage{lmodern}
\usepackage{euscript}
\usepackage{eufrak}
\usepackage{oldgerm}
\usepackage{lineno}

\numberwithin{equation}{section}

\newtheorem{remark}{\bf Remark}[section]

\newtheorem{example}{Example}[section]

\newtheorem{theorem}{Theorem}[section]

\begin{document}
\begin{frontmatter}
\title{ {Simple smoothness indicator and multi-level adaptive order WENO scheme for hyperbolic conservation laws}}

\date{}
% \author{Rakesh Kumar}
\author{
{\sc Rakesh Kumar} \thanks{ Email: rakeshiitb21@gmail.com, rakesh@tifrbng.res.in}, 
and {\sc Praveen Chandrashekar} \thanks{Corresponding Author. Email: praveen@tifrbng.res.in} \\[2pt]
Tata Institute of Fundamental Research\\
Centre For Applicable Mathematics \\
Bangalore-560065, India.
}
% \shortauthorlist{Rakesh \emph{et al.}}
 \journal{}
% \fntext[fn1]{ Email: rakeshiitb21@gmail.com, rakesh@tifrbng.res.in}
% 
\begin{abstract}
 In the present work, we propose two new variants of fifth order finite difference WENO schemes of adaptive order. 
We compare our proposed schemes with other variants of WENO  schemes with special emphasize on WENO-AO(5,3) scheme 
[Balsara, Garain, and Shu, {\it J. Comput. Phys.},  326 (2016), pp 780-804].
 The first algorithm (WENO-AON(5,3)), involves the construction of a new simple smoothness indicator which reduces the computational cost of WENO-AO(5,3) 
 scheme. Numerical experiments show that accuracy of WENO-AON(5,3) scheme is comparable to that of WENO-AO(5,3) scheme and 
 resolution of solutions involving shock or other discontinuities is comparable to that of WENO-AO(5,3) scheme. 
The second algorithm denoted as
 WENO-AO(5,4,3), involves
 the inclusion of an extra cubic polynomial reconstruction in the base WENO-AO(5,3) scheme, which leads to a more accurate 
scheme. Extensive 
 numerical experiments in 1D and 2D are performed, which shows that WENO-AO(5,4,3) scheme has better resolution near shocks or discontinuities
 among the considered WENO schemes with negligible increase in computational cost. 
\end{abstract}
\end{frontmatter}
\section{Introduction}
The aim of the present article is to develop efficient and accurate fifth order WENO schemes of adaptive order for the hyperbolic conservation laws given by
\begin{equation}\label{hyp.1}
 \frac{\partial \bold{u}}{\partial t}+\displaystyle{\sum_{\alpha=1}^d \frac{\partial \bold{f}_{\alpha}(\bold{u})}{\partial x_{\alpha}} =0 ,    }
\end{equation}
where $\bold{u}=(u_1,u_2,\ldots, u_m)\in \mathbb{R}^m$ are the conserved variables and $\bold{f}_{\alpha}:\mathbb{R}^m \rightarrow \mathbb{R}^m$,
 $\alpha=1,2,\ldots,d$ are the Cartesian components of flux. It is well-known that the classical solution of  \eqref{hyp.1} may cease to exist in 
 finite time, even if the initial data is
sufficiently smooth. The appearance of shocks, contact discontinuities and rarefaction waves in the solution profile makes it difficult to devise 
 stable and high order accurate numerical
schemes due to the development of spurious oscillations or even growing numerical instabilities. 

WENO schemes are one of the most successful higher order accurate schemes, which computes the solution accurately while maintaining high resolution
near the discontinuities in a non-oscillatory fashion. Initially, WENO schemes were constructed in finite volume framework by Liu  et al.
 \cite{liu-etal_94a}. They  have constructed a finite volume WENO scheme,  in which they choose a
convex combination of the reconstructions on all possible stencils used in  essentially non-oscillatory (ENO) schemes \cite{har-etal_97a}. The 
resulting scheme achieves one order more than ENO scheme.  Jiang and Shu \cite{jia-shu_96a} constructed the  fifth-order finite difference 
versions of WENO scheme
in higher dimensions for hyperbolic conservation laws. They also  designed a general framework for the smoothness indicators and the nonlinear
weights. Higher order extension of finite difference WENO can be found in  \cite{bal-shu_00a}, in which Balsara and Shu have devised  WENO schemes upto
 eleventh order. Further, Gerolymos et al. \cite{ger-etal_09a}  developed  very high-order WENO schemes upto the seventeenth order.  The WENO schemes
 were also developed in central framework  \cite{lev-etal_02a}, \cite{sam-etal_16a}, \cite{cra-sem_16a} and its extension to unstructured
 meshes can be found in  \cite{dum-kas_07a}, \cite{fri-98a}.
  WENO schemes have been successfully applied to many 
problems having both shocks and complicated smooth structures, for more details we refer \cite{shu_09a} and references therein.

In central WENO  \cite{lev-etal_99a, lev-etal_02a},  Levy et. al discussed the reconstruction problem where the linear weights may not exist. They provide a new way of assigning the linear 
weights to the stencils and devised the stable central WENO schemes. In \cite{lev-etal_02a},  Levy et. al developed a new third-order,
compact central WENO reconstruction, which is written as a convex combination of a quadratic and two linear polynomials.
Zhu and Qui \cite{zhu-qiu_16a} proposed a new fifth order finite difference WENO scheme, named as WENO-ZQ, for hyperbolic conservation laws. The
WENO-ZQ scheme is a  convex combination of a fourth degree polynomial and two linear polynomials. The scheme is shown to be more accurate 
as compared to 
classical WENO-JS scheme in $L^{\infty}$- and $L^1$- norms. Further, this new formulation is extended to finite volume framework \cite{zhu-qiu_16b} and 
 for solving Hamiton-Jacobi equation \cite{zhu-qui_16c}. Balsara et al. \cite{bal-etal_16a} proposed a new class of WENO 
schemes of adaptive order for hyperbolic conservation laws, denoted as WENO-AO. In \cite{bal-etal_16a}, fifth order WENO scheme of
adaptive order, denoted  as
WENO-AO(5,3), is a convex combination of three quadratic polynomials and a fourth degree polynomial. Balsara et al. \cite{bal-etal_16a} proposed a novel way to
compute smoothness indicators by expressing reconstruction polynomials in terms of Legendre polynomials.
 Even though WENO schemes of adaptive order are more accurate, it involves the calculation of an extra smoothness indicator over the largest stencils in 
comparison to classical WENO scheme. Huang and Chen \cite{hua-che_18a} showed that computational cost of WENO-AO(5,3) can be reduced
 by evaluating the smoothness indicator of bigger stencil as a  combination of smoothness indicators of the smaller 
 stencils. Further, numerical experiments were presented for fifth order WENO schemes of adaptive order with new smoothness indicators on 
  uniform and non-uniform meshes with demonstrable reduction in computational cost in comparison to WENO-AO(5,3) scheme. A similar idea is 
used by the same authors in developing efficient upwind and central WENO scheme of sixth order \cite{hua-che_18b}.

In this article, we have proposed two new variants of fifth order WENO schemes of adaptive order for hyperbolic conservation laws. In the first variant 
of WENO scheme of adaptive order, denoted as WENO-AON(5,3), we devised a simple  smoothness indicator defined over the smaller stencils,
that can  serve as
a smoothness indicator for the bigger stencil. The proposed simple smoothness indicator is closer to the original smoothness indicator (see \cite{bal-etal_16a})
  at the truncation error level as compared to existing indicators in literature and takes less time to compute solution. In devising the second variant of WENO scheme, our main concern is
 the accuracy and resolution of solutions
 near the shock or discontinuities. This variant uses a  convex combination of three quadratic polynomials, 
 one cubic polynomial and a fourth degree polynomial
 and the resulting scheme is named as WENO-AO(5,4,3). In fifth order WENO scheme of adaptive order given by \cite{bal-etal_16a} and \cite{hua-che_18a},
  solution is computed with quadratic polynomial, if the bigger stencil has  shocks or discontinuities. Whereas in the case of WENO-AO(5,4,3),
  once the discontinuities are detected over the bigger stencil, instead of computing solution directly using lower degree 
 quadratic polynomial, it may use the central approximation given by cubic polynomial. Hence WENO-AO(5,4,3) leads to a more accurate scheme 
 as compared to 
 other WENO schemes of order five in the presence of shocks or discontinuities. Using the approach discussed in \cite{Arb-etal_17a} (see also \cite{ar-etal_11a}, \cite{kol_14a} 
 for more details), we have shown
 that for both solutions, nonlinear weights approach the linear weights if solution is smooth over the bigger stencil. We have also shown that both types of 
 flux
 reconstruction
  are fifth order accurate under the assumption of smoothness of solution over the bigger stencil.

Further, we  compare the two proposed algorithms of WENO scheme of adaptive order with other fifth order WENO schemes present in the
literature. Here we consider the fifth order  WENO-JS \cite{jia-shu_96a}, WENO-Z \cite{bor-etal_08a}, WENO-AO(5,3)\cite{bal-etal_16a}, and WENO-ZQ
\cite{zhu-qiu_16a} schemes for comparison. Since WENO-AO(5,4,3) scheme contains the 
  additional computation of extra cubic polynomials, which leads to an increment in the computational cost, it is still interesting to see 
  resolution of 
  solutions computed with WENO-AO(5,4,3) scheme near shocks or discontinuities.
  The performance of proposed algorithms  are studied through numerical experiments. We observed that WENO-AON(5,3) scheme computes comparable 
   solution to WENO-AO(5,3) scheme with less computational cost. We also found WENO-AO(5,4,3) scheme resolves the solution near discontinuities 
   more accurately than the other considered WENO schemes with negligible increment in computational cost for one and two dimensional problems. 
   
The article is organized as follows: in section \ref{WENOAO}, we have detailed  the basic formulation of WENO scheme of adaptive order and discussed the 
construction of  WENO-AON(5,3) and WENO-AO(5,4,3) schemes.
The analysis of nonlinear weights and flux reconstruction has been discussed in Section \ref{SI}.
In the section \ref{nm}, The proposed algorithms are compared with other 
variants of WENO scheme for smooth test cases or test cases where solution has shocks or other discontinuities.
    At the end, conclusions are drawn in section \ref{ccl}.

%----------------Formulation of WENO scheme-------------------------------------
\section{Finite difference WENO schemes}\label{WENOAO}
In this section, we discuss the implementation of two new fifth order  WENO schemes of adaptive order for one-dimensional scalar hyperbolic conservation laws given by
\begin{equation}\label{hyp.scalar}
\left\{\begin{array}{ll}
                u_t+f(u)_x=0, ~~~~(x,t)\in[a,b]\times (0,T], \\
                u(x,0)=u_0(x).
                \end{array}\right.
\end{equation}
For simplicity, we distribute the computational domain into 
smaller cells $\mathbb{I}_i=[x_{i-\frac{1}{2}}, x_{i+\frac{1}{2}}]$ having a uniform mesh width $\Delta x= x_{i+\frac{1}{2}}-x_{i-\frac{1}{2}}$ with
$x_i=\frac{1}{2}{(x_{i-\frac{1}{2}}+x_{i+\frac{1}{2}}})$ denoting the center of the cell. 
The semi-discrete formulation of \eqref{hyp.scalar} over the cell 
$\mathbb{I}_i$
 can be given by
 \begin{equation}
  \frac{d u_i(t)}{dt}=-\frac{1}{\Delta x}\Big(\mathbb{F}_{i+\frac{1}{2}}-\mathbb{F}_{i-\frac{1}{2}}\Big),~~~~i=0,1,\ldots,N,
 \end{equation}
 where $u_i(t)=u(x_i,t)$ is the point value and $\mathbb{F}_{i+\frac{1}{2}}$ is numerical flux defined at the interface $x_{i+\frac{1}{2}}$. The 
 conservation property is assured by defining a function $\mathbb{H}(x)$ implicitly through the following equation \cite{shu-osh_89a} 
 \begin{equation}
  f(u(x))=\frac{1}{\Delta x}\int_{x-\frac{\Delta x}{2}}^{x+\frac{\Delta x}{2}} \mathbb{H}(\eta)d\eta.
 \end{equation}
 Differentiating with respect to $x$, we get
 \begin{equation}
  f(u)_x|_{x=x_i}=\frac{1}{\Delta x}\Big(\mathbb{H}(x_{i+\frac{1}{2}})-\mathbb{H}(x_{i-\frac{1}{2}})\Big),
 \end{equation}
where $\mathbb{H}(x_{i\pm\frac{1}{2}})$ is a approximation to the numerical flux $\mathbb{F}_{i\pm\frac{1}{2}}$ with an $r^{\mbox{th}}$ order of accuracy in the
 sense that
 \begin{equation}
 \mathbb{F}_{i\pm\frac{1}{2}}=\mathbb{H}(x_{i\pm\frac{1}{2}})+O(\Delta x^r).
 \end{equation}
 To ensure the numerical stability, the flux $f(u)$ is split into positive and negative parts as follows
\begin{equation}
 f(u)=f^{+}(u)+f^{-}(u),
\end{equation}
such that $\frac{df^{+}(u)}{du}\geq 0$ and $\frac{df^{-}(u)}{du}\leq 0$. There are many choices of splitting. Here we consider the Lax-Friedrich 
splitting given by
\begin{equation}
 f^{\pm}(u)=\frac{1}{2}(f(u)\pm\lambda u),
\end{equation}
where $\lambda =\max(|f'(u)|)$, such that maximum is taken over whole range of $u$ in the mesh. In case of the system of equations, $\lambda$ is the 
maximum of the absolute values of the largest eigenvalue of the Jacobian matrix  taken over $x$-axis. The numerical flux is also divided into positive
and negative fluxes, which
are obtained from positive and negative parts of $f(u)$, such that
\begin{equation}
 \mathbb{F}_{i\pm\frac{1}{2}}=\mathbb{F}^{+}_{i\pm\frac{1}{2}}+\mathbb{F}^{-}_{i\pm\frac{1}{2}}.
\end{equation}
Here, we define the WENO reconstruction procedure for positive part, since negative part is symmetric to the positive part with respect to
$x_{i+\frac{1}{2}}$. We drop the $'+'$ sign to avoid confusion.

\subsection{WENO-AON(5,3)}
The reconstruction algorithm of WENO scheme of adaptive order is designed to approximate the spatial derivative of the flux $f$, which involves 
the following steps
\begin{itemize}
 \item Step-1: We choose the bigger stencil $\mathbb{S}^5_0=\{{I}_{i-2},{I}_{i-1},{I}_{i},{I}_{i+1},{I}_{i+2}\}$
 and reconstruct a quartic polynomial $\mathbb{P}^4_0(x)$ satisfying
 \begin{equation}
  \frac{1}{\Delta x}\int_{x_{j-\frac{1}{2}}}^{x_{j+\frac{1}{2}}} \mathbb{P}_0^4(x)dx= f(u_j), ~~~~j=i-2,i-1,i,i+1,i+2.
 \end{equation}
 Similar to WENO-JS scheme \cite{jia-shu_96a}, we also consider three smaller stencils $$\mathbb{S}_k^3=\{I_{i+k-1},I_{i+k},I_{i+k+1}\},~~k=-1,0,1,$$
 and reconstruct quadratic polynomials
  $\mathbb{P}^2_k(x),~\mbox{for}~k=-1,0,1$ over each stencil satisfying
   \begin{equation}
  \frac{1}{\Delta x}\int_{x_{j-\frac{1}{2}}}^{x_{j+\frac{1}{2}}} \mathbb{P}^2_k(x)dx= f(u_j), ~~~~j=i+k-1,i+k,i+k+1,~~k=-1,0,1.
 \end{equation}

 \item Step-2: Define the fourth order polynomial $\mathbb{P}(x)$ as follows \cite{bal-etal_16a}
 \begin{equation}
 \mathbb{P}( x)=\frac{1}{\gamma^5_0}\Big[\mathbb{P}^4_0( x)-\displaystyle{\sum_{k=-1}^{1} \gamma_{k}^3
 \mathbb{P}_k^2( x)\Big]},
\end{equation}
where $\gamma_{-1}^3,~\gamma_0^3,~\gamma_1^3, ~\gamma_0^5$ are the positive linear weights corresponding to stencils $\mathbb{S}_{-1}^3,~\mathbb{S}_{0}^3
,~\mathbb{S}_{1}^3,\mathbb{S}_{0}^5$, respectively, and  satisfying the relation
\begin{equation}
 \gamma_{-1}^3+\gamma_{0}^3+\gamma_1^3+\gamma_0^5=1.
\end{equation}
Here, we choose the 
linear weights similar to 
\cite{bal-etal_16a}, as follows
\begin{align*}
%  &\gamma_{Hi} \in (0, 0.8], ~~~\gamma_{Lo}=1-\gamma_{Hi}\\
 \gamma_0^5 = \gamma_{Hi},~
 \gamma_{-1}^3 = (1-\gamma_{Hi})(1-\gamma_{Lo})/2,~
\gamma_0^3 = (1-\gamma_{Hi})\gamma_{Lo},~
 \gamma_1^3 =  (1-\gamma_{Hi})(1-\gamma_{Lo})/2.
\end{align*}
The range of values of $\gamma_{Hi}\in[0.85, 0.95]$ and $\gamma_{Lo}\in[0.85, 0.95]$.

\item Step-3:  Compute the smoothness indicators denoted by $\beta_k^m$, which measure the smoothness of the function (see \cite{jia-shu_96a} for 
more details) as given by
\begin{equation}
 \beta_k^m=\displaystyle{\sum_{l=1}^r \int_{x_{i-\frac{1}{2}}}^{x_{i+\frac{1}{2}}}\Delta x^{2l-1}\Big(\frac{d^l\mathbb{P}_k^m(x)}{dx^l} \Big)^2 dx},
 ~~ (k,m)=(-1,3), ~(0,3),~(1,3),~\mbox{and}~(0,5).
\end{equation}
In the quadratic case, smoothness indicators are given by
 \begin{align*}
  \beta_{-1}^3&=\frac{13}{12}(f_{i-2}-2f_{i-1}+f_{i})^2+\frac{1}{4}(f_{i-2}-4f_{i-1}+3f_{i})^2,\\
  \beta_0^3&=\frac{13}{12}(f_{i-1}-2f_{i}+f_{i+1})^2+\frac{1}{4}(f_{i-1}-f_{i+1})^2,\\
  \beta_1^3&=\frac{13}{12}(f_{i}-2f_{i+1}+f_{i+2})^2+\frac{1}{4}(3f_{i}-4f_{i+1}+f_{i+2})^2.
 \end{align*}
 The smoothness indicator over the stencil $\mathbb{S}^5_0$ is given in compact form as follows (see \cite{bal-etal_16a})
 \begin{equation}
  (\beta^5_0)^{(1)}=\Big(u_{51}+\frac{u_{53}}{10}\Big)^2+ \frac{13}{3}\Big(u_{52}+\frac{123}{455}u_{54}\Big)^2
  +\frac{781}{20}(u_{53})^2+\frac{1421461}{2275}(u_{54})^2,
 \end{equation}
 where
 \begin{equation}
 \left. 
\begin{array}{ll}
u_{51}&=\frac{1}{120}\Big(11 f_{i-2}-82f_{i-1}+82f_{i+1}-11f_{i+2}\Big),\\
u_{52}&=\frac{1}{56}\Big(-3 f_{i-2}+40 f_{i-1}-74f_i+40 f_{i+1}-3f_{i+2}\Big),\\
u_{53}&=\frac{1}{12}\Big(-f_{i-2}+2f_{i-1}-2f_{i+1}+f_{i+2}\Big),\\
u_{54}&=\frac{1}{24}\Big( f_{i+2}-4 f_{i-1}+6f_i-4 f_{i+1}+f_{i+2}\Big).
\end{array}
\right \}
\end{equation}
The superscript used for $\beta^5_0$ is to distinguish the numerical schemes. Here, $(\beta^5_0)^{(1)}$ correspond to smoothness indicators 
used in WENO-AO(5,3) scheme  
(see \cite{bal-etal_16a}).  Despite  being the more accurate scheme as compared to  WENO-JS and other WENO variants, WENO-AO(5,3) scheme involves
the computation of smoothness
indicator over the bigger stencil $\mathbb{S}_0^5$, 
which increases the computational cost of the scheme especially for systems of equations in multiple dimensions. In \cite{hua-che_18a}, Huang and Chen observed  that the smoothness indicator of bigger stencil
can be expressed as a linear 
combination of smoothness indicators defined over smaller stencils and they proposed the following smoothness indicator
\begin{equation}\label{SI:HC}
 (\beta^5_0)^{(2)}=\displaystyle{\frac{\epsilon+\beta_{-1}^3}{3\epsilon +\sum_{k=-1}^1 \beta_k^3} \beta_{-1}^3+
 \frac{\epsilon+\beta_0^3}{3\epsilon +\sum_{k=-1}^1 \beta_k^3} \beta_0^3+\frac{\epsilon+\beta_1^3}{3\epsilon +\sum_{k=-1}^1 \beta_k^3} \beta_1^3}
\end{equation}
For $\epsilon=0$, Taylor series expansion of $(\beta^5_0)^{(2)}$ at point $x_{i+\frac{1}{2}}$  shows that the first two terms  of $(\beta^5_0)^{(2)}$ 
coincide with that of $(\beta^5_0)^{(1)}$.

Further, in this article, we have proposed a new smoothness indicator denoted as $(\beta^5_0)^{(3)}$ for $\mathbb{S}^5_0$, which is constructed on an 
idea  similar to \cite{hua-che_18a} but is more efficient (shown numerically in subsequent sections)
\begin{equation}\label{SI:RP}
 (\beta^5_0)^{(3)}=\frac{1}{6}(\beta_{-1}^3+4\beta_0^3+\beta_1^3)+|\beta_{-1}^3-\beta_1^3|.
\end{equation}
An analysis of this smoothness indicator is given in later section.
\item Step-4: To avoid the loss of accuracy at critical points, we adopt the strategy discussed for  WENO-Z  schemes in
\cite{bor-etal_08a,cas-etal_11a,don-bor_13a}. We define
\begin{equation}
 \tau=\left(\frac{|(\beta^5_0)^{(l)}-\beta_{-1}^3|+|(\beta^5_0)^{(l)}-\beta_0^3|+|(\beta^5_0)^{(l)}-\beta_1^3|}{3}\right).
\end{equation}
We will get different reconstruction schemes by choosing the value of $l$. The value $l=1$ correspond to WENO-AO(5,3) \cite{bal-etal_16a}, $l=2$ 
gives WENO-AO-HC\cite{hua-che_18a}, and for $l=3$,
 we obtain a new reconstruction  of WENO  scheme named as WENO-AON(5,3).

The nonlinear weights can be defined using 
\begin{align*}
 \tilde{\omega}_{-1}^3&=\gamma_{-1}^3\Big(1+\frac{\tau^2}{(\epsilon+\beta_{-1}^3)^2}\Big),~~~~\tilde{\omega}_0^3=\gamma_0^3\Big(1+\frac{\tau^2}{(\epsilon+\beta_0^3)^2}\Big),\\
 \tilde{\omega}_1^3&=\gamma_1^3\Big(1+\frac{\tau^2}{(\epsilon+\beta_1^3)^2}\Big),~~~~\tilde{\omega}_0^5=\gamma_0^5\Big(1+\frac{\tau^2}{(\epsilon+(\beta_0^5)^{(l)})^2}\Big).
\end{align*}
where $\epsilon$ is a small positive real number to avoid the denominator to become zero. In numerical experiments, we take $\epsilon=10^{-12}$.
The normalized nonlinear weights are given by
\begin{equation}\label{nlw33444}
 \omega_k^l = \frac{\gamma_{k}^l\Big(1+\frac{\tau^2}{(\beta_k^l+\epsilon)^2}\Big)}{\displaystyle{\sum_{j=-1}^1 \gamma_j^3
 \Big(1+\frac{\tau^2}{(\beta_j^3+\epsilon)^2}\Big)+  \gamma_0^5\Big(1+\frac{\tau^2}{(\beta_0^5+\epsilon)^2}\Big)}}~~\mbox{for}~~(l,k)=\{(3,-1),
 (3,0), (3,1), (5,0)\}.
\end{equation}
\item Step-5: Finally, the nonlinear reconstruction is defined as
\begin{equation}
 \mathbb{F}_{i+\frac{1}{2}}=\displaystyle{ \omega_0^5 \mathbb{P}( x_{i+\frac{1}{2}})+
 \sum_{k=-1}^{1} \omega_k^3 \mathbb{P}_{k}^2( x_{i+\frac{1}{2}})   }.
\end{equation}
Similarly, the negative part can be calculated.
\end{itemize}

\subsection{WENO-AO(5,4,3)}
In this section, we describe the WENO-AO(5,4,3) reconstruction for the positive part of the flux.  The  reconstruction procedure  of WENO-AO(5,4,3) is 
similar
 to WENO-AO(5,3) except for the addition of an extra stencil $\mathbb{S}^4_0$. The WENO-AO(5,4,3) scheme involves the following steps
\begin{itemize}
\item Step-1: Apart from the stencils $\mathbb{S}_{-1}^3$, $\mathbb{S}_{0}^3$, $\mathbb{S}^3_{1}$, and $\mathbb{S}^5_0$, we consider another central 
stencil $\mathbb{S}^4_0$,
which contains four intervals, given by 
\begin{align}
  \mathbb{S}_0^4&=\{I_{i-1},I_i,I_{i+1},I_{i+2}\}.
\end{align}
Let $\mathbb{P}^3_0$ denote the cubic polynomial defined over the stencil $\mathbb{S}^4_0$ satisfying
 \begin{equation}
  \frac{1}{\Delta x}\int_{x_{j-\frac{1}{2}}}^{x_{j+\frac{1}{2}}} \mathbb{P}^3_0(x)dx= f(u_j), ~~~~j=i-1,i,i+1,i+2.
 \end{equation}
\item Step-2: (Linear Weights) The linear weights are assigned to stencils as follows
  \begin{equation}\label{lw:543}
 \left. 
\begin{array}{ll}
\gamma_{-1}^3&=\frac{1}{2}(1-\gamma_{Hi})(1-\gamma_{Avg})(1-\gamma_{Lo}),\\
 \gamma_0^3&=(1-\gamma_{Hi})(1-\gamma_{Avg})\gamma_{Lo},\\
 \gamma_1^3&=\frac{1}{2}(1-\gamma_{Hi})(1-\gamma_{Avg})(1-\gamma_{Lo}),\\
 \gamma_0^4&=(1-\gamma_{Hi})\gamma_{Avg},\\
 \gamma_0^5&=\gamma_{Hi}.
\end{array}
\right \}
\end{equation}
where $\gamma_0^4$ denotes the linear positive weights corresponding to the stencils $\mathbb{S}_0^4$, satisfying $\gamma_{-1}^3+\gamma_0^3+\gamma_1^3
+\gamma_0^4
+\gamma_0^5=1$. The $\gamma_{Hi}$,
$\gamma_{Avg}$, and $\gamma_{Lo}$ are positive numbers and $\gamma_{Hi}, \gamma_{Avg},\gamma_{Lo} \in (0,1)$.

Then the linear reconstruction can be defined as
\begin{align}
 \mathbb{P}(x)=\frac{1}{\gamma_0^5}\Big(\mathbb{P}^4_0(x)-\gamma_{-1}^3 \mathbb{P}^2_{-1}(x)-\gamma_{0}^3 \mathbb{P}^2_0(x)
 -\gamma_1^3 \mathbb{P}^2_1(x)-\gamma_0^4 \mathbb{P}^3_0(x)\Big).
\end{align}
\item Step-3: (Smoothness Indicator) The smoothness indicators $\beta_{-1}^3$, $\beta_0^3$, $\beta_1^3$, and $(\beta_0^5)^{(l)}$ can be obtained using
relations
similar to that shown for WENO-AON(5,3) case. The smoothness indicator corresponding to stencil $\mathbb{S}^4_0$  denoted by $\beta_0^4$, can be given by 
(see \cite{bal-etal_16a})
\begin{equation}\label{SI2}
 \beta_0^4=\Big(u_{41}+\frac{u_{43}}{10}\Big)^2+\frac{13}{3}u_{42}^2+\frac{781}{20}u_{43}^2,
\end{equation}
where 
\begin{equation}
 \left. 
\begin{array}{ll}
u_{41}&=\frac{1}{60}\Big(- 19f_{i-1}-33 f_{i}+63f_{i+1}-11 f_{i+2}\Big),\\
u_{42}&=\frac{1}{2}\Big( f_{i-1}-2f_{i}+ f_{i+1}\Big),\\
u_{43}&=\frac{1}{6}\Big( -f_{i-1}+3f_{i}-3 f_{i+1}+f_{i+2}\Big).
\end{array}
\right \}
\end{equation}

\item Step-4: (Non-linear Weights) We define the parameter $\tau$ as
\begin{equation}
 \tau=\Big(\frac{|(\beta_0^5)^{(l)}-\beta_{-1}^3|+|(\beta_0^5)^{(l)}-\beta_0^3|+|(\beta_0^5)^{(l)}-\beta_1^3|+|(\beta_0^5)^{(l)}-\beta_0^4|}{4}\Big).
\end{equation}
Then, nonlinear weights can be defined using 
% \begin{align*}
%  \omega_{-1}^3&=\gamma_{-1}^3\Big(1+\frac{\tau}{(\epsilon+\beta_{-1}^3)^2}\Big),~~~~\omega_0^3=\gamma_0^3\Big(1+\frac{\tau}{(\epsilon+\beta_0^3)^2}\Big),~~
%  \omega_1^3=\gamma_1^3\Big(1+\frac{\tau}{(\epsilon+\beta_1^3)^2}\Big),\\~~~~\omega_0^4&=\gamma_0^4\Big(1+\frac{\tau}{(\epsilon+\beta_0^4)^2}\Big),~~
%  \omega_0^5=\gamma_0^5\Big(1+\frac{\tau}{(\epsilon+\beta_0^5)^2}\Big).
% \end{align*}
\begin{equation}\label{nlw22}
 \omega_k^l = \frac{\gamma_{k}^l\Big(1+\frac{\tau^2}{(\beta_k^l+\epsilon)^2}\Big)}{\displaystyle{\sum_{j=-1}^1 \gamma_j^3
 \Big(1+\frac{\tau^2}{(\beta_j^3+\epsilon)^2}\Big)+\gamma_0^4\left(1+\frac{\tau^2}{(\beta_0^4+\epsilon)^2}\right)+  \gamma_0^5\Big(1+\frac{\tau^2}{(\beta_0^5+\epsilon)^2}\Big)}},
\end{equation}
where $(l,k)=\{(3,-1), (3,0), (3,1), (4,0), (5,0)\}.$
\item Step-5: The nonlinear reconstruction is defined by
\begin{equation}
 \mathbb{F}_{i+\frac{1}{2}}=\displaystyle{ \omega_0^5 \mathbb{P}( x_{i+\frac{1}{2}})+\omega_0^4 \mathbb{P}^3_0( x_{i+\frac{1}{2}})+
 \sum_{k=-1}^{1} \omega_k^3 \mathbb{P}_{k}^2( x_{i+\frac{1}{2}})   }.
\end{equation}
Similarly, the negative part can be calculated.
\end{itemize}
\begin{remark}{\rm In developing WENO-AO(5,4,3) scheme, we consider the addition of extra stencil $\mathbb{S}_0^4$ for $f_{i+\frac{1}{2}}^{+}$ which help us to improve the resolution of the scheme across the discontinuities in comparison to WENO-AO(5,3) scheme. The resolutions of the scheme can also be improved using 
the left biased stencil $\mathbb{S}_{-1}^4= \{I_{i-2},I_{i-1},I_i,I_{i+1}\}$ instead of taking the 
central stencil. A numerical scheme similar to WENO-AO(5,4,3) can be developed which includes the stencil $\mathbb{S}_{-1}^4$ in place of $\mathbb{S}_0^4$. The cubic polynomial approximation corresponding to the stencil $\mathbb{S}_{-1}^4$ at the interface $x_{i+\frac{1}{2}}$ is
 given by
\[
 \mathbb{P}_{-1}^3\left(x_{i+\frac{1}{2}}\right)= \frac{1}{12}f_{i-2} -\frac{5}{12}f_{i-1} + \frac{13}{12}f_i +\frac{1}{4}f_{i+1}.
\]
The smoothness indicator corresponding  to stencil $\mathbb{S}_{-1}^4$ can be written in compact form (see \cite{bal-etal_16a}) as follow
\begin{equation*}
 \beta_{-1}^4=\left(u_{41}^{'}+\frac{u_{43}^{'}}{10}\right)^2+\frac{13}{3}(u_{42}^{'})^2+\frac{781}{20}(u_{43}^{'})^2.
\end{equation*}
where
\begin{equation*}
 \left. 
\begin{array}{ll}
u_{41}^{'}&=\frac{1}{60}\left(11f_{i-2}-63 f_{i-1}+33f_{i}+19 f_{i+1}\right),\\
u_{42}^{'}&=\frac{1}{2}\left( f_{i-1}-2f_{i}+ f_{i+1}\right),\\
u_{43}^{'}&=\frac{1}{6}\left( -f_{i-2}+3f_{i-1}-3 f_{i}+f_{i+1}\right).
\end{array}
\right \}
\end{equation*}
We named the scheme WENO-AOL(5,4,3), and it produces numerical solutions comparable to 
WENO-AO(5,4,3) scheme as demonstrated in examples (Section \ref{nm}). 
}
\end{remark}

\begin{remark}{\rm
 Here we consider the smoothness indicator $(\beta_0^5)^{(1)}$ for more accurate solutions, but WENO-AO(5,4,3) scheme can be made efficient by using the 
 smoothness indicators
 $(\beta_0^5)^{(2)}$ and $(\beta_0^5)^{(3)}$.}
\end{remark}

\begin{remark}{\rm
 The WENO-AO(5,4,3) scheme is  identical to WENO-AO(5,3) scheme  for $\gamma_{Avg}=0$.}
\end{remark}

\begin{remark}{\rm
 Both of the above schemes can be extended to higher dimensions using dimension by dimension approach.}
\end{remark}

\begin{remark}{\rm
 Both the algorithms can be extended to system of equations using the characteristics-wise approach \cite{shu-97d}. We choose $\lambda$ as the absolute value of 
 largest eigenvalue of the system.}
\end{remark}

%-------------------Analysis of smoothness indicators--------------------
\section{Analysis of new Smoothness Indicator}\label{SI}
WENO schemes of adaptive order  \cite{bal-etal_16a} are shown to be more accurate than classical WENO schemes, however that 
accuracy comes at the cost of computing the extra smoothness indicator for the bigger stencil. In classical WENO-JS scheme \cite{jia-shu_96a}, we need to
 compute smoothness indicators $\beta_k^3$ ($k=-1,0,1$) corresponding to stencils $\mathbb{S}_k^3$, whereas in  case of WENO schemes of adaptive order,
 apart from calculating $\beta_k^3$ for stencils $\mathbb{S}_k^3$, we also need to calculate smoothness indicator $(\beta_0^5)^{(1)}$ for the bigger stencil
 $\mathbb{S}_0^5$ whose computational cost is comparable with sum of computational costs of $\beta_k^3$. To reduce the computational cost, Huang and Chen 
\cite{hua-che_18a} proposed a simple smoothness indicator $(\beta_0^5)^{(2)}$ for the stencil $\mathbb{S}_0^5$ which is obtained from the combination
of lower order smoothness indicators $\beta_k^3$. Further, they have shown that $(\beta_0^5)^{(2)}$ can serve the purpose of smoothness indicator
$(\beta_0^5)^{(1)}$. By using this result computational cost of WENO-AO(5,3) scheme can be reduced  \cite{hua-che_18a}. In lieu of Huang and Chen,
 we also proposed a simple smoothness indicator $(\beta_0^5)^{(3)}$ which is equal to $(\beta_0^5)^{(1)}$ upto truncation error of magnitude
 $O(\Delta x^5)$.
 
In this section, we have compared the proposed smoothness indicator $(\beta_0^5)^{(3)}$ with $(\beta_0^5)^{(1)}$, $(\beta_0^5)^{(2)}$ at the truncation error level of
Taylor series expansion. Using smoothness indicator $(\beta_0^5)^{(3)}$, the numerical scheme WENO-AON(5,3) converges to analytical solution with convergence
rate approximately equal to five for a  sufficient assumption on smoothness of flux $f$.

\subsection{Comparison of smoothness indicator $(\beta_0^5)^{(2)}$ and $(\beta_0^5)^{(3)}$}
Using the Taylor series expansion of $f$ at $x_{i+\frac{1}{2}}$, we can obtain 
\begin{equation}\label{t.si}
 \left. 
\begin{array}{ll}
\beta_{-1}^3&=(f^{'})^2 \Delta x^2 - f^{'} f^{''} \Delta x^3 + \Big(\frac{4}{3}(f^{''})^2-\frac{1}{3}f^{'} f^{''}\Big) \Delta x^4+ O(\Delta x^5),\\
\beta_{0}^3&=(f^{'})^2  \Delta x^2- f^{'} f^{''}  \Delta x^3 + \Big(\frac{4}{3}(f^{''})^2+\frac{2}{3}f^{'} f^{''}\Big) \Delta x^4+ O(\Delta x^5) ,\\
\beta_{1}^3&=(f^{'})^2  \Delta x^2 - f^{'} f^{''} \Delta x^3 + \Big(\frac{4}{3}(f^{''})^2-\frac{1}{3}f^{'} f^{''}\Big) \Delta x^4+ O(\Delta x^5),
\end{array}
\right \}
\end{equation}
whereas for smoothness indicator $(\beta^5_0)^{(1)}$, we have 
\begin{equation}\label{t.rp}
(\beta^5_0)^{(1)}=(f^{'})^2 \Delta x^2 - f^{'} f^{''} \Delta x^3 + \Big(\frac{4}{3}(f^{''})^2+\frac{1}{3}f^{'} f^{''}\Big) \Delta x^4+ O(\Delta x^5).
\end{equation}
The Taylor series expansion of smoothness indicator $(\beta_0^5)^{(3)}$ is given as
\begin{equation}\label{t0.rp}
 (\beta_0^5)^{(3)}=(f^{'})^2 \Delta x^2 - f^{'} f^{''} \Delta x^3 + \Big(\frac{4}{3}(f^{''})^2+\frac{1}{3}f^{'} f^{''}\Big) \Delta x^4+ O(\Delta x^5),
\end{equation}
which is exactly same as $(\beta^5_0)^{(1)}$ upto truncation error level of order $O(\Delta x^5)$, whereas Taylor series expansion of the
smoothness indicator $(\beta_0^5)^{(2)}$ \cite{hua-che_18a}
for $\epsilon =0$  is  
%exact upto  to $(\beta^5_0)^{(1)}$ of order $O(\Delta x^4)$ for $\epsilon =0$ (see \cite{hua-che_18a} for more details).
given by
\begin{equation}\label{t1.rp}
 (\beta_0^5)^{(2)}=(f^{'})^2 \Delta x^2 - f^{'} f^{''} \Delta x^3 + \frac{4}{3}(f^{''})^2 \Delta x^4+ O(\Delta x^5)
\end{equation}
which is exact to $(\beta^5_0)^{(1)}$ of order $O(\Delta x^4)$.
\subsection{Non-linear Weights}
In this subsection, we will show that the nonlinear weights converge to the linear weights when solution is smooth on $\mathbb{S}_0^5$. In the proof,
we have followed the approach discussed in \cite{Arb-etal_17a, ar-etal_11a,kol_14a},  where rigorous analysis is done 
for more general WENO schemes of adaptive order with both WENO-JS and WENO-Z weights.
Here we adapt the notation used in \cite{Arb-etal_17a}. For a function $g$, we have
\begin{equation}
 g={O}(\Delta x^r)~~~~\Leftrightarrow~~~~~ |\Delta x^{-r} g|\leq C ~~~\mbox{as}~~ \Delta x\rightarrow 0
\end{equation}
where $C$ is positive constants independent of $\Delta x$. 
\begin{theorem}{\rm
 If flux $f$ is smooth on the stencil $\mathbb{S}_0^5$, then following results hold in case of WENO-AON(5,3) reconstruction for nonlinear weights
 \begin{equation}
 \omega_k^l = \gamma_{k}^l +
                O(\Delta x^4) 
\end{equation}
where $\beta_0^5$ can be chosen from any one of the  smoothness indicators $(\beta_0^5)^{(1)}$, $(\beta_0^5)^{(2)}$, $(\beta_0^5)^{(3)}$.
}
\end{theorem}
{\it Proof}. The nonlinear weights in case of WENO-AON(5,3) reconstruction are given by
\begin{equation}\label{nlw}
 \omega_k^l = \frac{\gamma_{k}^l\Big(1+\frac{\tau^2}{(\beta_k^l+\epsilon)^2}\Big)}{\displaystyle{\sum_{j=-1}^1 \gamma_j^3
 \Big(1+\frac{\tau^2}{(\beta_j^3+\epsilon)^2}\Big)+  \gamma_0^5\Big(1+\frac{\tau^2}{(\beta_0^5+\epsilon)^2}\Big)}}~~\mbox{for}~~(l,k)=\{(3,-1),
 (3,0), (3,1), (5,0)\}.
\end{equation}
Before proving the results for the nonlinear weights, we derived some results for the smoothness indicators and parameter $\tau$. 
From \eqref{t.si}, \eqref{t.rp}, \eqref{t0.rp}  and \eqref{t1.rp}, we can  easily conclude that 
\begin{equation}
 |\beta_0^5-\beta_k^3|=O(\Delta x^4), ~~~k=-1,0,1.
\end{equation}
Using the above, we can get $\tau =O(\Delta x^4)$. Further, the estimates of smoothness indicators is given as
 \begin{equation}
{(\beta_k^l+\epsilon)^2}=
                O(\Delta x^4),
\end{equation}
and 
 \begin{equation}\label{tau.es}
\frac{\tau^2}{(\beta_k^l+\epsilon)^2}=
                O(\Delta x^4),
\end{equation}
Since $\gamma_0^5+\sum_{k=-1}^1 \gamma_k^3 =1$,  equation \eqref{nlw} can be written as 
\begin{equation}\label{nlw1}
 \omega_k^l = \frac{\gamma_{k}^l\Big(1+\eta_k^l\Big)}{1+\displaystyle{\sum_{j=-1}^1 \gamma_j^3
 \eta_k^3+  \gamma_0^5\eta_0^5}}~~\mbox{for}~~(l,k)=\{(3,-1),
 (3,0), (3,1), (5,0)\},
\end{equation}
where $\eta_k^l=\frac{\tau^2}{(\beta_k^l+\epsilon)^2}$.
Since $\eta_k^l\rightarrow 0$ as $\Delta x\rightarrow 0$, we have
\begin{equation}
({1+\displaystyle{\sum_{j=-1}^1 \gamma_j^3
 \eta_k^3+  \gamma_0^5\eta_0^5}}) = 1+
                O(\Delta x^4).
\end{equation}
Further, we have 
\begin{align}
\frac{1}{ ({1+\displaystyle{\sum_{j=-1}^1 \gamma_j^3
 \eta_k^3+  \gamma_0^5\eta_0^5}})}
 &=  \frac{1}{1+O(\Delta x^4)}
 = 1- O(\Delta x^4) 
\end{align}
Finally, using \eqref{tau.es} we get
 \begin{equation*}
 \omega_k^l = \gamma_{k}^l +
                O(\Delta x^4).
\end{equation*}
\begin{remark}
 {\rm  The parameter $\epsilon$ is used to prevent division by zero and must be chosen to be a small value, so that it does not dominate the smoothness indicators $\beta$.
 As mentioned before we use $\epsilon = 10^{-12}$ but another way to choose this is $\epsilon = k \Delta x^2$ for small $k$.

 }
\end{remark}

\begin{theorem}\label{thm.n53}{\rm
 If the flux $f$ is smooth on the stencil $\mathbb{S}_0^5$, then following results hold in case of WENO-AO(5,4,3) reconstruction for the 
 nonlinear weights
 \begin{equation}
 \omega_k^l = \gamma_{k}^l +
                O(\Delta x^4).
\end{equation}
}
\end{theorem}
{\it Proof}. Similar to the case of WENO-AON(5,3) reconstruction.
\subsection{Accuracy of reconstruction schemes}
Now we prove the accuracy of WENO-AON(5,3) and WENO-AO(5,4,3) reconstruction, when flux is smooth in the bigger stencil $\mathbb{S}_0^5$, 
following the
 approach discussed in \cite{Arb-etal_17a} (also see \cite{ar-etal_11a}, \cite{kol_14a} for more details).
\begin{theorem}\label{thm.543}{\rm
If flux function is smooth on $\mathbb{S}_0^5$, then WENO-AON(5,3) and WENO-AO(5,4,3) reconstruction are fifth order accurate at the interface
\begin{equation}
 \mathbb{F}_{i+\frac{1}{2}}- \mathbb{H}(x_{i+\frac{1}{2}})= O(\Delta x^5).
\end{equation}

}
\end{theorem}
{\it Proof}. In the case of WENO-AON(5,3) reconstruction, we have
 \begin{align*}
 \mathbb{F}_{j+\frac{1}{2}}- \mathbb{H}(x_{i+\frac{1}{2}}) & = \displaystyle{ \frac{\omega^5_0}{\gamma^5_0}\Big( \mathbb{P}^4_0(x_{i+\frac{1}{2}})-
 \sum_{k=-1}^1 \gamma_k^3
 \mathbb{P}_k^2(x_{i+\frac{1}{2}})\Big)   +\sum_{k=-1}^1 \omega_k^3\mathbb{P}_k^2(x_{i+\frac{1}{2}})-\mathbb{H}(x_{i+\frac{1}{2}}}) \\      
 &=  \frac{\omega^5_0}{\gamma^5_0}\Big(\mathbb{P}^4_0(x_{i+\frac{1}{2}}) - \mathbb{H}(x_{i+\frac{1}{2}}) +\mathbb{H}(x_{i+\frac{1}{2}})-
 \sum_{k=-1}^1\gamma_k^3 \mathbb{P}^2_k(x_{i+\frac{1}{2}})\Big)\\
 &+\sum_{k=-1}^1 \omega_k^3\mathbb{P}_k^2(x_{i+\frac{1}{2}})-(\omega_0^5+ \sum_{k=-1}^1 \omega_k^3)\mathbb{H}(x_{i+\frac{1}{2}})  \\
 &=  \frac{\omega^5_0}{\gamma^5_0}\Big(\mathbb{P}^4_0(x_{i+\frac{1}{2}}) - \mathbb{H}(x_{i+\frac{1}{2}}) + (\gamma_0^5+ \sum_{k=-1}^1 \gamma_k^3)\mathbb{H}(x_{i+\frac{1}{2}})\\
 &-\sum_{k=-1}^1\gamma_k^3 \mathbb{P}^2_k(x_{i+\frac{1}{2}})\Big)
 +\sum_{k=-1}^1 \omega_k^3\mathbb{P}_k^2(x_{i+\frac{1}{2}})
-(\omega_0^5+ \sum_{k=-1}^1 \omega_k^3)\mathbb{H}(x_{i+\frac{1}{2}})
\end{align*}
 \begin{align*}
 \mathbb{F}_{i+\frac{1}{2}}- \mathbb{H}(x_{i+\frac{1}{2}}) 
 &=  \frac{\omega^5_0}{\gamma^5_0}\Big((\mathbb{P}^4_0(x_{i+\frac{1}{2}}) - \mathbb{H}(x_{i+\frac{1}{2}})) 
 -\sum_{k=-1}^1\Big( \frac{\omega_0^5-\gamma_0^5}{\gamma_0^5} \gamma_k^3-(\omega_k^3-\gamma_k^3)\Big)(\mathbb{P}_k^2(x_{i+\frac{1}{2}})\\
 &-\mathbb{H}(x_{i+\frac{1}{2}}))
\end{align*}
The sufficient condition to achieve fifth order accurate reconstruction using   WENO-AON(5,3) reconstruction is
\begin{equation}
 \sum_{k=-1}^1\Big( \frac{\omega_0^5-\gamma_0^5}{\gamma_0^5} \gamma_k^3-(\omega_k^3-\gamma_k^3)\Big)= O(\Delta x^2),
\end{equation}
which can be  observed from Theorem \ref{thm.n53}.
In the case of WENO-AO(5,4,3) reconstruction, we have
\begin{align}\label{rec:543}
 \mathbb{F}_{i+\frac{1}{2}}-\mathbb{H}(x_{i+\frac{1}{2}}) & = \frac{\omega_0^5}{\gamma_0^5}(\mathbb{P}_0^4(x_{i+\frac{1}{2}})-\mathbb{H}(x_{i+\frac{1}{2}}))-\Big(\frac{\omega_0^5- 
 \gamma_0^5}{\gamma_0^5} \gamma_0^4- (\omega_0^4-\gamma_0^4)\Big) (\mathbb{P}_0^3(x_{i+\frac{1}{2}})-\mathbb{H}(x_{i+\frac{1}{2}}))\nonumber\\
 &- \sum_{k=-1}^1 \Big(\frac{\omega_0^5-\gamma_0^5}{\gamma_0^5} \gamma_k^3 -(\omega_k^3 -\gamma_k^3)\Big)(\mathbb{P}_k^2(x_{i+\frac{1}{2}})-\mathbb{H}(x_{i+\frac{1}{2}})) 
\end{align}
Using similar arguments as in the case of WENO-AON(5,3) reconstruction, we can also conclude that WENO-AO(5,4,3) reconstruction is fifth order accurate.

\begin{theorem}\label{th:lww}{\rm
If the solution contains a discontinuity in $\mathbb{S}_0^5$ and is smooth in $\mathbb{S}_0^4$, then nonlinear weights in WENO-AO(5,4,3) reconstruction  satisfy 
\begin{align*}
 \omega_0^5 &=O(\Delta x^4), ~~~\omega_{-1}^3 =O(\Delta x^4),\\
 \omega_k^l &= \frac{\gamma_k^l}{\gamma_0^4+\gamma_0^3 +\gamma_1^3}+O(\Delta x), ~~~(l,k) =\{(4,0), (3,0), (3,1)\}.
\end{align*}
}
\end{theorem}
{\it Proof}. The nonlinear weights in case of WENO-AO(5,4,3) reconstruction given by \eqref{nlw22} can be written as
\begin{align}\label{eq:wi2}
 \omega_k^l 
 &=\frac{\gamma_{k}^l\Big(1+\frac{\tau^2}{(\beta_k^l+\epsilon)^2}\Big)}{1+\left(\displaystyle{\sum_{j=-1}^1 \gamma_j^3
 \frac{\tau^2}{(\beta_j^3+\epsilon)^2}+\gamma_0^4\frac{\tau^2}{(\beta_0^4+\epsilon)^2}+  \gamma_0^5\frac{\tau^2}{(\beta_0^5+\epsilon)^2}}\right)}.
\end{align}
We have, for $\epsilon \leq \Delta x^2$
 \begin{align*}
  \frac{\tau^2}{(\beta_0^4+\epsilon)^2} &= \frac{\tau^2}{((f^{'})^2\Delta x^2- f^{'} f^{'''}\Delta x^3 +O(\Delta x^4))^2}\\
          &= \frac{\tau^2}{(f^{'})^4 \Delta x^4}\left(1- \frac{f^{''}}{f^{'}}\Delta x + O(\Delta x^2)\right)^{-2}=\frac{\tau^2}{(f^{'})^4 \Delta x^4}+O(\Delta x^{-3}).
 \end{align*}
 Similarly
 \begin{equation}
  \frac{\tau^2}{(\beta_0^3+\epsilon)^2} = \frac{\tau^2}{(f^{'})^4 \Delta x^4}+O(\Delta x^{-3}),~~\frac{\tau^2}{(\beta_1^3+\epsilon)^2} = \frac{\tau^2}{(f^{'})^4 \Delta x^4}+O(\Delta x^{-3}),~~
 \end{equation}
 and 
 \begin{equation}\label{eq:wi0}
  \frac{\tau^2}{(\beta_0^5+\epsilon)^2} = \frac{\tau^2}{O(1)},~~~~~~\frac{\tau^2}{(\beta_{-1}^3+\epsilon)^2} = \frac{\tau^2}{O(1)}.
 \end{equation}
 Further, the denominator of \eqref{eq:wi2} can be estimated as
\begin{align}\label{eq:wi1}
 1+&\left(\displaystyle{\sum_{j=-1}^1 \gamma_j^3
 \frac{\tau^2}{(\beta_j^3+\epsilon)^2}+\gamma_0^4\frac{\tau^2}{(\beta_0^4+\epsilon)^2}+  \gamma_0^5\frac{\tau^2}{(\beta_0^5+\epsilon)^2}}\right)\nonumber\\ &=1+ \gamma_{-1}^3 \frac{\tau^2}{O(1)}+\gamma_0^3\frac{\tau^2}
 {(f^{'})^4 \Delta x^4}
 +\gamma_{1}^3 \frac{\tau^2}{(f^{'})^4 \Delta x^4} +\gamma_{0}^4 \frac{\tau^2}{(f^{'})^4 \Delta x^4}+\gamma_{0}^5 \frac{\tau^2}{O(1)}+O(\Delta x^{-3})\nonumber\\
 &=\frac{\tau^2 (\gamma_0^4+\gamma_0^3+\gamma_1^3)}
 {(f^{'})^4 \Delta x^4}{(1+O(\Delta x))}.
\end{align}
\begin{itemize}
 \item Case-I: When $(l,k)=\{(5,0),(3,-1)\}$, then using \eqref{eq:wi0} and \eqref{eq:wi1}, we have
 \begin{equation}
  \omega_k^l =O(\Delta x^4).
 \end{equation}

 \item Case-II: For $(l,k)=(4,0)$, we have
\begin{align*}
 \omega_0^4 &= \frac{\gamma_0^4 \left(1+\frac{\tau^2}{(f^{'})^4 \Delta x^4}+O(\Delta x^{-3})\right)(f^{'})^4 \Delta x^4}{\tau^2 (\gamma_0^4+\gamma_0^3+\gamma_1^3)(1 +O(\Delta x))}\\
            &= \left(\frac{\gamma_0^4}{\gamma_0^4+\gamma_0^3+\gamma_1^3}\right) (1+O(\Delta x))(1-O(\Delta x))\\
            & = \frac{\gamma_0^4}{\gamma_0^4+\gamma_0^3+\gamma_1^3}+O(\Delta x).
\end{align*}
Similarly, we can prove for $\omega_0^3$ and $\omega_1^3$.
\end{itemize}

\begin{theorem}\label{thm:543w}{\rm
If the solution contains a discontinuity in $\mathbb{S}_0^5$ and is smooth in $\mathbb{S}_0^4$, then WENO-AO(5,4,3) reconstruction is fourth order accurate at the interface 
\begin{equation}
 \mathbb{F}_{i+\frac{1}{2}}- \mathbb{H}(x_{i+\frac{1}{2}})= O(\Delta x^4),
\end{equation}
provided $\gamma_0^3=\gamma_1^3$.
}
\end{theorem}
{\it Proof}. The WENO-AO(5,4,3) flux reconstruction at the interface is given by \eqref{rec:543}. Under the assumption that flux is smooth over $\mathbb{S}_0^4$ but contains a discontinuity in $\mathbb{S}_0^5$, we have
% \begin{align*}
%  \mathbb{F}_{i+\frac{1}{2}}-\mathbb{H}(x_{i+\frac{1}{2}}) & = \frac{\omega_0^5}{\gamma_0^5}(\mathbb{P}_0^4(x_{i+\frac{1}{2}})-
%  \mathbb{H}(x_{i+\frac{1}{2}}))\\&+ \Big(\frac{\omega_0^5- 
%  \gamma_0^5}{\gamma_0^5} \gamma_0^4- (\omega_0^4-\gamma_0^4)\Big) (\mathbb{P}_0^3(x_{i+\frac{1}{2}})-\mathbb{H}(x_{i+\frac{1}{2}}))\\
%  &- \sum_{k=-1}^1 \Big(\frac{\omega_0^5-\gamma_0^5}{\gamma_0^5} \gamma_k^3 -(\omega_k^3 -\gamma_k^3)\Big)(\mathbb{P}_k^2(x_{i+\frac{1}{2}})-\mathbb{H}(x_{i+\frac{1}{2}})) 
% \end{align*}
\begin{equation*}
 \omega_0^5 = O(\Delta x^4) ~~\mbox{and}~~\omega_{-1}^3 =O(\Delta x^4).
\end{equation*}
Then the flux reconstruction \eqref{rec:543} reduce to
\begin{align*}
 \mathbb{F}_{i+\frac{1}{2}}-\mathbb{H}(x_{i+\frac{1}{2}})  = & - \omega_0^4 (\mathbb{P}_0^3(x_{i+\frac{1}{2}})-\mathbb{H}(x_{i+\frac{1}{2}}))
 +  \omega_0^3 (\mathbb{P}_0^2(x_{i+\frac{1}{2}})-\mathbb{H}(x_{i+\frac{1}{2}})) \\
 &+ \omega_1^3 (\mathbb{P}_1^2(x_{i+\frac{1}{2}})-\mathbb{H}(x_{i+\frac{1}{2}})).
\end{align*}
From Theorem \ref{th:lww}, we have
\begin{equation}
 \omega_0^3= \frac{\gamma_0^3}{\gamma_0^4+\gamma_0^3+\gamma_1^3}+O(\Delta x),~~~\omega_1^3= \frac{\gamma_1^3}{\gamma_0^4+\gamma_0^3+\gamma_1^3}+O(\Delta x).
\end{equation}
Now if we take $\gamma_0^3=\gamma_1^3 =\gamma_0$, then we have
 \begin{align*}
 \mathbb{F}_{i+\frac{1}{2}}-\mathbb{H}(x_{i+\frac{1}{2}})  = & - \omega_0^4 (\mathbb{P}_0^3(x_{i+\frac{1}{2}})-\mathbb{H}(x_{i+\frac{1}{2}}))
 +  (\frac{\gamma_0}{\gamma_0^4+2\gamma_0} +O(\Delta x))  (\mathbb{P}_0^2(x_{i+\frac{1}{2}})-\mathbb{H}(x_{i+\frac{1}{2}})) \\
 &+ (\frac{\gamma_0}{\gamma_0^4+2\gamma_0} +O(\Delta x))(\mathbb{P}_1^2(x_{i+\frac{1}{2}})-\mathbb{H}(x_{i+\frac{1}{2}})),
\end{align*}
which can be further written as 
 \begin{align}\label{rec:5431}
 \mathbb{F}_{i+\frac{1}{2}}-\mathbb{H}(x_{i+\frac{1}{2}})  = & - \omega_0^4 (\mathbb{P}_0^3(x_{i+\frac{1}{2}})-\mathbb{H}(x_{i+\frac{1}{2}}))\nonumber\\
 &+  (\frac{\gamma_0}{\gamma_0^4+2\gamma_0} +O(\Delta x)) (\mathbb{P}_0^2(x_{i+\frac{1}{2}})+\mathbb{P}_1^2(x_{i+\frac{1}{2}})-2\mathbb{H}(x_{i+\frac{1}{2}})).
\end{align}
We can easily observe that 
\begin{equation}\label{eq:sum}
 \frac{1}{2}(\mathbb{P}_0^2(x_{i+\frac{1}{2}})+ \mathbb{P}_1^2(x_{i+\frac{1}{2}}))= \mathbb{P}_0^3(x_{i+\frac{1}{2}}),
\end{equation}
where
\begin{equation*}
 \left. 
\begin{array}{ll}
 \mathbb{P}_0^2(x_{i+\frac{1}{2}})&=-\frac{1}{6}f_{i-1}+ \frac{5}{6}f_i +\frac{1}{3}f_{i+1},\\
 \mathbb{P}_1^2(x_{i+\frac{1}{2}})&= \frac{1}{3}f_{i}+\frac{5}{6} f_{i+1}-\frac{1}{6}f_{i+2},\\
 \mathbb{P}_0^3(x_{i+\frac{1}{2}})&=-\frac{1}{12}f_{i-1}+\frac{7}{12}f_{i}+\frac{7}{12}f_{i+1}-\frac{1}{12}f_{i+2}.
\end{array}
\right \}
\end{equation*}
% \begin{align*}
%  \mathbb{P}_0^2(x_{i+\frac{1}{2}})&=-\frac{1}{6}f_{i-1}+ \frac{5}{6}f_i +\frac{1}{3}f_{i+1}\\
%  \mathbb{P}_1^2(x_{i+\frac{1}{2}})&= \frac{1}{3}f_{i}+\frac{5}{6} f_{i+1}-\frac{1}{6}f_{i+2}\\
%  \mathbb{P}_0^3(x_{i+\frac{1}{2}})&=-\frac{1}{12}f_{i-1}+\frac{7}{12}f_{i}+\frac{7}{12}f_{i+1}-\frac{1}{12}f_{i+2}
% \end{align*}
Using \eqref{eq:sum} in \eqref{rec:5431}, we get the required result.

\begin{remark}{\rm 
As per the above theorem if flux is smooth over stencil $\mathbb{S}_0^4$ but not over $\mathbb{S}_0^5$, then WENO-AO(5,4,3) reconstruction achieves the fourth order at the interface $x_{i+\frac{1}{2}}$ 
under the assumption $\gamma_0^3=\gamma_1^3$. But our numerical experiments suggest that giving more weight to stencil $\mathbb{S}_0^3$ in comparison to stencils $\mathbb{S}_{-1}^3$, $\mathbb{S}_1^3$ (similar to WENO-AO(5,3) scheme 
\cite{bal-etal_16a}), yields better resolution 
of solution across shocks or discontinuities. This is illustrated in examples (Section \ref{nm}). This is the reason of choosing linear weights like \eqref{lw:543}, where we have given more weightage to stencil $\mathbb{S}_0^3$ than $\mathbb{S}_{-1}^3$, $\mathbb{S}_1^3$ .
} 
\end{remark}

%---------------------SSPRK3------------------------------------------------------------------------
% \section{Time-Integration}\label{ssprk3}
% For the time integration, we consider the explicit Runge-Kutta method of order three \cite{shu-osh_88a}. The Runge-Kutta method 
%  considered here, can be defined
%   for the following ordinary differential equation 
%   \begin{equation}
%    \frac{du}{dt}=\mbox{RHS}(u)
%   \end{equation}
% for vector $u^n$ as follows 
%  \begin{align*}
%   u_{0}&=u^n\\
%   u_{1}& =u_0+\Delta t \mbox{RHS}(u_0),\\
%   u_{2}& =\frac{3}{4}u_0+\frac{1}{4}(u_1+\Delta t \mbox{RHS}(u_{1})),\\
%   u_{3}& =\frac{1}{3}u_0+\frac{2}{3}(u_2+\Delta t \mbox{RHS}(u_{2})),\\
%   u^{n+1}&=u_{3}.
%  \end{align*}
%  When the solution $u^n$ is known, then it can be updated to the next level using the above stages.
%  
%--------------------Numerical experiments----------------------------------------------------------

\begin{table}
\begin{center}
\small
\begin{tabular}{|l| c|c| c| c| c|  c| c| r|}
\hline 
 $N$ & WENO-AO(5,3) & Order & WENO-AON(5,3) & Order & WENO-AO(5,4,3) & Order \\
\hline
20  &  1.7343e-03 & --     &1.7462e-03 &-- &1.734265e-03 &\\ 
\hline 
40  &  5.6930e-05 & 4.93     &5.6971e-05 &4.94  &5.693340e-05 & 4.93    \\
\hline
80  &  1.8762e-06 & 4.92     &1.8763e-06 &4.92  &1.876227e-06 & 4.92   \\
\hline
160  &  6.2731e-08 & 4.90     &6.2731e-08 &4.90  &6.273129e-08 & 4.90    \\
\hline
320  &  2.1399e-09 & 4.87     &2.1399e-09 &4.87  &2.139861e-09 & 4.87    \\
\hline
640  &  9.4846e-11 & 4.50     &9.4848e-11 &4.50  &9.484731e-11 & 4.50   \\
\hline
\end{tabular}
\caption{ Comparison of $L^{\infty}$-error for WENO-AO(5,3), WENO-AON(5,3), and
WENO-AO(5,4,3)  schemes  along with their convergence rates  for Example \ref{example.1} at time $T=10$.
}
\label{Table.advl8}
\end{center}
\end{table}

\begin{table}
\centering
\small
\begin{tabular}{|l| c|c| c| c| c|  c| c| r|}
\hline 
 $N$ & WENO-AO(5,3) & Order   & WENO-AON(5,3) & Order & WENO-AO(5,4,3) & Order \\
\hline
20  &  2.2065e-03 & --       &2.2064e-03    &-- &2.2065e-03 &\\ 
\hline 
40  &  7.2469e-05 & 4.93     &7.2469e-05    &4.93  &7.2468e-05 & 4.93    \\
\hline
80  &  2.3888e-06 & 4.92     &2.3888e-06    &4.92  &2.3888e-06 & 4.92   \\
\hline
160  &  7.9873e-08 & 4.90    &7.9873e-08    &4.90  &7.9873e-08 & 4.90   \\
\hline
320  &  2.7247e-09 & 4.87    &2.7247e-09    &4.87  &2.7247e-09 & 4.87    \\
\hline
640  &  1.2075e-10 & 4.50    &1.2075e-10    &4.50  &1.2075e-10 & 4.50    \\
\hline
\end{tabular}
\caption{   Comparison of $L^{1}$-error for WENO-AO(5,3), WENO-AON(5,3), and
WENO-AO(5,4,3)  schemes along with their convergence rates  for Example \ref{example.1} at time $T=10$.}
\label{Table.advl1}
\end{table}
\begin{figure}[t!]
 \includegraphics[width=0.48\textwidth]{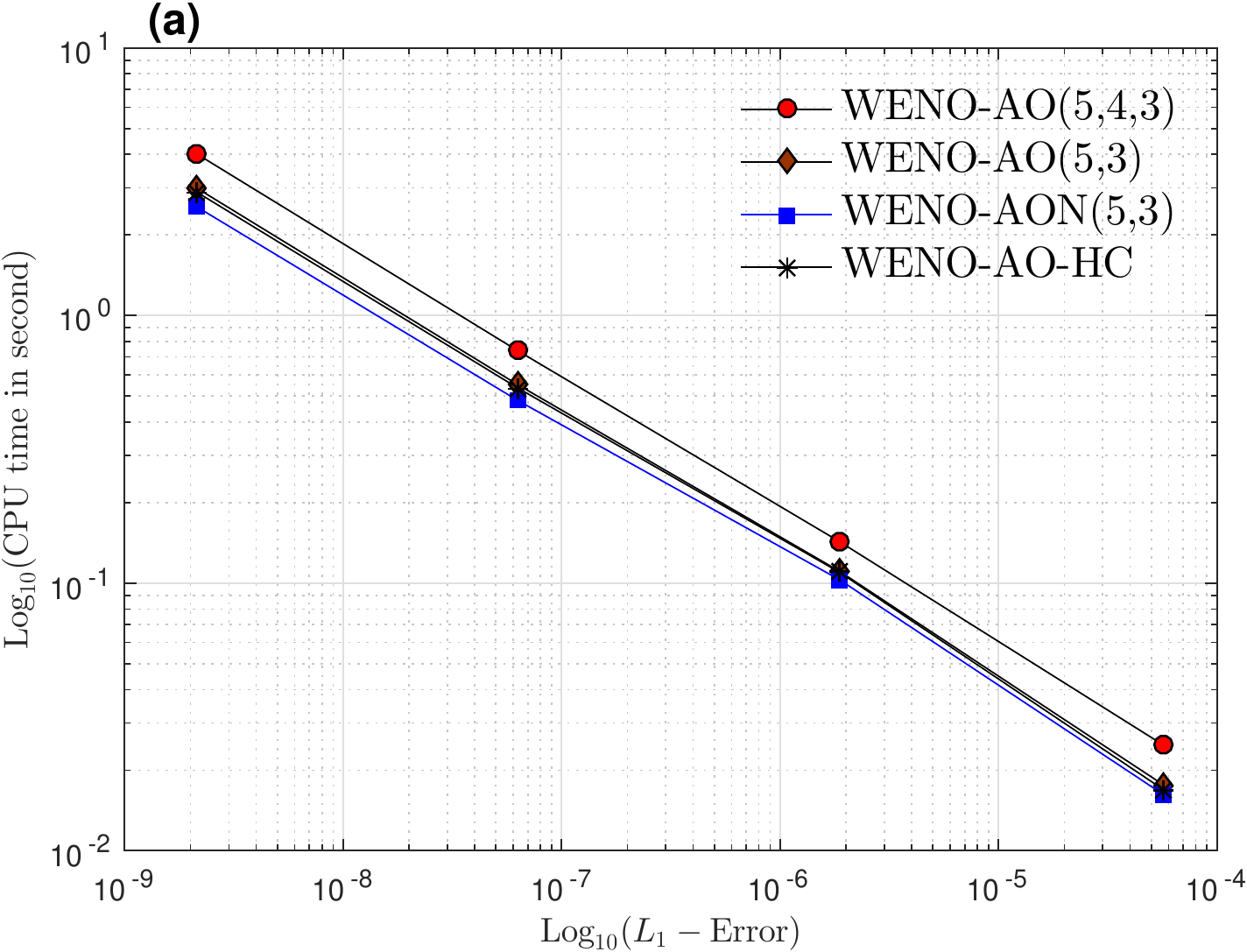}
 \includegraphics[width=0.48\textwidth]{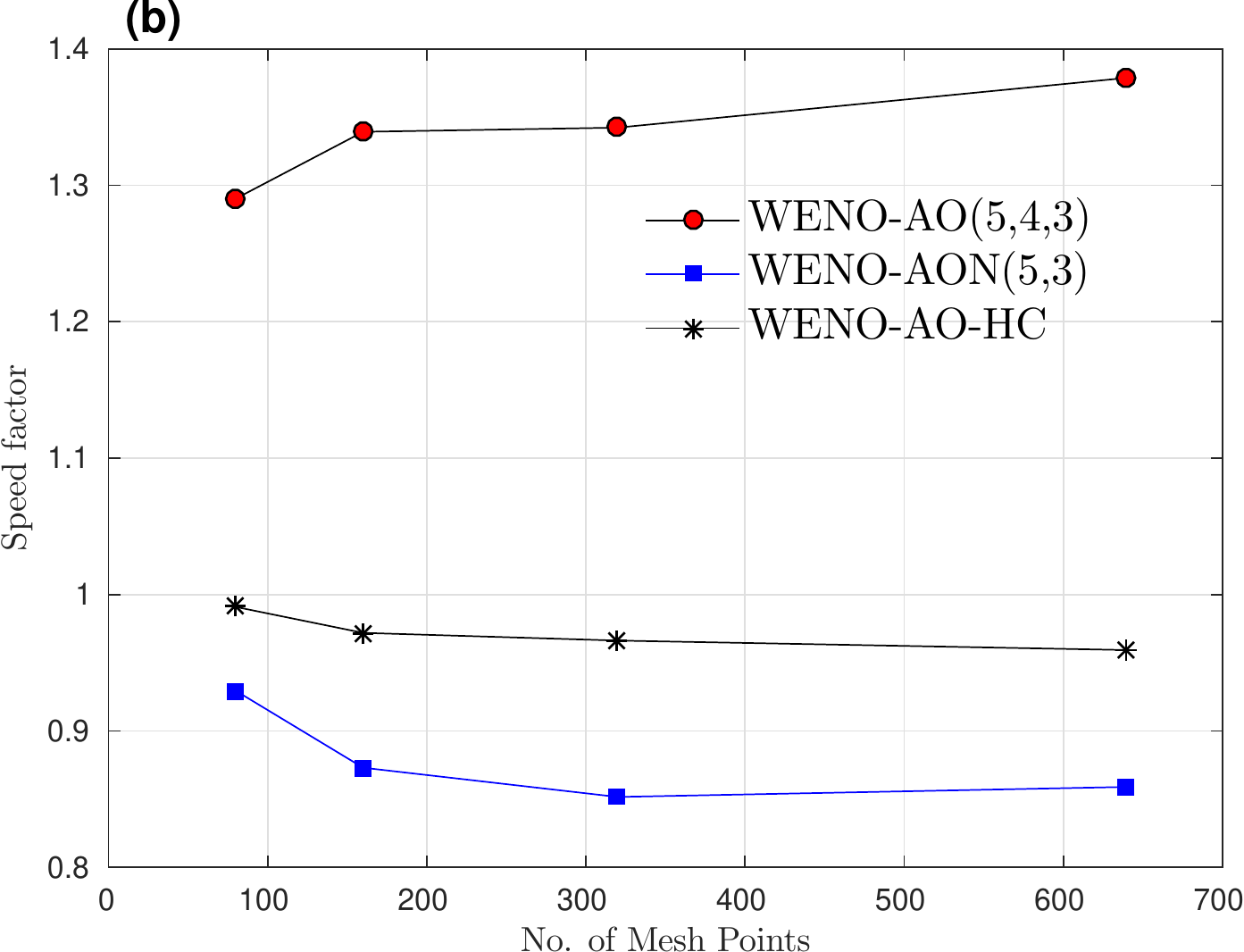}
 \caption{(a) Comparison of WENO schemes for Example \ref{example.1} in computing  $L^1-$error and CPU time, (b) relative CPU time taken by
 WENO-AO-HC, WENO-AON(5,3), and
WENO-AO(5,4,3) schemes with respect to WENO-AO(5,3) for different number of mesh points.}
 \label{speed.adv}
\end{figure}

\section{Numerical Experiments}\label{nm}
For the time integration, we consider the explicit Runge-Kutta method of order three \cite{shu-osh_88a}. The Runge-Kutta method 
 considered here, can be defined
  for the following ordinary differential equation 
  \begin{equation}
   \frac{du}{dt}=\mbox{RHS}(u)
  \end{equation}
for vector $u^n$ as follows 
 \begin{align*}
  u_{0}&=u^n\\
  u_{1}& =u_0+\Delta t \mbox{RHS}(u_0),\\
  u_{2}& =\frac{3}{4}u_0+\frac{1}{4}(u_1+\Delta t \mbox{RHS}(u_{1})),\\
  u_{3}& =\frac{1}{3}u_0+\frac{2}{3}(u_2+\Delta t \mbox{RHS}(u_{2})),\\
  u^{n+1}&=u_{3}.
 \end{align*}
 When the solution $u^n$ is known, then it can be updated to the next level using the above stages.
 
In this section, our aim is to compare the two proposed algorithms WENO-AON(5,3) and WENO-AO(5,4,3) with other variants of WENO schemes of order five.
Here, we have used 
WENO5-JS \cite{jia-shu_96a}, WENO-Z \cite{bor-etal_08a}, WENO-ZQ \cite{zhu-qiu_16a}, WENO-AO(5,3) \cite{bal-etal_16a}, and WENO-AO-HC \cite{hua-che_18a} 
schemes for the comparison purpose. In order to have a fair comparison, in the WENO reconstruction procedure, basis space of reconstruction polynomials are 
spanned by Legendre polynomials. This allows us to express smoothness indicators in the compact form,
thereby reducing the computational costs of old variants of WENO schemes  \cite{bal-etal_09a}, \cite{bal-etal_16a}. All WENO
schemes considered here have same framework as given in the
Section \ref{WENOAO} and they differ from each other at the 
reconstruction level only. We have compared the proposed algorithms with other WENO algorithms in terms of accuracy and computational cost. 
Further, numerical experiments are  performed to show how well the WENO schemes resolve the various discontinuities for the scalar and system of 
nonlinear problems.  The CFL number is taken to be 0.95 in the 1D test cases and 0.5 for 2D test cases or it will be mentioned separately.

 The relative speedup factors are calculated for each WENO schemes with respect to WENO-AO(5,3) scheme.  The speedup factor of  the 
 WENO-AO(5,3) algorithm  gets normalized to unity. If the speedup number of the scheme is less than unity, it indicates a faster 
 scheme  in
 comparison to that of WENO-AO(5,3) scheme and vice versa.

The values of parameters $\gamma_{Hi}$ and $\gamma_{Lo}$ are taken as 0.85  \cite{bal-etal_16a} and \cite{hua-che_18a} for all WENO schemes of adaptive order. 
In case of WENO-ZQ \cite{zhu-qiu_16a} scheme, we have taken the weight of 0.98 for the large central stencil. In WENO-AO(5,4,3) scheme, we  take
$\gamma_{Hi}=0.85,~\gamma_{Avg}=0.85,~\mbox{and}~\gamma_{L0}=0.7$ for all test-cases.
\subsection{Accuracy Test}
We test the accuracy and convergence rate  of WENO schemes for linear as well as a nonlinear equations.
The accuracy of the schemes are measured in $L^{\infty}$-, and  $L^1$-error norms, which are defined for error $e$ over the domain $[a,b]$ as follows
 \begin{align*}
 \|e\|_{\infty}&=\displaystyle{\max_{j} |u_j-(u_{\Delta x})_j|},\\ \|e\|_1&=\displaystyle{\frac{(b-a)}{N+1} \sum_{j}|u_j-(u_{\Delta x})_j|},
\end{align*}
where $N$ denotes the number of subdivisions of the domain and $u_j$ and $(u_{\Delta x})_j$ denote the exact and approximate solutions 
(corresponding to $\Delta x=(b-a)/N$) at point $x_j$, respectively.

Comparison of the accuracy of the  WENO-AO-HC, WENO-AO(5,3), WENO-JS and WENO-ZQ schemes can be found in \cite{bal-etal_16a} \cite{hua-che_18a}.
We have devised WENO-AON(5,3) scheme using an
approach similar to \cite{hua-che_18a}.  Therefore, we intend to compare the computational cost of WENO-AON(5,3) with that of WENO-AO-HC. In order to save space, we have compared only the computational
cost of the WENO-AO(5,3), WENO-AO-HC,
WENO-AON(5,3), and WENO-AO(5,4,3) schemes and accuracy of  WENO-AO(5,3), 
WENO-AON(5,3), and WENO-AO(5,4,3) schemes.
\begin{example}\label{example.1}{\rm (Linear Advection)
 Consider the linear advection equation given by
 \begin{equation}\label{lin.adv}
  u_t+u_x=0,\;\;\;\;(x,t)\in (-1,1)\times (0,T]
 \end{equation}
with the initial data $u(x,0)=\sin(\pi x)$ and periodic boundary conditions. Numerical solutions are computed using WENO-AO(5,3), WENO-AON(5,3), and
WENO-AO(5,4,3) schemes at time $T=10$. To ensure  negligible effect of time-discretization, we take time step 
$\Delta t= 0.5\Delta x^{1.5}$.
The $L^{\infty}$ and $L^1$-errors and  convergence rates of the considered
 schemes are shown in Tables \ref{Table.advl8} and \ref{Table.advl1}, respectively. The magnitude of the errors and their
 corresponding convergence rates are comparable
 for each of the WENO schemes of adaptive order. As the grid is refined, WENO-AO(5,3), WENO-AON(5,3), and
WENO-AO(5,4,3) schemes converge to solution with rate almost equal to five.
 Further, we observe that, on the same mesh, 
comparing with the WENO-AO(5,3) scheme, the WENO-AON(5,3) and WENO-AO-HC scheme uses  less computational time, whereas WENO-AO(5,4,3) uses more time than
WENO-AO(5,3) scheme. We can see from Figure \ref{speed.adv}, WENO-AO(5,4,3) takes  approximately $34\%$ percent more time, whereas WENO-AO-HC and WENO-AON(5,3) takes
 approximately $4\%$ percent  and $12\%$ percent less 
 computational time, respectively, than  WENO-AO(5,3) scheme.

}
\end{example}
\begin{table}[b]
\centering 
\small
\begin{tabular}{|l| c|c| c| c| c|  c| c| r|}
\hline 
$N$ &WENO-AO(5,3) & Order & WENO-AON(5,3) & Order & WENO-AO(5,4,3) & Order \\
\hline
20  &  5.0311e-03 & --     &5.0211e-03 &-- &5.0162e-03 &\\ 
\hline 
40  &  3.7929e-04 & 3.73     &3.7923e-04 &3.73  &3.7921e-04 & 3.73    \\
\hline
80  &  1.5331e-05 & 4.69    &1.5331e-05 &4.69  &1.5330e-05 & 4.69    \\
\hline
160  &  4.6722e-07 & 5.03    &4.6722e-07 &5.04  &4.6722e-07 & 5.04   \\
\hline
320  &  1.3543e-08 & 5.11     &1.3543e-08 &5.11  &1.3543e-08 & 5.11    \\
\hline
\end{tabular}
\caption{  Comparison of $L^{\infty}$-errors for WENO-AO(5,3), WENO-AON(5,3), and
WENO-AO(5,4,3) schemes  along with their convergence rates  for Example \ref{bur.ex} at time $T=1/\pi$.}
\label{table.burl}
\end{table}
\begin{table}
\centering 
\small
\begin{tabular}{|l| c|c| c| c| c|  c| c| r|}
\hline 
 $N$ &WENO-AO(5,3) & Order & WENO-AON(5,3) & Order & WENO-AO(5,4,3) & Order \\ 
\hline
20  &  1.3928e-03 & --     &1.3654e-03 &-- &1.405610e-03 &\\ 
\hline 
40  &  7.0010e-05 & 4.31     &7.0933e-05 &4.27  &6.9975e-05 & 4.33    \\
\hline
80  &  2.3757e-06 & 4.87     &2.3795e-06 &4.88  &2.3754e-06 & 4.87   \\
\hline
160  &  6.9663e-08 & 5.09    &6.9663e-08 &5.09  &6.9663e-08 & 5.09    \\
\hline
320  &  2.1142e-09 & 5.04     &2.1142e-09 &5.04  &2.1142e-09 & 5.04  \\
\hline
\end{tabular}
\caption{  Comparison of $L^{1}$-errors for  WENO-AO(5,3), WENO-AON(5,3), and
WENO-AO(5,4,3) schemes  along with their convergence rates  for Example \ref{bur.ex} at time $T=1/\pi$.}
\label{table.bur1}
\end{table}
\begin{example}\label{bur.ex}{\rm (Burgers' Equation)
Consider the Burgers' equation given by
\begin{equation}
 u_t+\Big(\frac{u^2}{2}\Big)_x=0,~~~~x\in(-1,1)
\end{equation}
subject to the initial data $u(x,0)=0.25+0.5\sin(\pi x)$ with periodic boundary conditions. The time step in this case is taken to be 
$\Delta t=\mbox{0.5}\Delta x^{1.25}$. The numerical solutions are computed at time $T=1/\pi$. In Tables \ref{table.burl}, \ref{table.bur1}, we 
 have shown the $L^{\infty}$ and $L^1$-errors respectively of the WENO-AO(5,3), WENO-AON(5,3), and
WENO-AO(5,4,3) schemes  along with their convergence rates. The magnitude of errors of the schemes are comparable in the both norms and all 
 considered 
 schemes converge to solution with convergence rate five.

} 
\end{example}
\begin{table}
\centering 
\small
\begin{tabular}{|l| c|c| c| c| c|  c| c| r|}
\hline 
 $N$ & WENO-AO(5,3) & Order & WENO-AON(5,3) & Order & WENO-AO(5,4,3) & Order \\
\hline
20  &  2.326380e-05 & --     &2.3544e-05 &-- &2.3264e-05 &\\ 
\hline 
40  &  7.4180e-07 & 4.97     &7.4265e-07 &4.99 &7.4180e-07 & 4.97   \\
\hline
80  &  2.3343e-08 & 4.99     &2.3345e-08 &4.99  &2.3343e-08 & 4.99   \\
\hline
160  &  7.3390e-10 & 4.99    &7.3390e-10 &4.99  &7.3389e-10 & 4.99   \\
\hline
320  &  2.3084e-11 & 4.99    &2.3086e-11 &4.99  &2.3086e-11 & 4.99   \\
\hline
\end{tabular}
\caption{ Comparison of $L^{\infty}$-errors of the WENO-AO(5,3), WENO-AON(5,3), and
WENO-AO(5,4,3) schemes  along with their convergence rates  for Example \ref{euler.1d} at time $T=1$.}
\label{table.eul1d1}
\end{table}
\begin{table}[t!]
\centering 
\small
\begin{tabular}{|l| c|c| c| c| c|  c| c| r|}
\hline 
 $N$ & WENO-AO(5,3) & Order & WENO-AON(5,3) & Order & WENO-AO(5,4,3) & Order \\
\hline
20  &  9.3105e-05 & --     &9.3069e-05 &-- &9.3105e-05 &\\ 
\hline 
40  &  2.9632e-06 & 4.97    &2.9632e-06 &4.97  &2.9632e-06 & 4.97   \\
\hline
80  &  9.3446e-08 & 4.99    &9.3446e-08 &4.99  &9.3446e-08 & 4.99    \\
\hline
160  &  2.9355e-09 & 4.99    &2.9355e-09 &4.99  &2.9355e-09 & 4.99    \\
\hline
320  &  9.2337e-11 & 4.99    &9.2335e-11 &4.99  &9.2336e-11 & 4.99  \\
\hline
\end{tabular}
\caption{ Comparison of $L^{1}$-errors of the WENO-AO(5,3), WENO-AON(5,3), and
WENO-AO(5,4,3) schemes  along with their convergence rates  for Example \ref{euler.1d} at time $T=1$.}
\label{table.eul1d2}
\end{table}
\begin{figure}[t!]
 \includegraphics[width=0.48\textwidth]{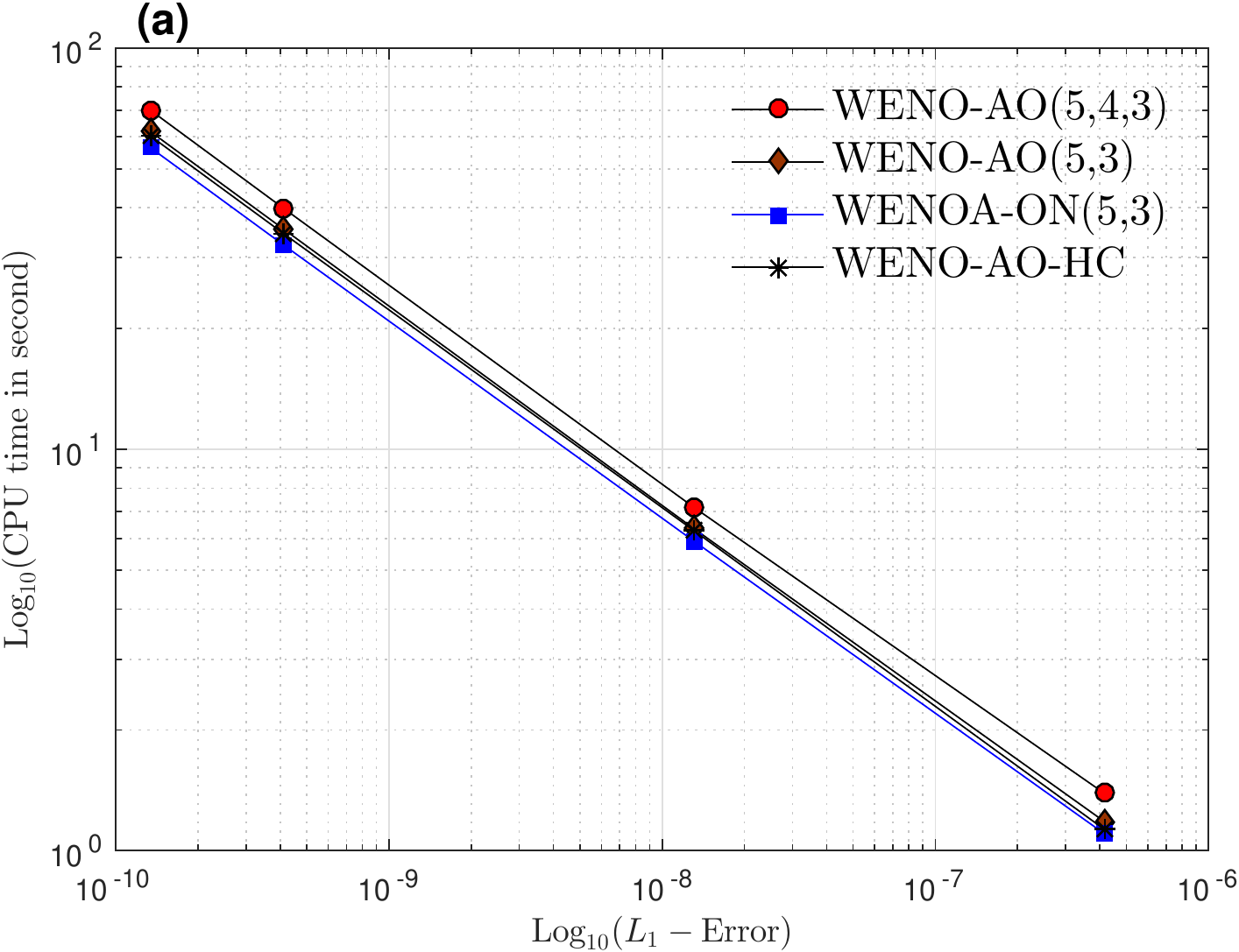}
 \includegraphics[width=0.48\textwidth]{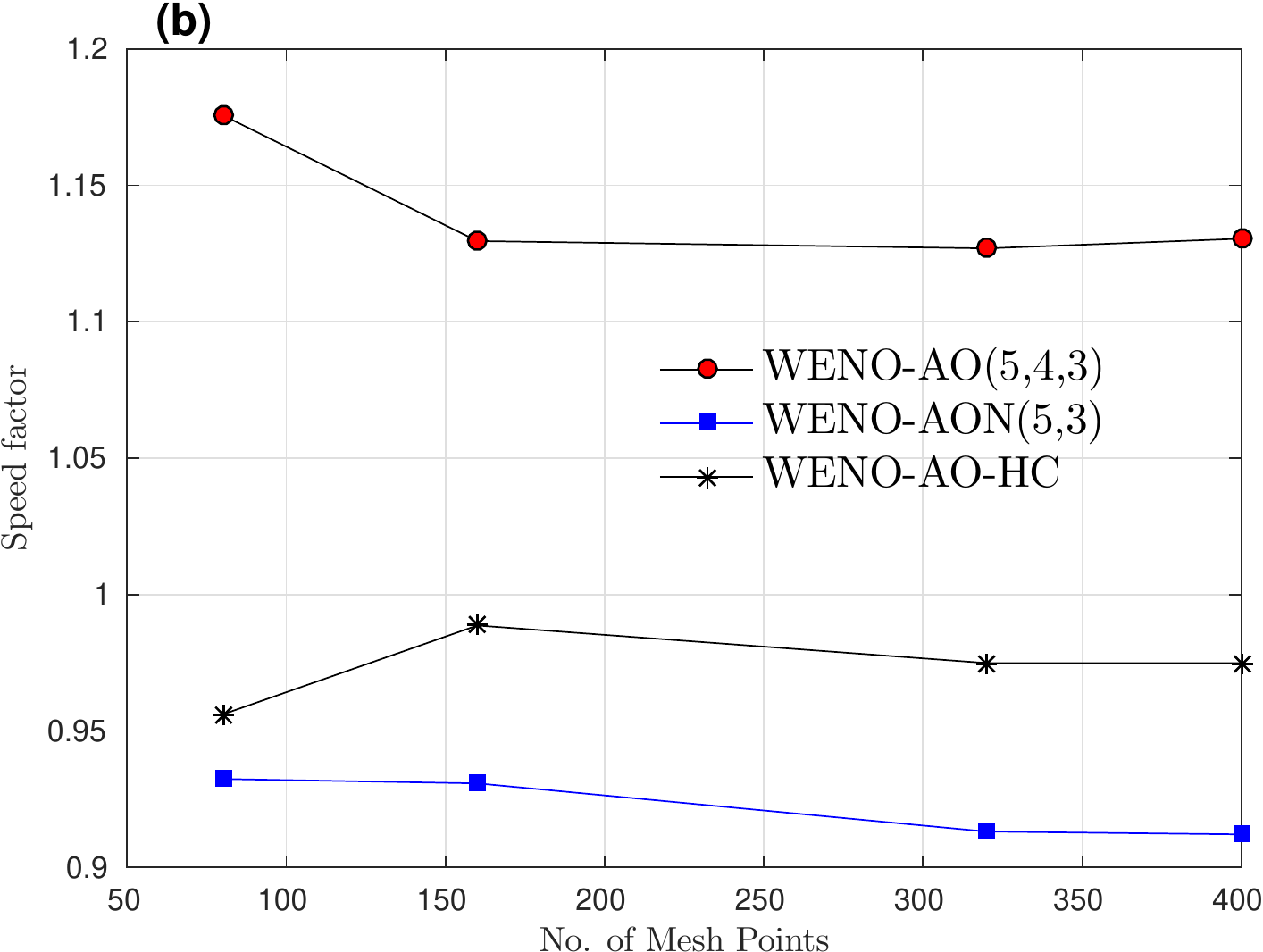}
 \caption{(a) Comparison of WENO schemes for Example \ref{euler.1d} in computing  $L^1-$error and CPU time, (b) relative CPU time taken by
 WENO-AO-HC, WENO-AON(5,3), and
WENO-AO(5,4,3) schemes with respect to WENO-AO(5,3) for different number of mesh points.}
 \label{speed.euler1d}
\end{figure}
\begin{example}\label{euler.1d}{\rm (1d Euler Equations)
The one-dimensional Euler system of equations  in the conservation form is
\begin{equation}\label{Euler.system}
\bold{U}_t+\bold{F(U)}_x=0,\; \; (x,t) \in \mathbb{R}\times(0,T],
\end{equation}
where $\bold{U}$ and $\bold{F(U)}$ are vector of conservative variables and fluxes, respectively  given by\\
$$\textbf{U}=\begin{bmatrix}
    \rho \\
\rho u \\
E
   \end{bmatrix},
\hspace{1.0 cm}
\textbf{F(U)}=\begin{bmatrix}
       \rho u\\
\rho u^2+p\\
(E+p)u
      \end{bmatrix}.$$
       Here, the dependent variables $\rho$, $p$, and $u$ denote density, pressure and velocity, respectively. The energy $E$ is given by the relation
 \begin{equation}
E=\frac{p}{(\gamma-1)} +\frac{1}{2}\rho u^2,
\end{equation}
where $\gamma$ is the ratio of specific heats and taken as 1.4 in test cases or it will be mentioned separately.
}
\end{example}
In order to measure the accuracy of the schemes, we consider the Euler system of equations with the initial data $\rho(x,0) =1+0.2 \sin(x)$, $u(x,0)=1$,
and $p(x,0)=1$ over the domain $[0, 2\pi]$. Numerical solutions are computed at time $T=1$. The time step is
$\Delta t=0.5\Delta x^{1.5}$. In Table \ref{table.eul1d1}, \ref{table.eul1d2}, we have depicted the  $L^{\infty}$- and $L^1$-errors
and convergence rates of the 
numerical schemes. We can easily observe from the Tables that schemes are comparable in terms of accuracy and convergence rate. In Figure 
\ref{speed.euler1d} (a), 
 we have shown the loglog plot of $L^1$- errors and time taken to compute the error. In Figure \ref{speed.euler1d}(b), we have shown the speed factor
  of the scheme. The WENO-AO(5,4,3) scheme takes 12 percent more time to compute solution. The WENO-AON(5,3) and  WENO-AO-HC schemes take less 
  computational time in comparison to WENO-AO(5,3) scheme.
\begin{table}
\centering 
\small
\begin{tabular}{|c| c|c| c| c| c|  c| c| c|}
\hline 
 $N_x\times N_y$ & WENO-AO(5,3) & Order & WENO-AON(5,3) & Order & WENO-AO(5,4,3) & Order \\
\hline
$10 \times 10$  &  1.5001e-03 & --     &1.5267e-03 &-- &1.4680e-03 &\\ 
\hline 
$20 \times 20$  &  5.0152e-05 & 4.90    &5.0652e-05 &4.91 &5.0128e-05 & 4.87   \\
\hline
$40 \times 40$ &  1.6676e-06 & 4.91     &1.6692e-06 &4.92 &1.6676e-06 & 4.91   \\
\hline
$80 \times 80$  &  5.4836e-08 & 4.93    &5.4839e-08 &4.93  &5.4836e-08 & 4.93  \\
\hline
$160 \times 160$ &  1.8285e-09 & 4.91    &1.8286e-09 &4.91  &1.8285e-09 & 4.91   \\
\hline
\end{tabular}
\caption{  Comparison of $L^{\infty}$-errors for WENO schemes   along with their convergence rates  for Example \ref{2deuler.crate} 
over the domain $[0,2\pi]\times [0,2\pi]$ at time $T=2$.}
\label{table.eul2d1}
\end{table}

\begin{table}
\centering 
\small
\begin{tabular}{|c| c|c| c| c| c|  c| c| c|}
\hline 
 $N_x\times N_y$ & WENO-AO(5,3) & Order & WENO-AON(5,3) & Order & WENO-AO(5,4,3) & Order \\
\hline
$10 \times 10$  &  1.5001e-03 & --     &1.5267e-03 &-- &1.4680e-03 &\\ 
\hline 
$20 \times 20$  &  1.2825e-03 & 4.92    &1.2818e-03 &4.94  &1.2797e-03 & 4.88     \\
\hline
$40 \times 40$  &  4.1891e-05 & 4.93    &4.1891e-05 &4.94  &4.1886e-05 & 4.93    \\
\hline
$80 \times 80$  &  1.3791e-06 & 4.92     &1.3791e-06 &4.92 &1.3791e-06 & 4.92    \\
\hline
$160 \times 160$ &  4.5952e-08 & 4.91     &4.5952e-08 &4.91  &4.5952e-08 & 4.91  \\
\hline
\end{tabular}
\caption{ Comparison of $L^{1}$-errors for WENO schemes   along with their convergence rates  for Example \ref{2deuler.crate} 
over the domain $[0,2\pi]\times [0,2\pi]$ at time $T=2$.}
\label{table.eul2d2}
\end{table}

\begin{figure}[t!]
 \includegraphics[width=0.48\textwidth]{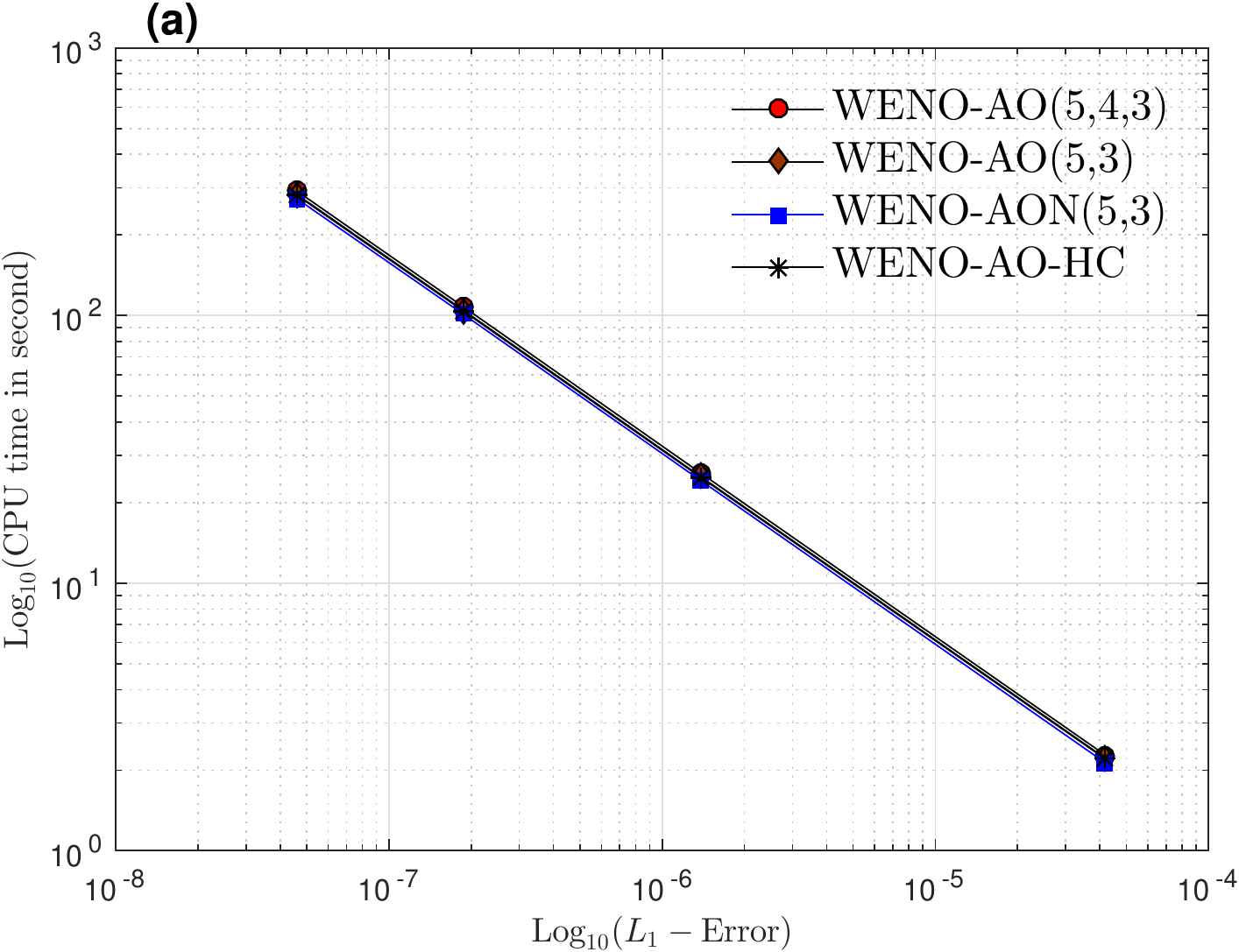}
 \includegraphics[width=0.48\textwidth]{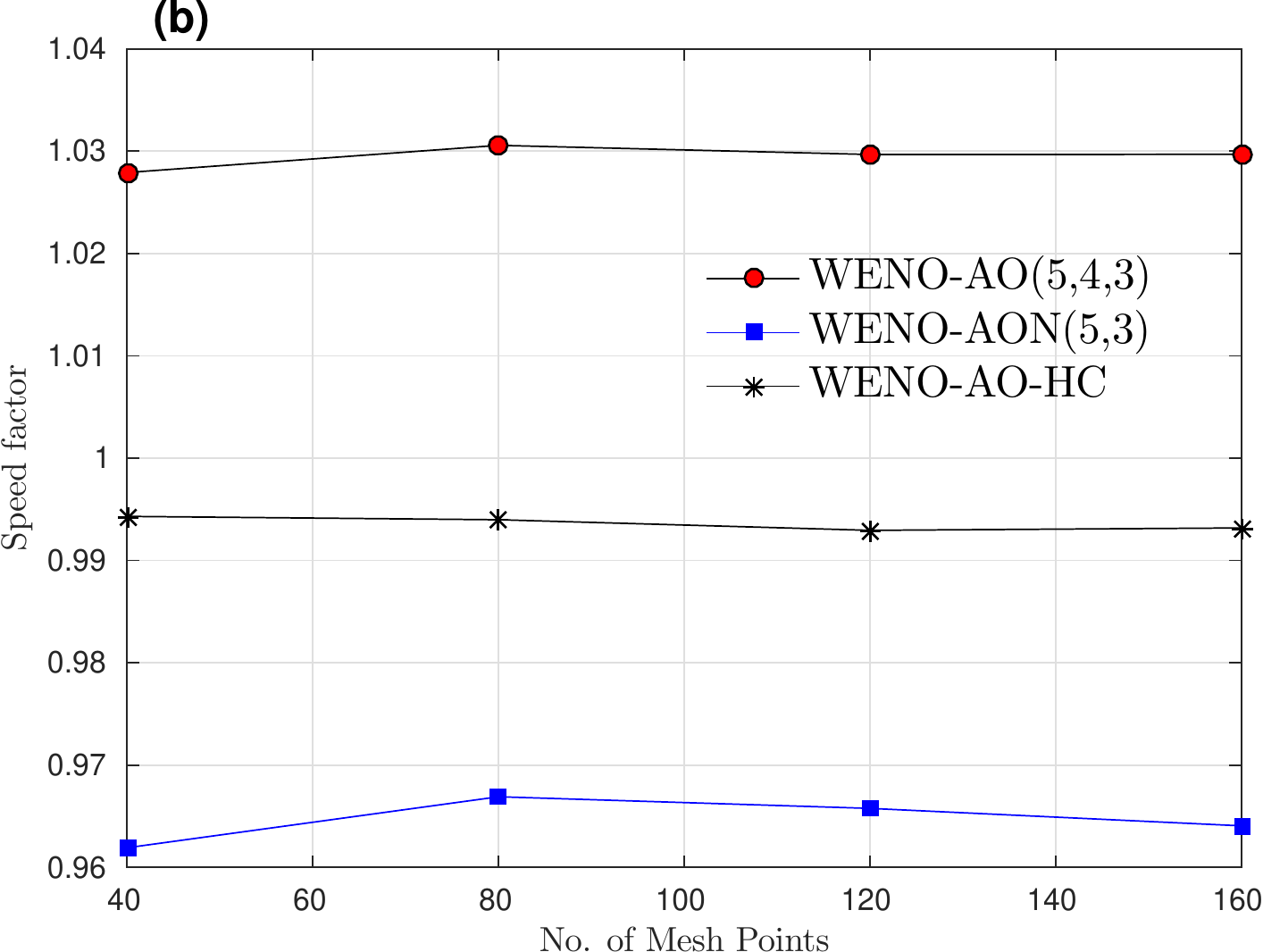}
\caption{(a) Comparison of WENO schemes for Example \ref{2deuler.crate} in computing  $L^1-$error and CPU time, (b) relative CPU time taken by
 WENO-AO-HC, WENO-AON(5,3), and
WENO-AO(5,4,3) schemes with respect to WENO-AO(5,3) for different number of mesh points.}
 \label{speed.euler2d}
\end{figure}

\begin{example}\label{2deuler.crate}{\rm (2d-Euler Equations)
 The two-dimensional Euler system of equations in the conservation form is
\begin{equation}\label{Euler.system2d}
\bold{U}_t+\bold{F(U)}_x+\bold{G(U)}_y =0,\; \; (x,y,t) \in \Omega \times(0,T],
\end{equation}
where $\bold{U}$,  $\bold{F(U)}$, and $\bold{G(U)}$ are vectors of conservative variables and fluxes, given respectively by
$$\bold{U}=\begin{bmatrix}
    \rho \\
\rho u \\
\rho v \\
E
   \end{bmatrix},
\hspace{1.0 cm}
\bold{F(U)}=\begin{bmatrix}
       \rho u\\
\rho u^2+p\\
\rho u v\\
(E+p)u
      \end{bmatrix},
      \hspace{1.0 cm}
\bold{G(U)}=\begin{bmatrix}
       \rho v\\
       \rho u v\\
\rho v^2+p\\
(E+p)u
      \end{bmatrix}.$$
      
 }
\end{example}
To measure the accuracy in this case, we consider  equation \eqref{Euler.system2d} with the initial data 
\begin{align}
 \rho(x,y,0)=1+0.2\sin(x+y), ~~~u(x,y,0)=v(x,y,0)=p(x,y,0)=1
\end{align}
over the domain $[0,2\pi]\times [0,2\pi] $. Numerical solutions are computed at time $T=2$ and errors are depicted in Tables \ref{table.eul2d1} $\&$
\ref{table.eul2d2}. All the schemes achieve their designed order of accuracy and their errors are comparable in both the norms. The loglog plot 
between $L^1$-errors of schemes and CPU time to compute solution is depicted in Figure \ref{speed.euler2d}(a), whereas in Figure \ref{speed.euler2d} (b),
 we have shown the speedup factor with number of mesh points. The WENO-AO(5,4,3) uses approximately 3$\%$ more computational time than WENO-AO(5,3),
  whereas computational cost of WENO-AON(5,3), WENO-AO-HC are comparable with WENO-AO(5,3). 
  
 We can observe from these test cases that  errors of WENO-AO(5,3), WENO-AON(5,3), WENO-AO(5,4,3) schemes are comparable. We can also observe the drastic change in
  speedup factor of WENO-AO(5,4,3) scheme as we move from scalar to system of equations in 1d and 2d.

\begin{figure}
\begin{center}
 \begin{tabular}{cc}
 \includegraphics[width=0.48\textwidth]{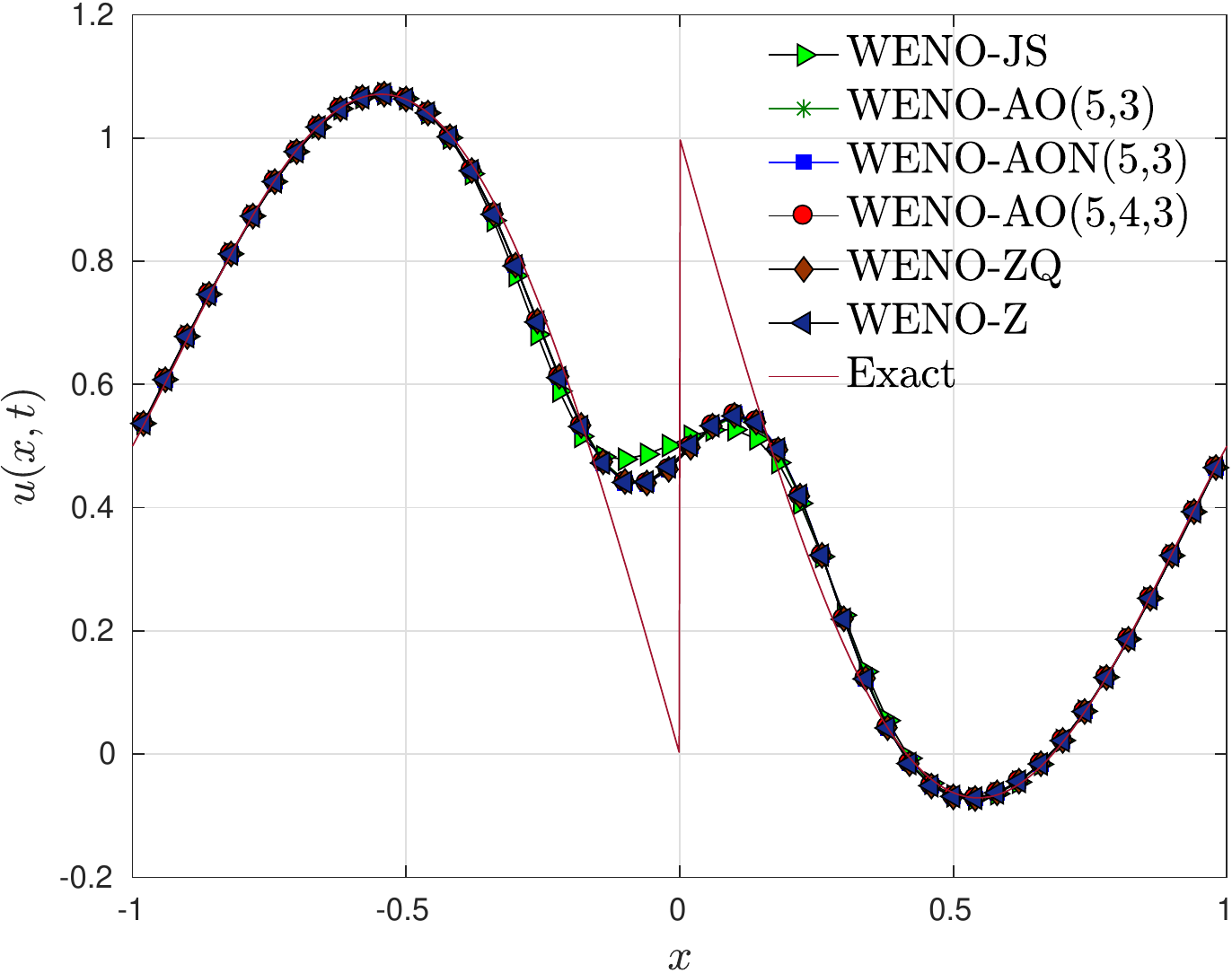} &
 \includegraphics[width=0.48\textwidth]{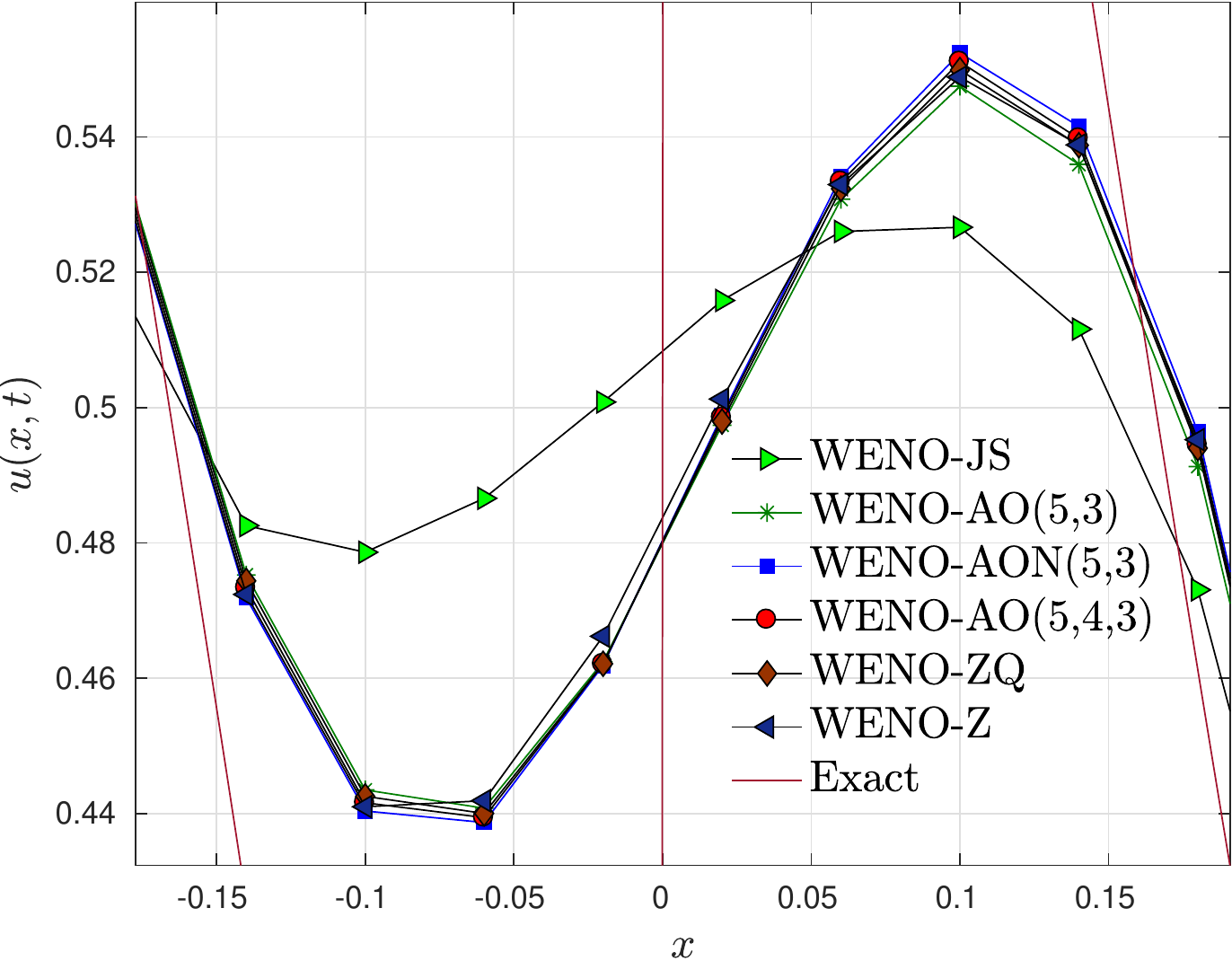}\\
 (a)  & (b) 
 \end{tabular}
\end{center}
  \caption{(a) Numerical solution of linear advection equation \eqref{lin.adv} with discontinuous initial data \eqref{dis.adv} at $T=8$ for 
 the WENO-JS, WENO-AO(5,3), WENO-AON(5,3), WENO-AO(5,4,3), WENO-ZQ, and WENO-Z scheme using 50 mesh points. The exact solution is shown in a solid line.
  (b) zoom version of (a) near discontinuity.}
 \label{fig.dis50}
\end{figure}

\begin{figure}
\begin{center}
 \begin{tabular}{cc}
 \includegraphics[width=0.48\textwidth]{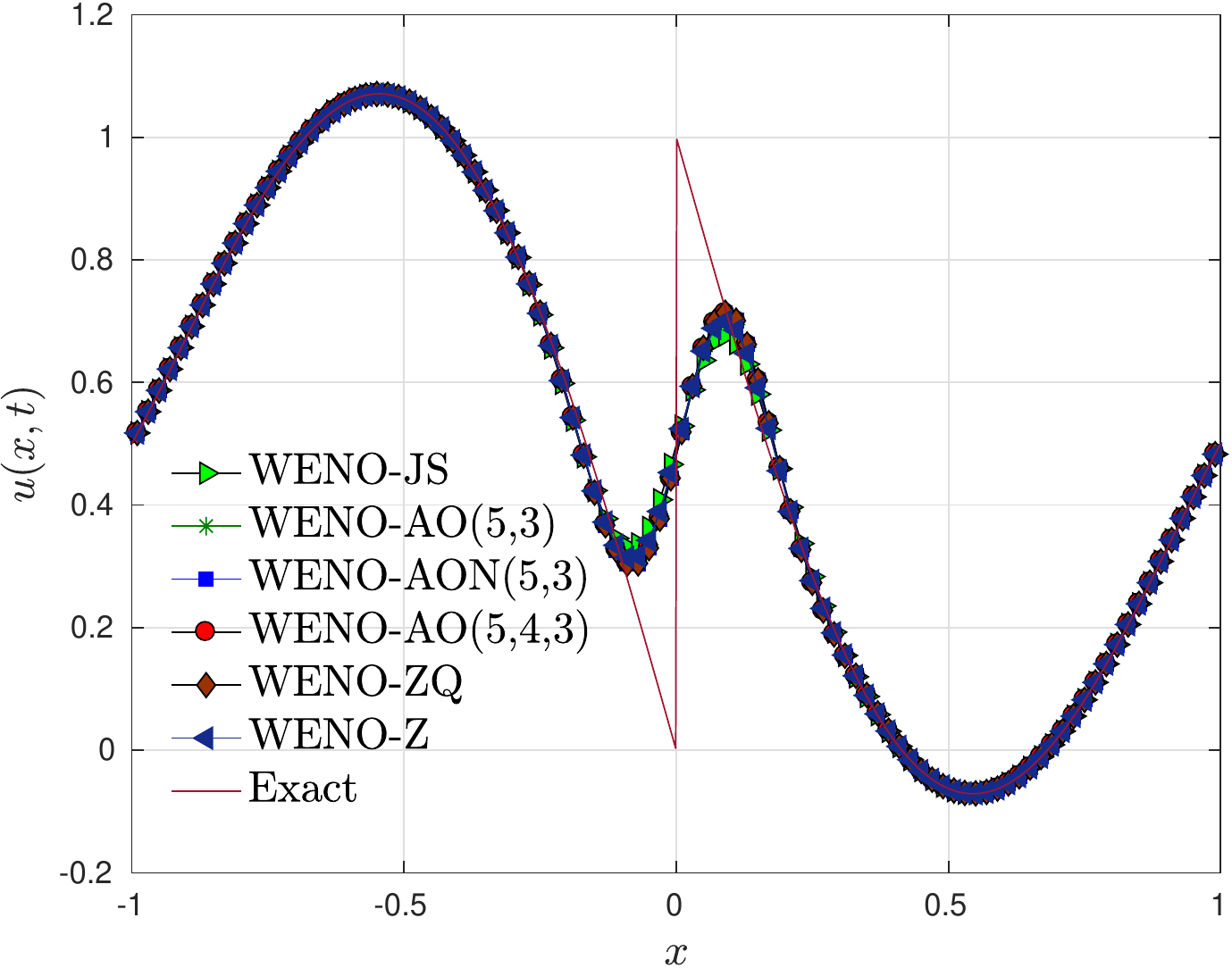} &
 \includegraphics[width=0.48\textwidth]{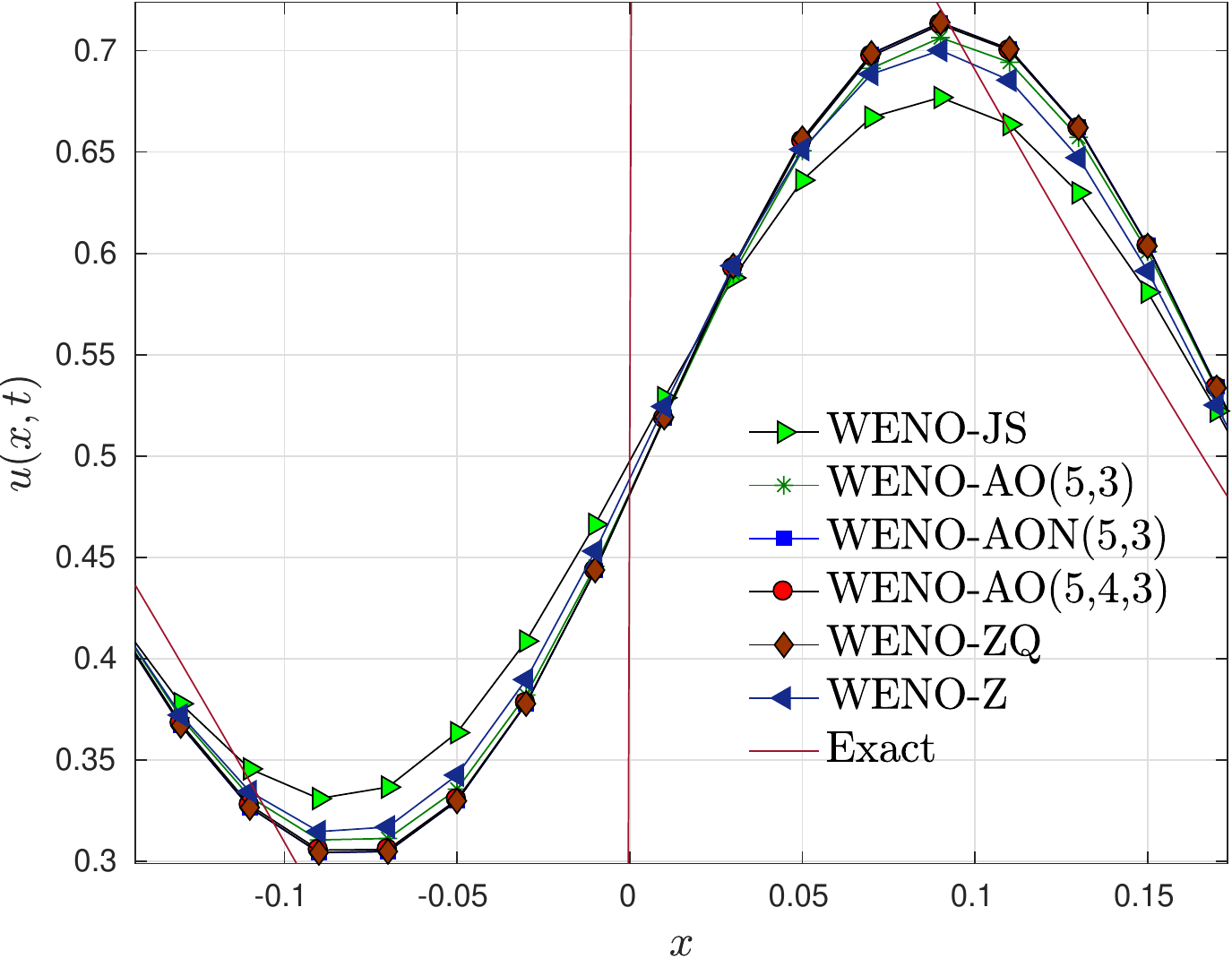}\\
 (a)  & (b) 
 \end{tabular}
\end{center}
 \caption{(a) Numerical solution of linear advection equation \eqref{lin.adv} with discontinuous initial data \eqref{dis.adv} at $T=8$ for 
 the WENO-JS, WENO-AO(5,3), WENO-AON(5,3), WENO-AO(5,4,3), WENO-ZQ, and WENO-Z scheme using 100 mesh points. The exact solution is shown in a solid line.
  (b) zoom version of (a) near discontinuity.}
 \label{fig.dis100}
\end{figure}

\begin{figure}
\begin{center}
 \begin{tabular}{cc}
 \includegraphics[width=0.48\textwidth]{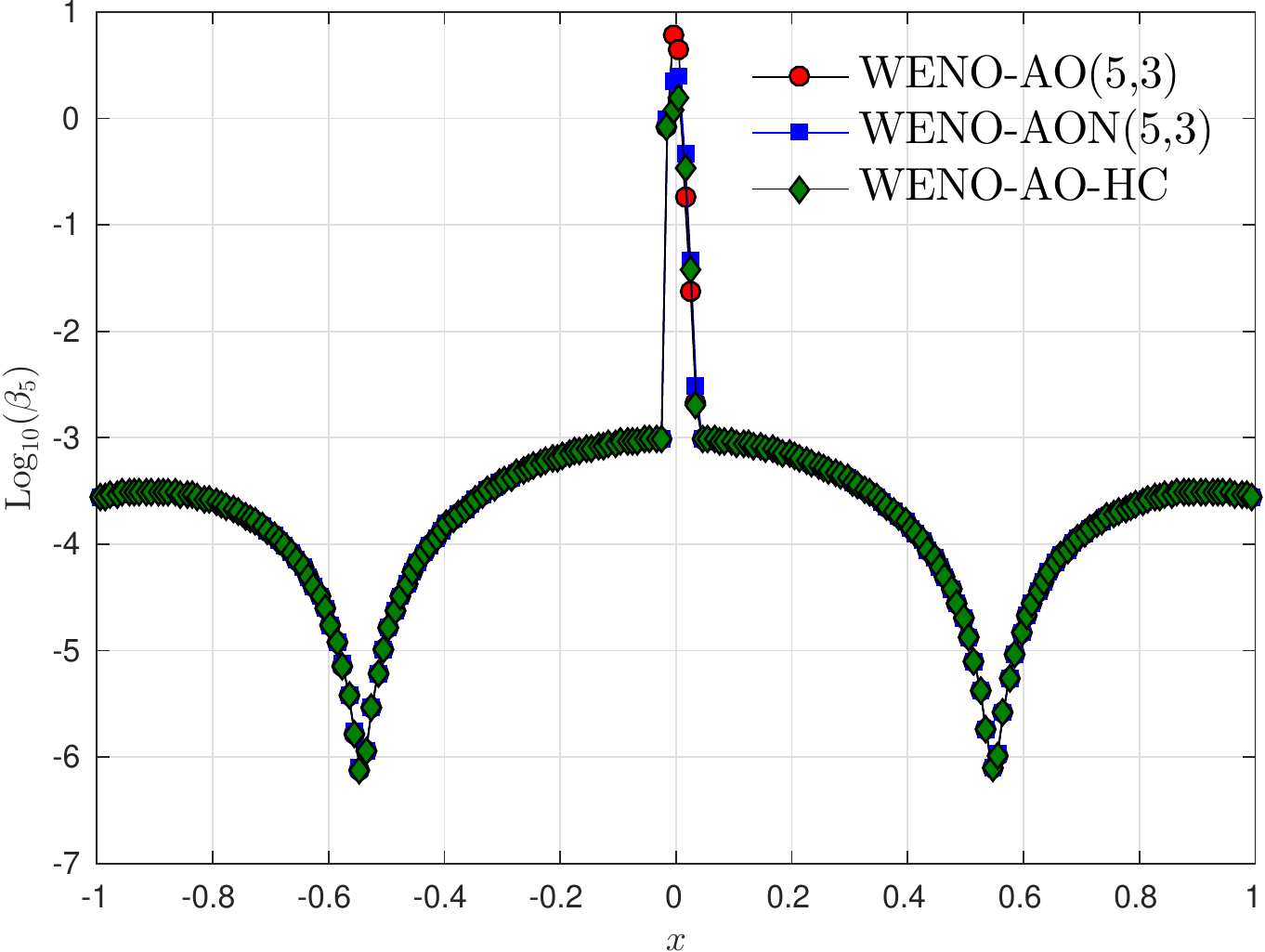} &
 \includegraphics[width=0.48\textwidth]{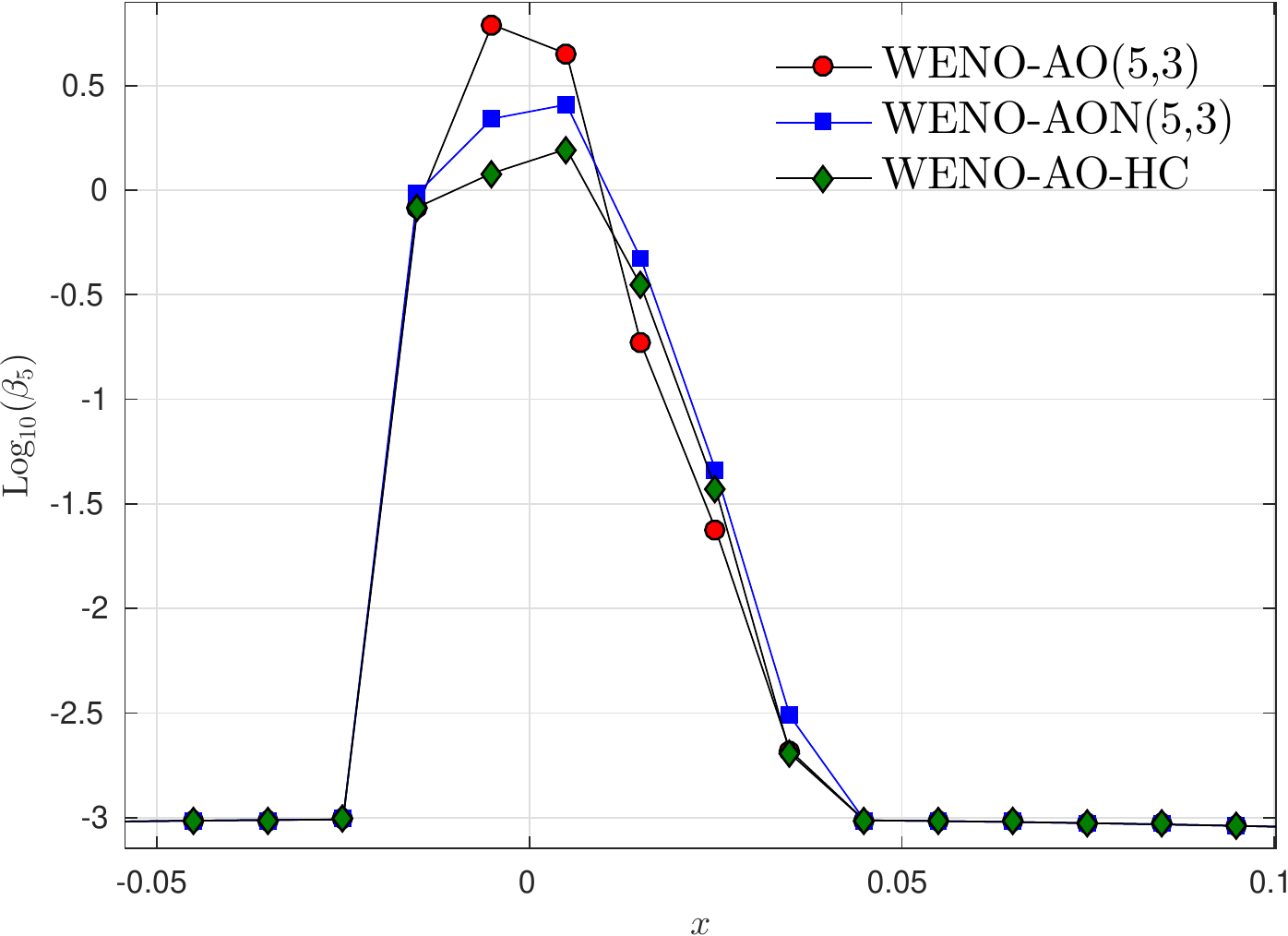}\\
 (a)  & (b) 
 \end{tabular}
\end{center}
\caption{(a) Comparison of values of $\beta_5$ obtained in case of WENO-AO(5,3), WENO-AON(5,3), and WENO-AO-HC schemes with $\epsilon = 10^{-40}$. 
 (b) Enlarge portion of (a) around discontinuity.}
 \label{fig.beta}
\end{figure}

\begin{figure}
\begin{center}
 \begin{tabular}{cc}
 \includegraphics[width=0.48\textwidth]{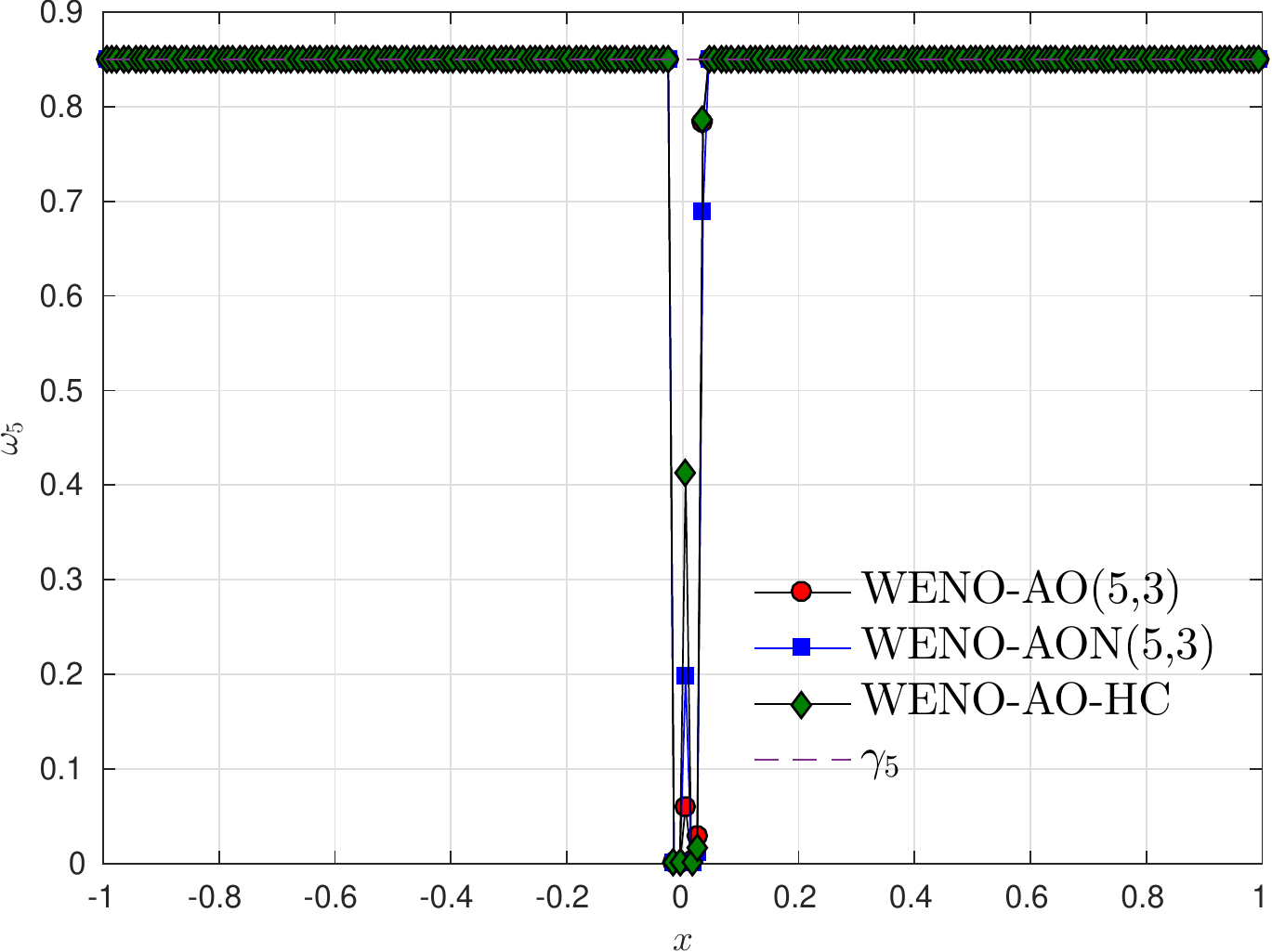} &
 \includegraphics[width=0.48\textwidth]{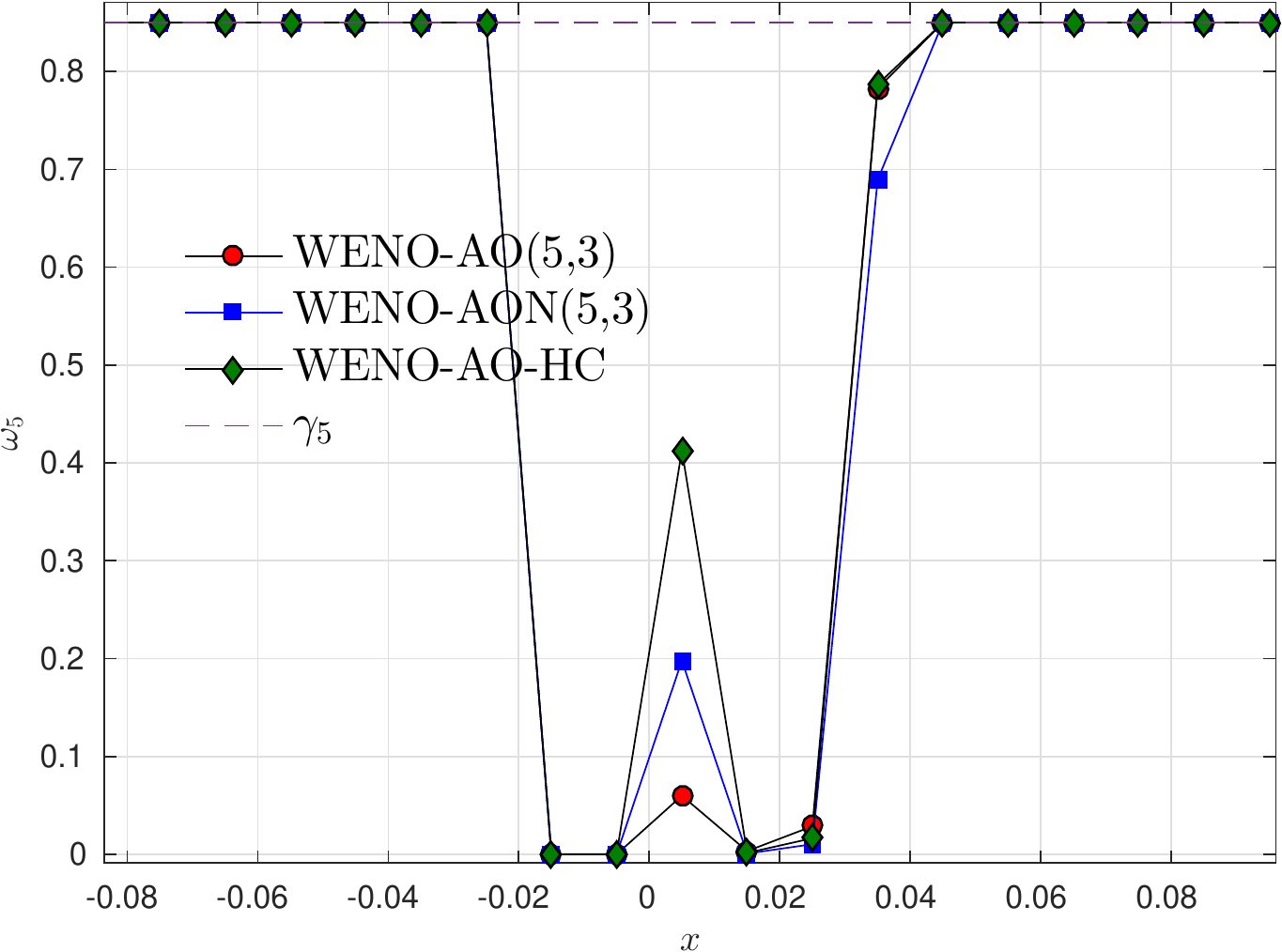}\\
 (a)  & (b) 
 \end{tabular}
\end{center}
\caption{(a)The distribution of nonlinear weights $\beta_5$  and linear weigths $\gamma_5$ for WENO-AO(5,3),  WENO-AON(5,3), and  WENO-AO-HC schemes at the
  first step of numerical solution of the linear advection equation \eqref{lin.adv} with discontinuous periodic initial condition \eqref{dis.adv}.
  (b) enlarge portion of (a) around discontinuity.}
 \label{fig.gamma}
\end{figure}
\subsection{Test Cases with Discontinuities}\label{sec:shock}
Now we test the algorithms of the  WENO-AON(5,3) and WENO-AO(5,4,3) schemes on  test cases having discontinuities and compare them  with the WENO-JS,
WENO-Z, WENO-ZQ, and 
WENO-AO(5,3) schemes.
\begin{example}{\rm
  Consider the linear advection equation \eqref{lin.adv} with following initial data \cite{bor-etal_08a}
 \begin{equation}\label{dis.adv}
u(x,0)=\left\{\begin{array}{ll}
                -(\sin(\pi x)+\frac{x^3}{2}), ~~~~   -1\leq x<0, \\
                -(\sin(\pi x)+\frac{x^3}{2})+1, ~~  0\leq x\leq 1.
                \end{array}\right.
\end{equation}
The problem contains a discontinuity at $x=0$ which advects, to the right. Numerical solutions are computed at $T=8$ with 50 and 100
 grid points using WENO schemes and comparison is shown in Figure \ref{fig.dis50}, \ref{fig.dis100}, respectively. From Figure \ref{fig.dis50}(a)-(b),
  we can observe that WENO schemes of adaptive order are comparable and resolves the discontinuity in least dissipative manner. The WENO-AON(5,3) perform 
  better than other scheme on coarse mesh  and  it was replaced by WENO-AO(5,4,3) scheme as the grid refined. Figure \ref{fig.beta} shows the value of smoothness
   indicator $\beta_5$ used in the WENO-AO(5,3), WENO-AON(5,3), and WENO-AO-HC  schemes computed at first step of simulation. We can observe from Figure
    \ref{fig.beta} (b)  that values of smoothness indicator used in WENO-AON(5,3) scheme is closer to that of WENO-AO(5,3) scheme than that of WENO-AO-HC scheme.
    In Figure \ref{fig.gamma}, we have compared the nonlinear weight $\omega_5$ for the  WENO-AO(5,3), WENO-AON(5,3), and WENO-AO-HC  schemes along
     with the linear weight $\gamma_5=0.85$. In all schemes, value of the nonlinear weight $\omega_5$ remains close to the linear weight $\gamma_5$ in smooth part and 
      tend to zero in the presence of shock. The value of $\omega_5$ in WENO-AON(5,3) scheme is  closer to that of WENO-AO(5,3) scheme than 
     value of $\omega_5$ in WENO-AO-HC scheme.
  
   }
\end{example}

\begin{figure}[t!]
\begin{center}
 \begin{tabular}{cc}
\includegraphics[width=0.48\textwidth]{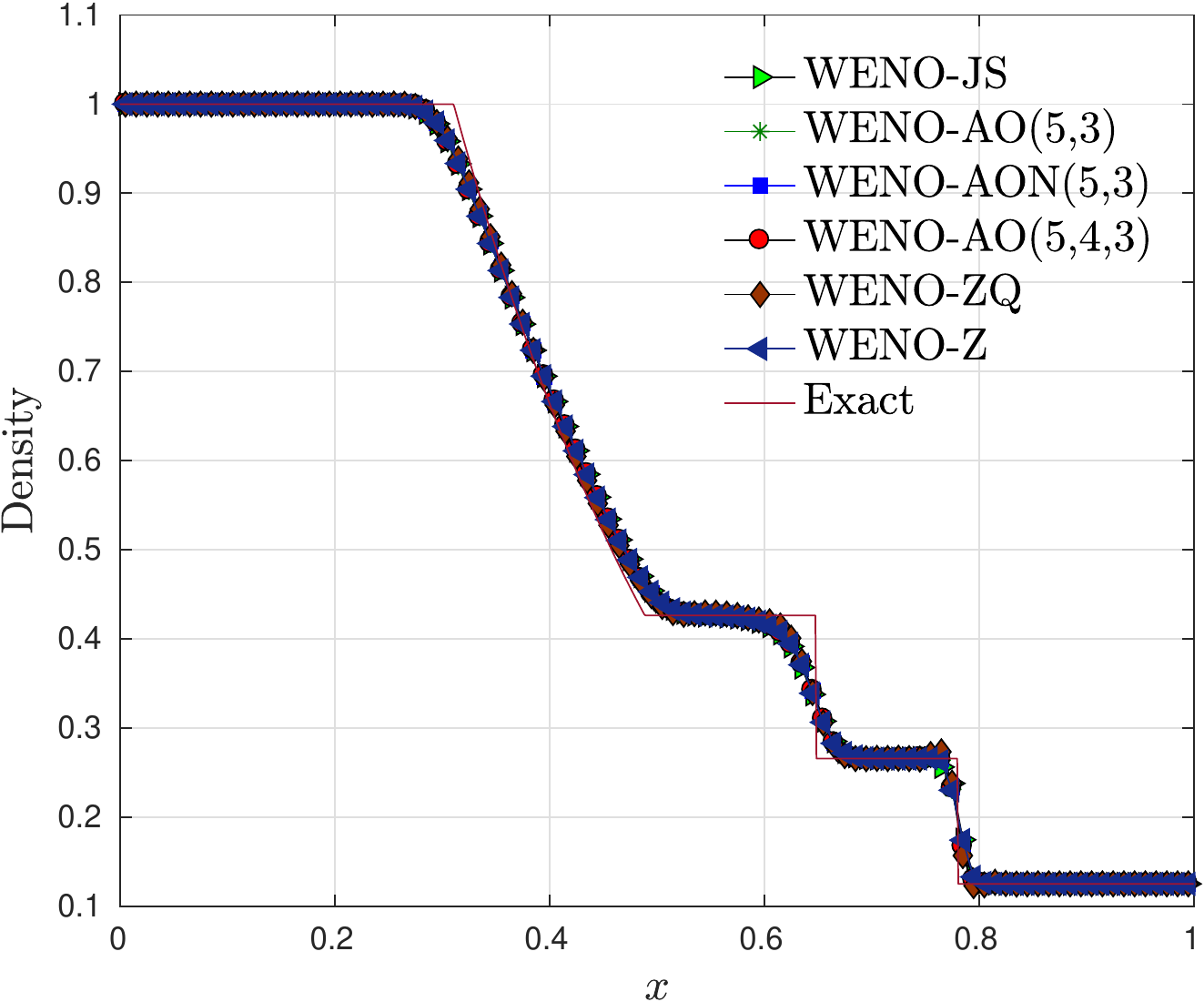}
\includegraphics[width=0.48\textwidth]{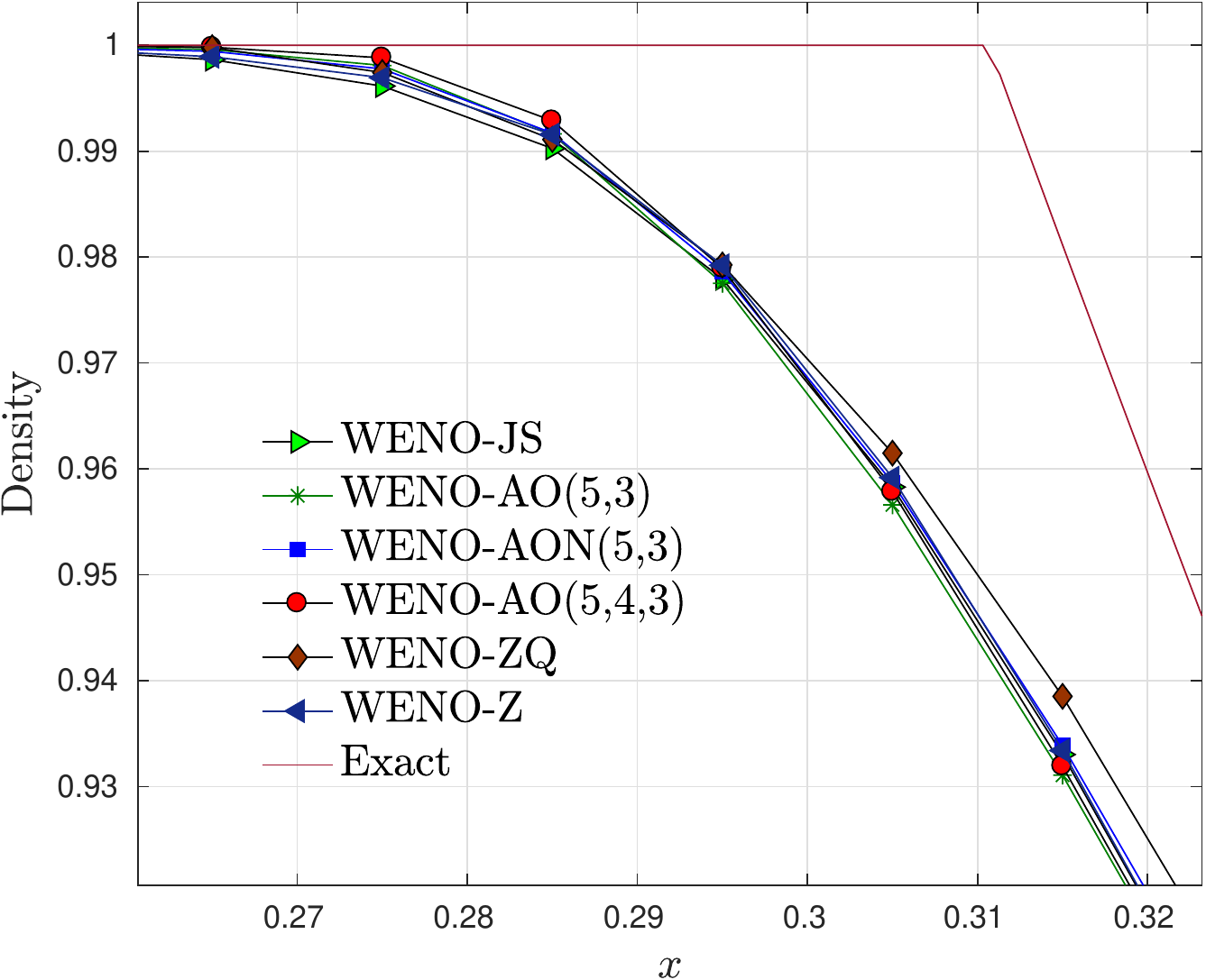}\\
\includegraphics[width=0.48\textwidth]{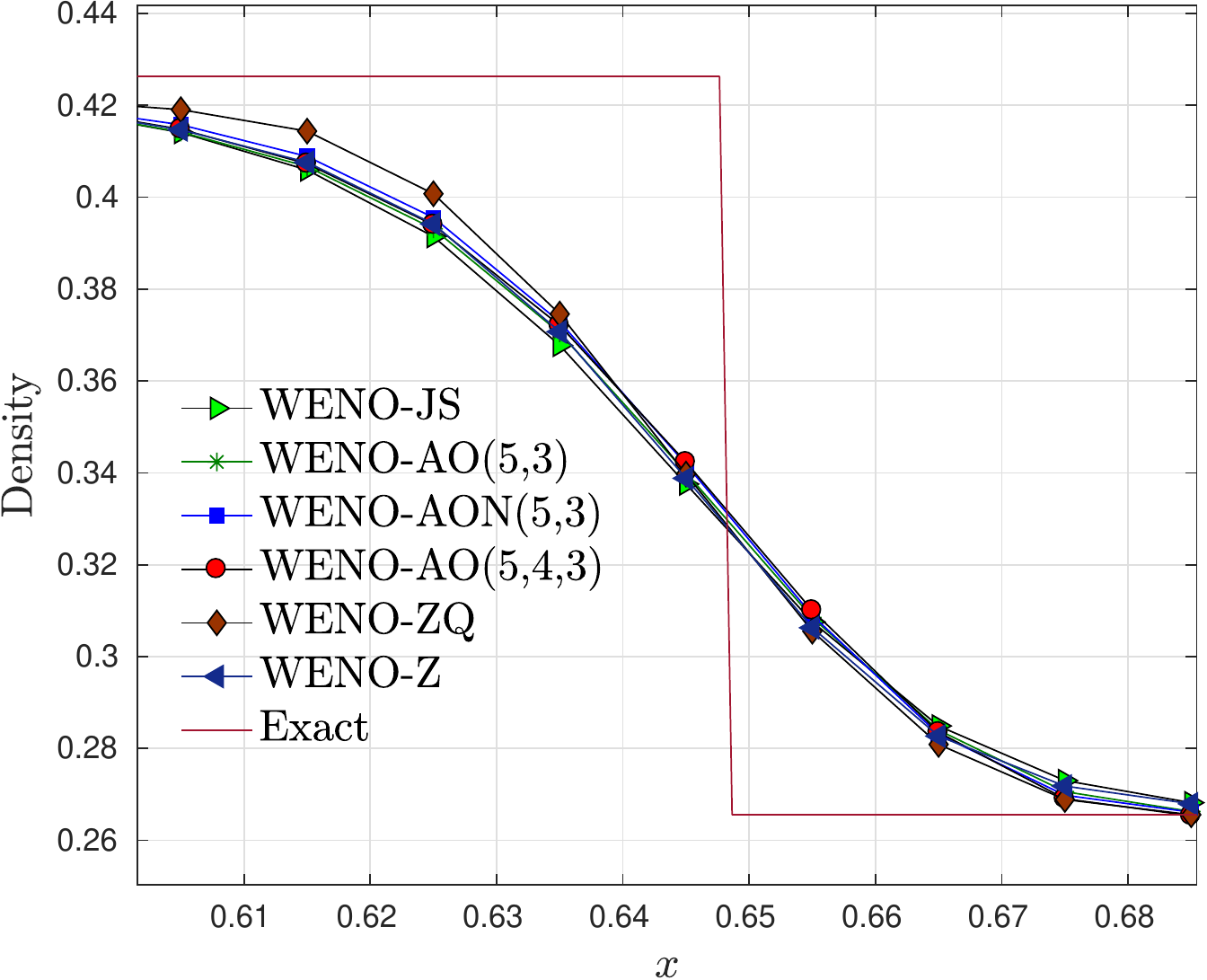}
\includegraphics[width=0.48\textwidth]{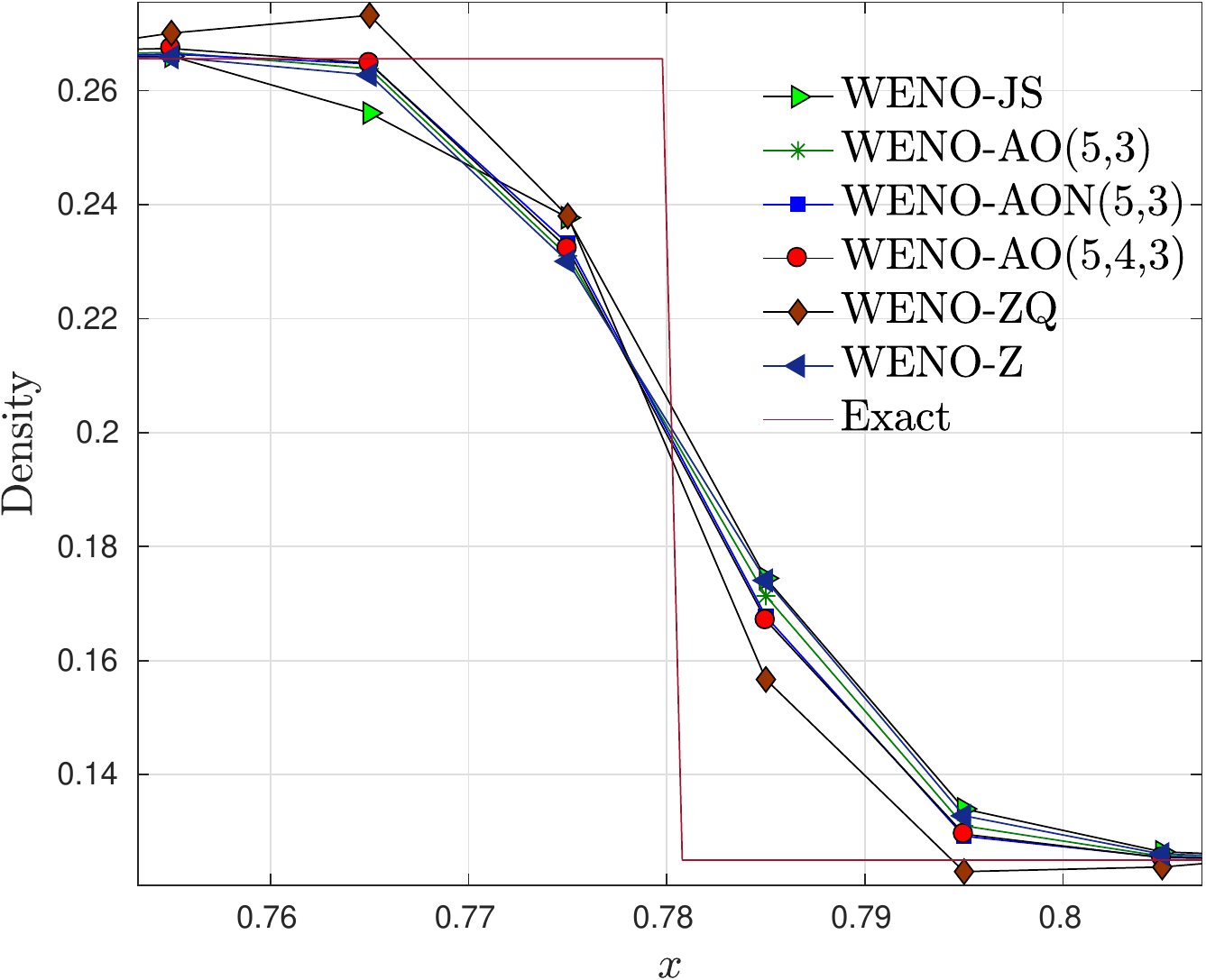}
 \end{tabular}
\end{center}
 \caption{Comparison of the WENO-JS, WENO-Z, WENO-ZQ, WENO-AO(5,3), WENO-AON(5,3), and WENO-AO(5,4,3) schemes for the Example \ref{shock_test}
 at time $T=0.16$ over a domain $[0,1]$ with a uniform mesh of having 100 points.}
 \label{Figure.Sod}
\end{figure}

\begin{figure}[t!]
\begin{center}
 \begin{tabular}{cc}
\includegraphics[width=0.48\textwidth]{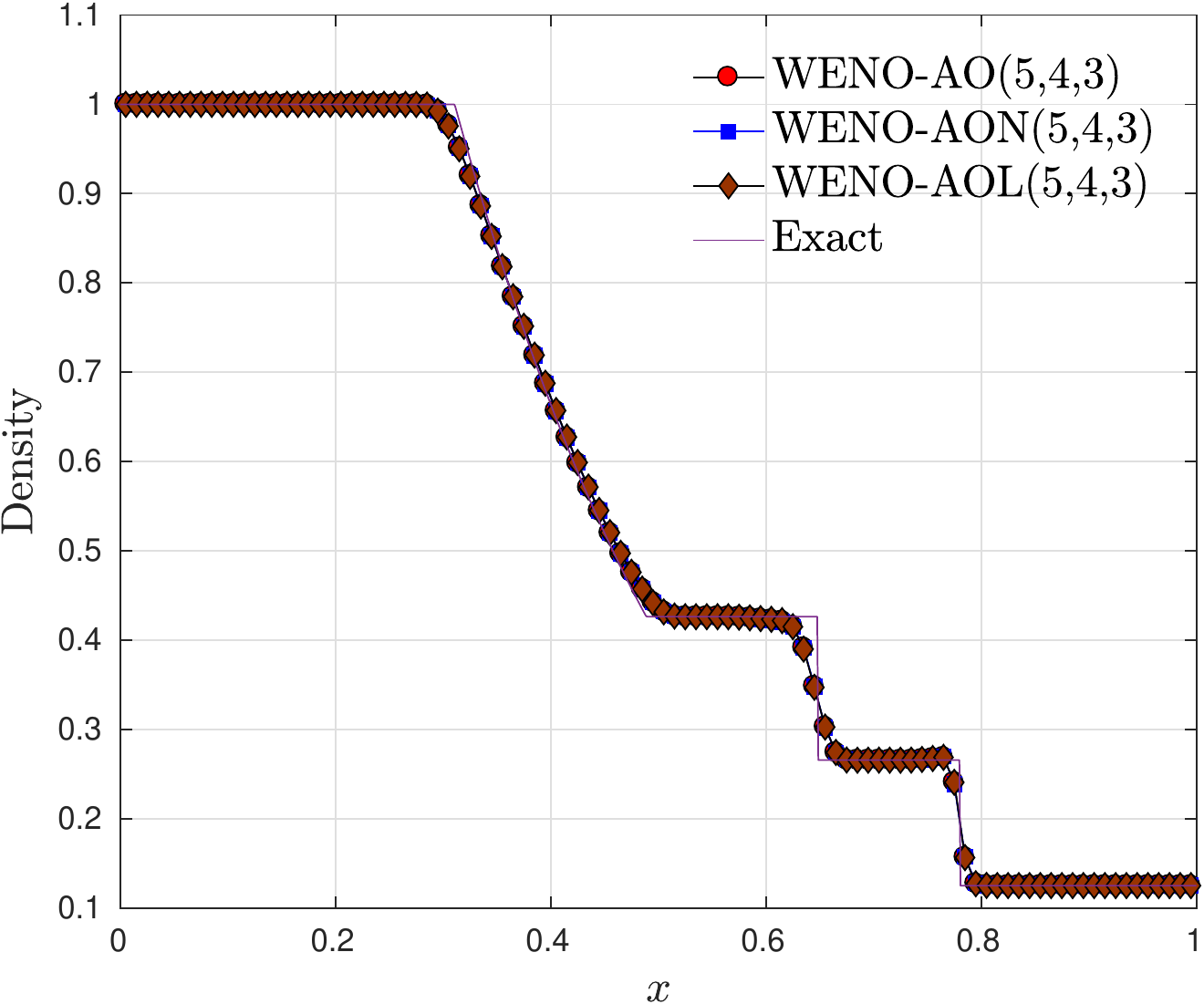}
\includegraphics[width=0.48\textwidth]{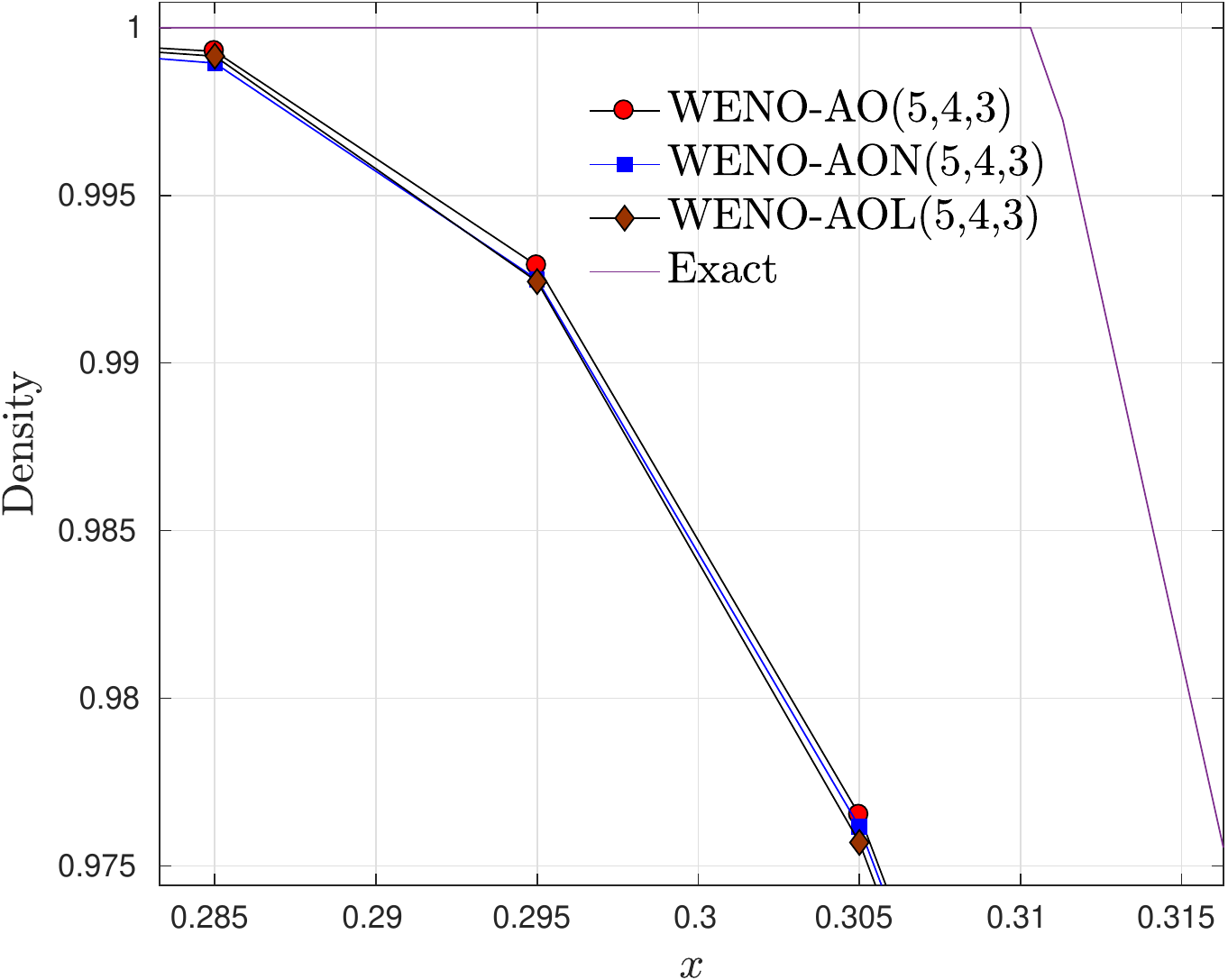}\\
\includegraphics[width=0.48\textwidth]{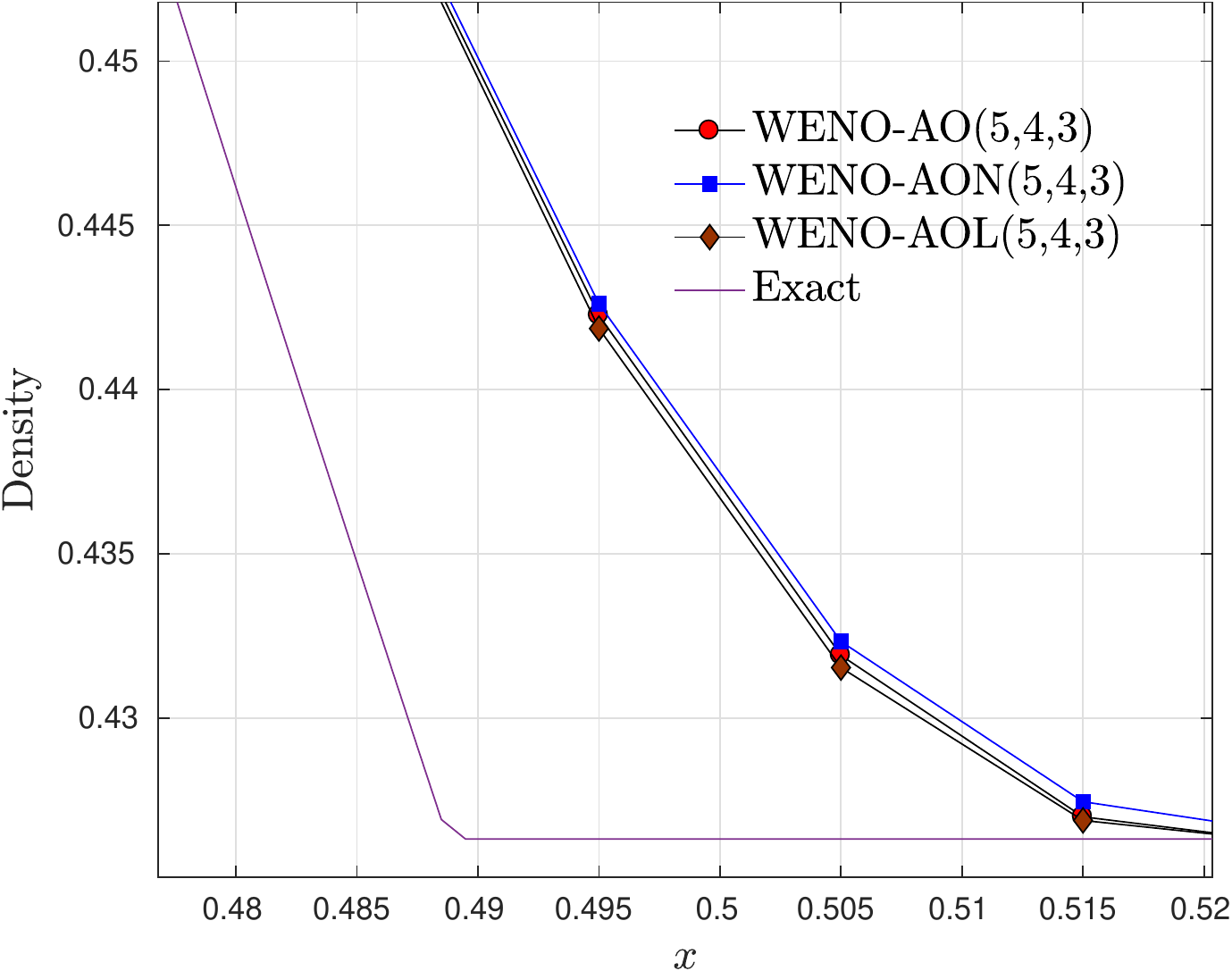}
\includegraphics[width=0.48\textwidth]{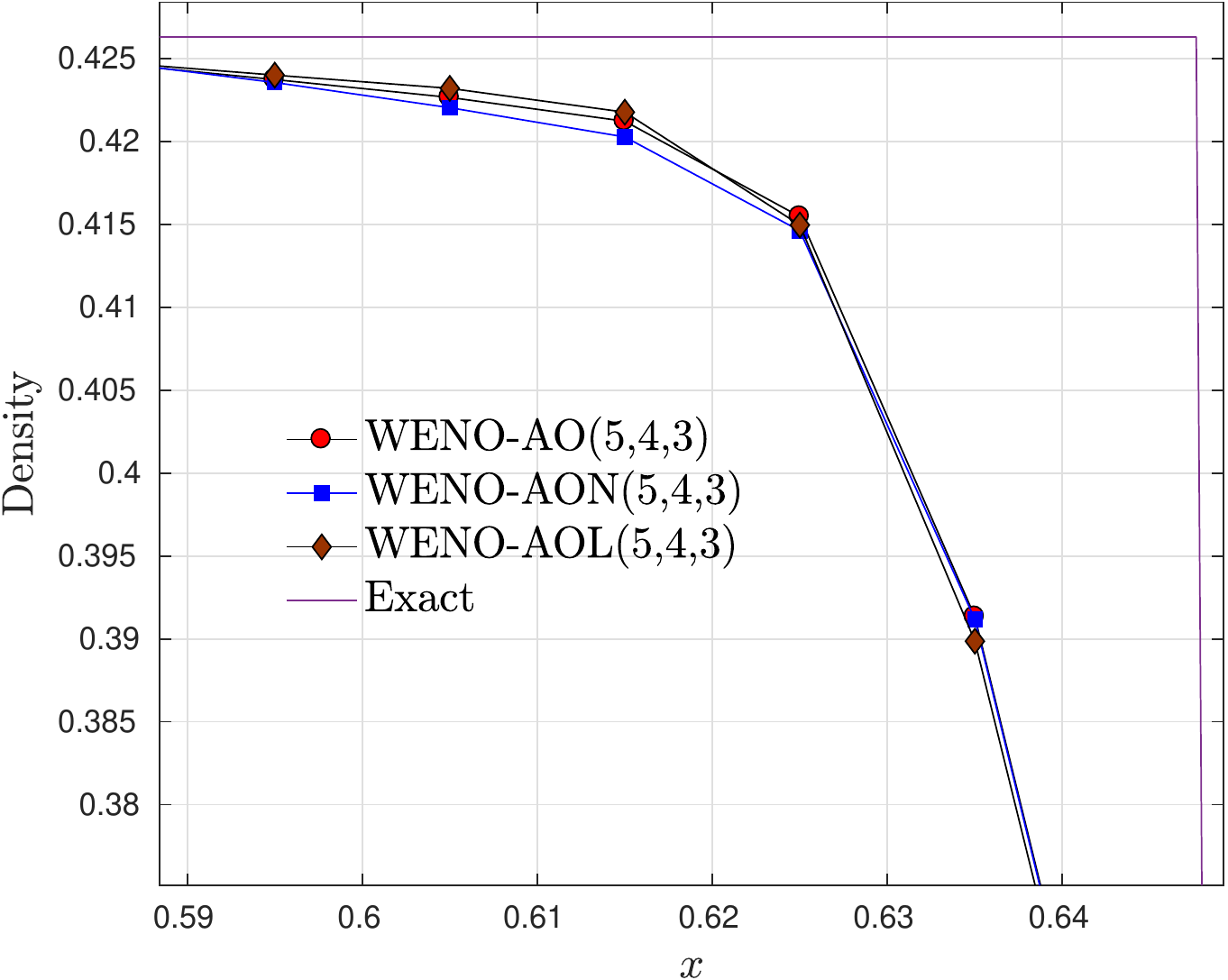}
 \end{tabular}
\end{center}
 \caption{Comparison of the WENO-AO(5,4,3), WENO-AOL(5,4,3), and WENO-AON(5,4,3) schemes for the Example \ref{shock_test}
 at time $T=0.16$ over a domain $[0,1]$ with a uniform mesh of having 100 points.}
 \label{Figure.Sod_comp}
\end{figure}

\begin{example}\label{shock_test}{\rm (Shock tube test)
This test consists of two constant states $\bold{U}_l$ and $\bold{U}_r$ separated with an initial discontinuity and  involves the formation of 
a rarefaction wave, a
contact discontinuity, and a shock wave  \cite{sod_78a}. 
 Here, we consider Euler equations \eqref{Euler.system} with the following initial conditions
\begin{equation}\label{Sod.system1}
(\rho, u, p)(x,0)=\left\{\begin{array}{ll}
                (1.0,0.0,1.0) ~~~~   x<0.5, \\
                (0.125,0.0,0.1) ~  x>0.5.
                \end{array}\right.
\end{equation}
over the domain $[0,1]$ with transmissive boundary conditions. Numerical solution are computed using WENO schemes at time $T=0.16$ with 100 mesh
points. The numerical densities along with zoomed views near the discontinuities are
depicted in Figure \ref{Figure.Sod}. The simulations of WENO-AO(5,4,3) yield better resolution of solution near the discontinuities than other WENO 
 variants and it gives least dissipative results without over and undershoot at the shock position. Form Table \ref{Table.sod}, it can observed that 
 the $L^1$-error of WENO-AO(5,4,3) is least among the considered WENO
 scheme.  The solution computed with WENO-AO(5,3) and WENO-AON(5,3) has comparable accuracy
  and resolution at discontinuities.
  
   In WENO-AO(5,4,3), we have given twice the linear weight to the stencil $\mathbb{S}_0^3$ in comparison to stencils 
  $\mathbb{S}_{-1}^3$, $\mathbb{S}_1^3$. Now we construct the scheme based on Theorem \ref{thm:543w}, where equal weights are given to stencils $\mathbb{S}_{0}^3$, $\mathbb{S}_1^3$. The linear weights are
  chosen as follows
    \begin{equation}
 \left. 
\begin{array}{ll}
\gamma_{-1}^3&=\frac{1}{3}(1-\gamma_{Hi})(1-\gamma_{Lo}),\\
 \gamma_0^3&=\frac{1}{3}(1-\gamma_{Hi})(1-\gamma_{Lo}),\\
 \gamma_1^3&=\frac{1}{3}(1-\gamma_{Hi})(1-\gamma_{Lo}),\\
 \gamma_0^4&=(1-\gamma_{Hi})\gamma_{Lo},\\
 \gamma_0^5&=\gamma_{Hi},
\end{array}
\right \}
\end{equation}
 and corresponding scheme is named as WENO-AON(5,4,3). The values of $\gamma_{Hi}$ and $\gamma_{Lo}$ are taken to be 0.85. 
  In Figure \ref{Figure.Sod_comp}, we have compared the resolution of WENO-AO(5,4,3)  with that of WENO-AON(5,4,3) and WENO-AOL(5,4,3) scheme.  
  From the enlarged portion of Figure \ref{Figure.Sod_comp}, we observed that resolution of WENO-AO(5,4,3) and WENO-AOL(5,4,3) schemes are slightly better than WENO-AON(5,4,3) scheme across the shock and discontinuity. Whereas WENO-AO(5,4,3) and WENO-AOL(5,4,3) schemes are comparable to each other.
 }
\end{example}
\begin{table}[t!]
\centering 
\small
\begin{tabular}{|c| c|c| c| c| c|  c| c| r|}
\hline 
 $N$ & WENO-JS & WENO-AO(5,3) & WENO-AON(5,3) & WENO-AO(5,4,3) & WENO-ZQ & WENO-Z \\
\hline
$200$  &  3.5686e-03 & 2.9433e-03     &2.8900e-03 &2.8172e-03 &2.9151e-03 &3.2170e-03\\ 
\hline 
$400$  &  1.8130e-03 & 1.4768e-03     &1.4541e-03 &1.4180e-03  &1.4835e-03 & 1.6194e-03  \\
\hline
$800$  &  9.7134e-04 & 7.9350e-04     &7.8250e-04 &7.6496e-04  &8.0507e-04 &8.6793e-04    \\
\hline
\end{tabular}
\caption{ Comparison of $L^1$-error of WENO schemes for Example \ref{shock_test} at time $T=0.16$ over the domain $[0,1] $.}
\label{Table.sod}
\end{table}

\begin{table}[t!]
% \centering 
\small
\begin{tabular}{|c| c|c| c| c| c|  c| c| r|}
\hline 
 $N$ & WENO-JS & WENO-AO(5,3) & WENO-AON(5,3) & WENO-AO(5,4,3) & WENO-ZQ & WENO-Z \\
\hline
$200$  &  1.0773e-01 & 8.7228e-02    &8.6492e-02 &8.3750e-02 &8.8088e-02 &9.7515e-02\\ 
\hline 
$400$  &  5.2252e-02 & 4.0127e-02     &3.9965e-02 &3.8542e-02  &4.1212e-02 &4.5822e-02  \\
\hline
$800$  &  2.9815e-02 & 2.3262e-02     &2.3119e-02 &2.2765e-02  &2.4506e-02 &2.6248e-02     \\
\hline
\end{tabular}
\caption{ Comparison of $L^1$-error of WENO schemes for Example \ref{laxp} at time $T=1.3$ over the domain $[-4,4] $.}
\label{table.lax}
\end{table}


\begin{figure}[t!]
\begin{center}
 \begin{tabular}{cc}
\includegraphics[width=0.48\textwidth]{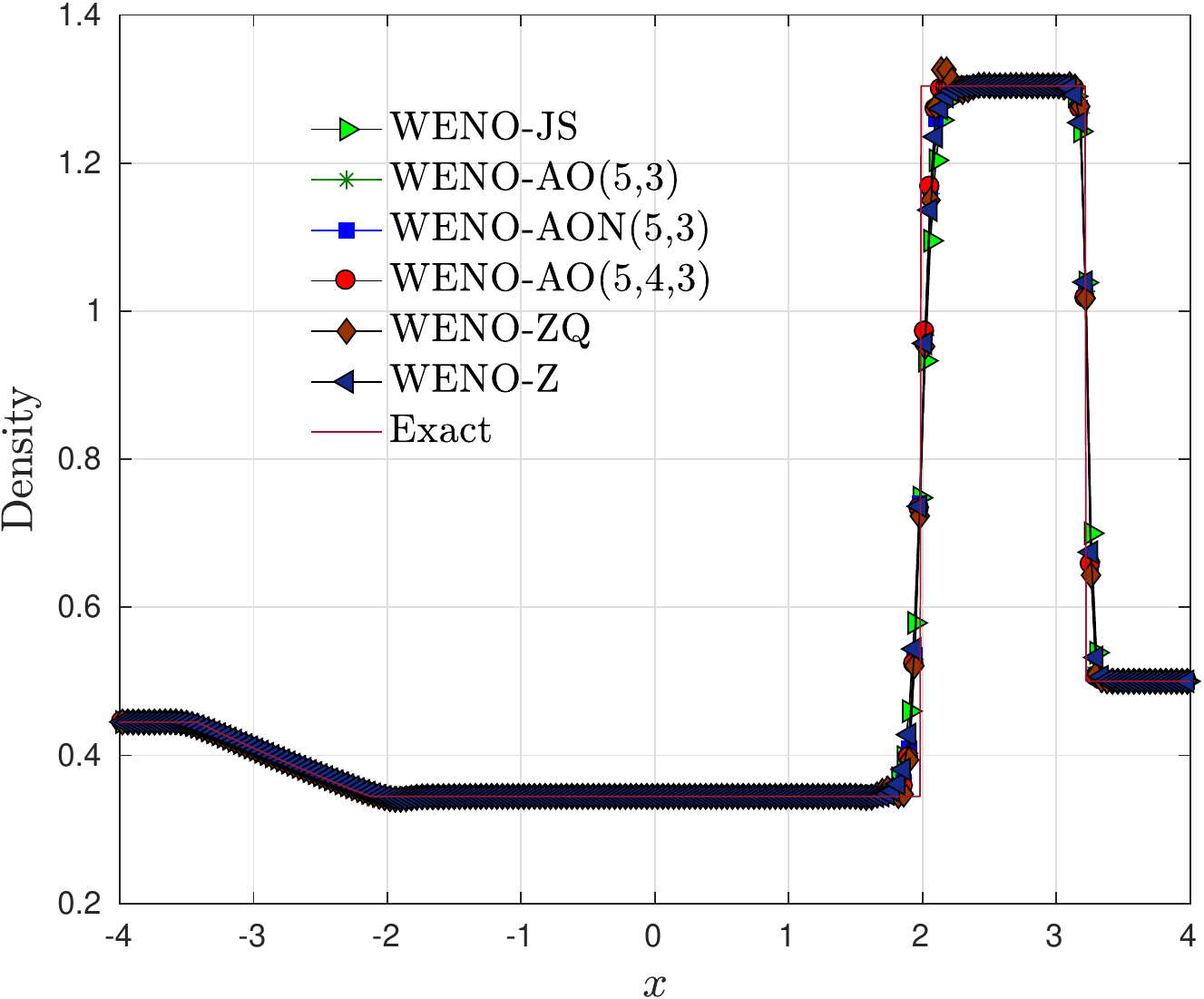}
\includegraphics[width=0.48\textwidth]{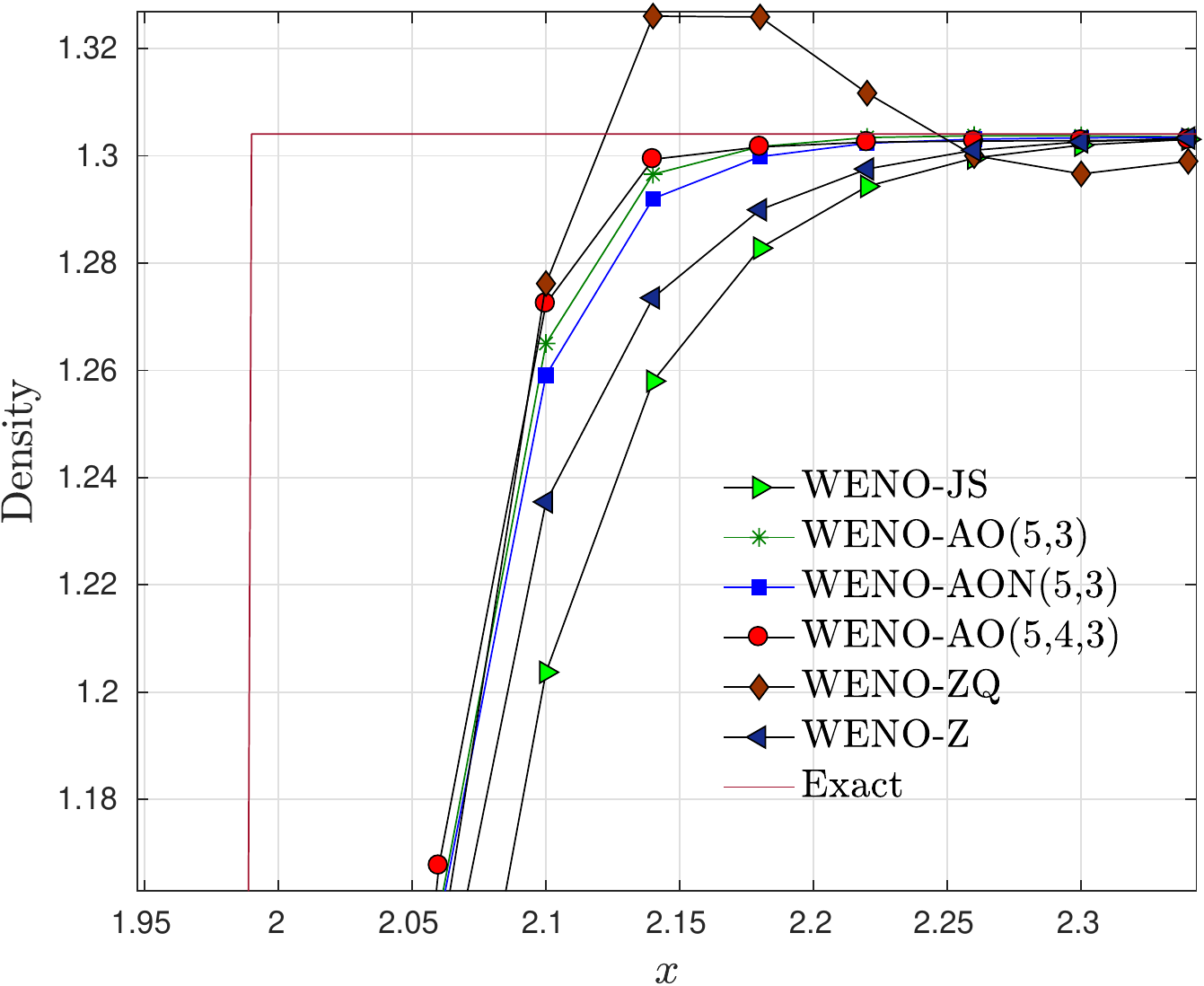}\\
\includegraphics[width=0.48\textwidth]{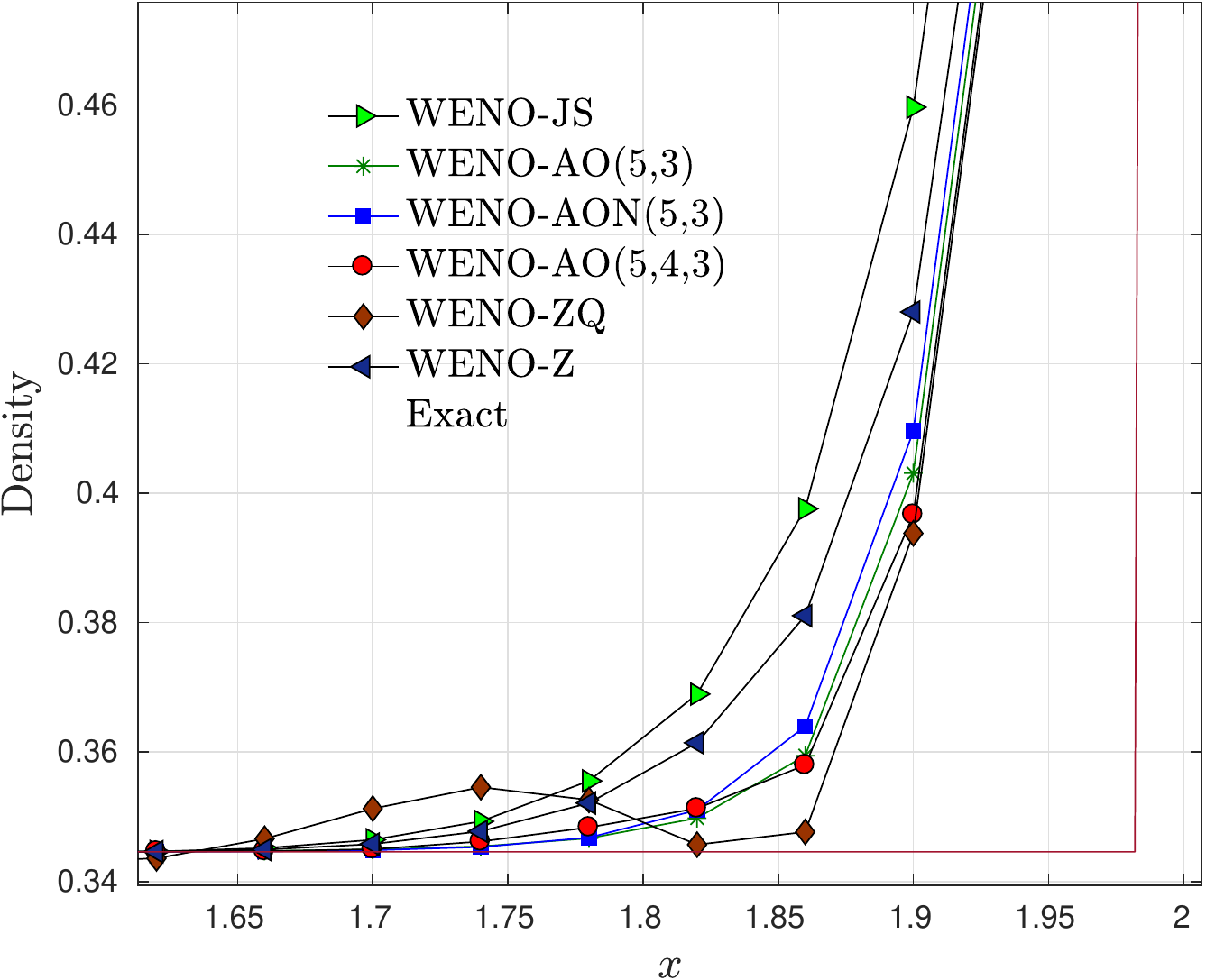}
\includegraphics[width=0.48\textwidth]{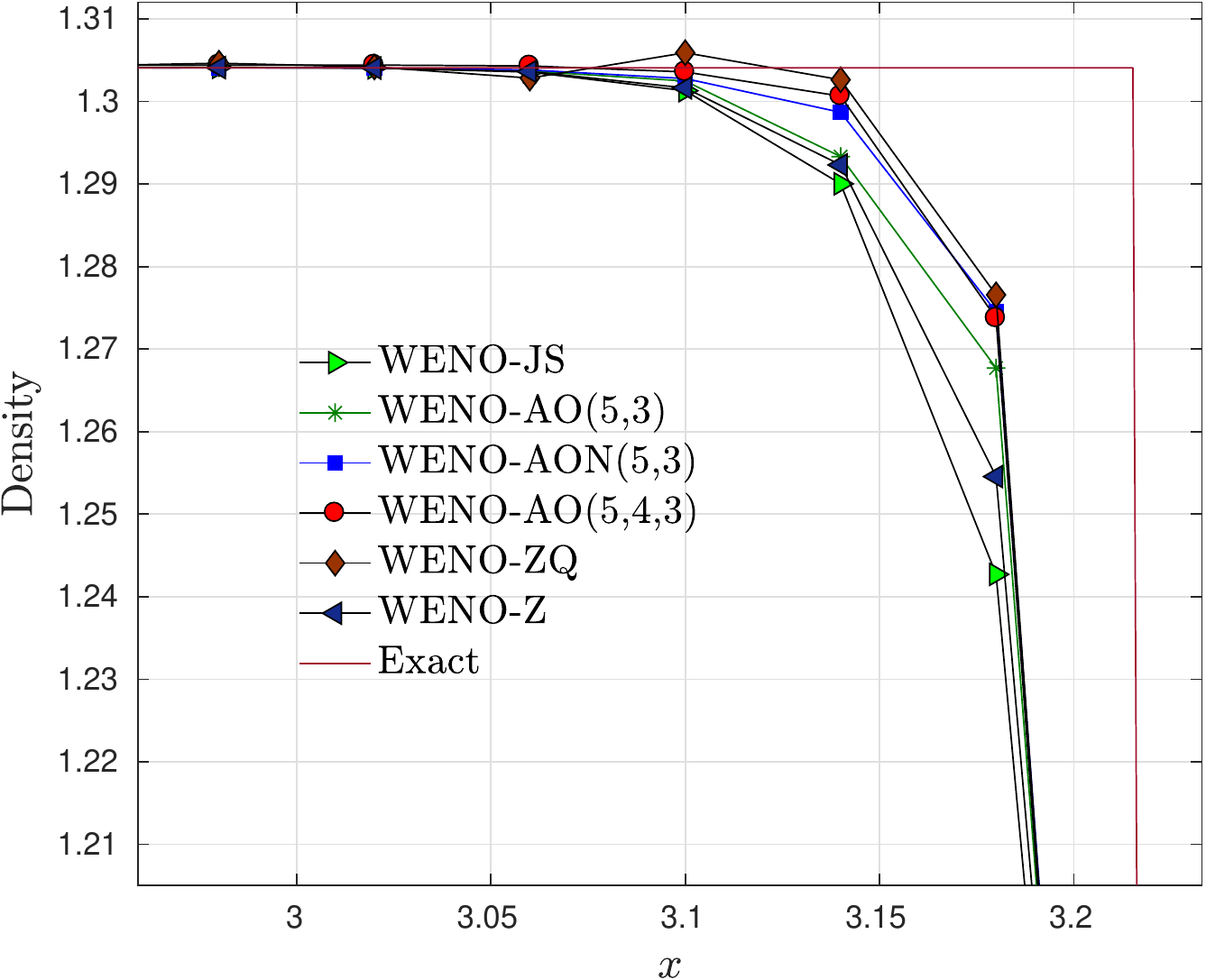}
 \end{tabular}
\end{center}
 \caption{ Comparison of the WENO-JS, WENO-Z, WENO-ZQ, WENO-AO(5,3), WENO-AON(5,3), and WENO-AO(5,4,3) schemes for the Example \ref{laxp}
 at time $T=1.3$ over a domain $[-4,4]$ with a uniform mesh of having 200 points.}
 \label{lax_fig}
\end{figure}

\begin{figure}[t!]
\begin{center}
 \begin{tabular}{cc}
\includegraphics[width=0.48\textwidth]{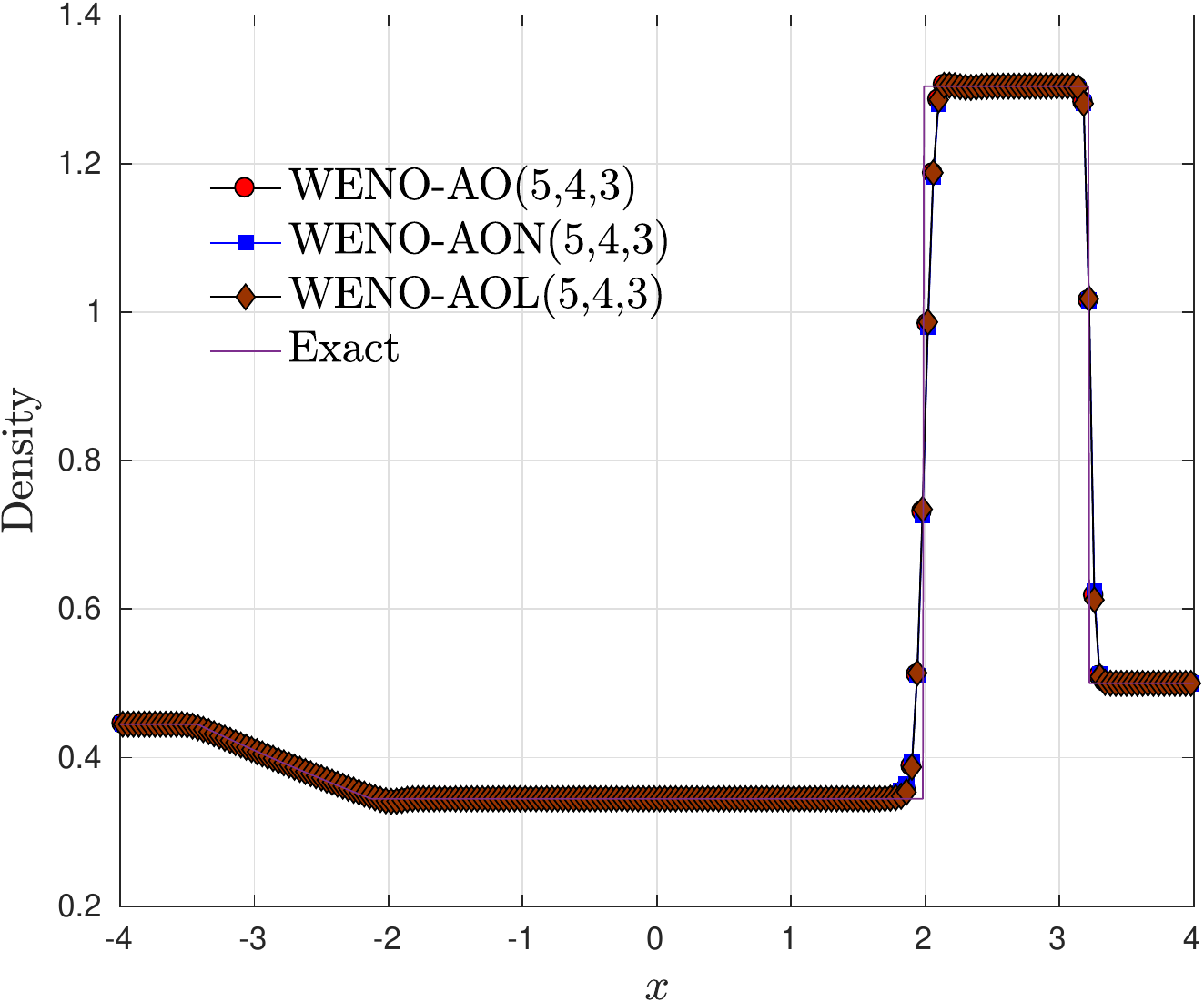}
\includegraphics[width=0.48\textwidth]{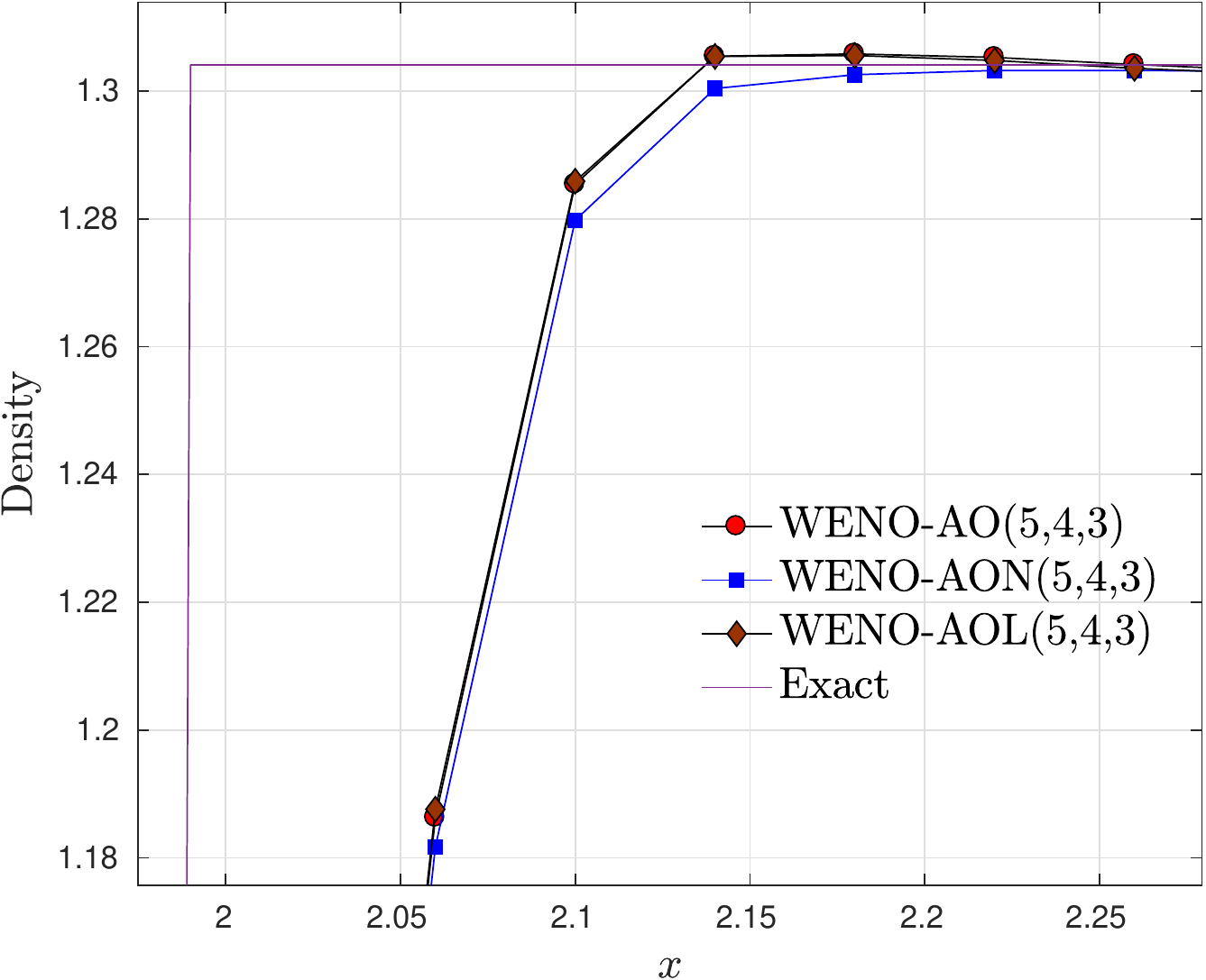}\\
\includegraphics[width=0.48\textwidth]{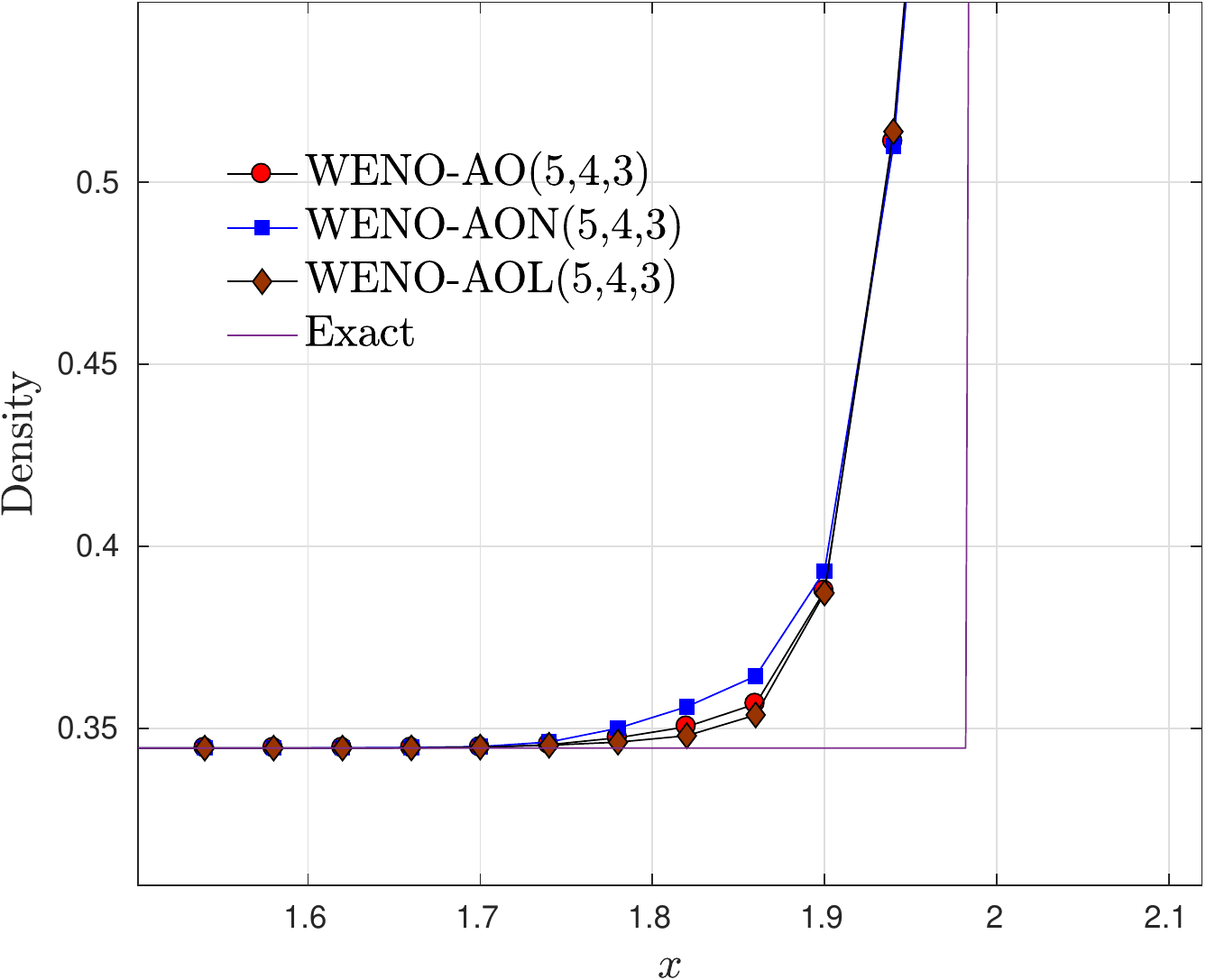}
\includegraphics[width=0.48\textwidth]{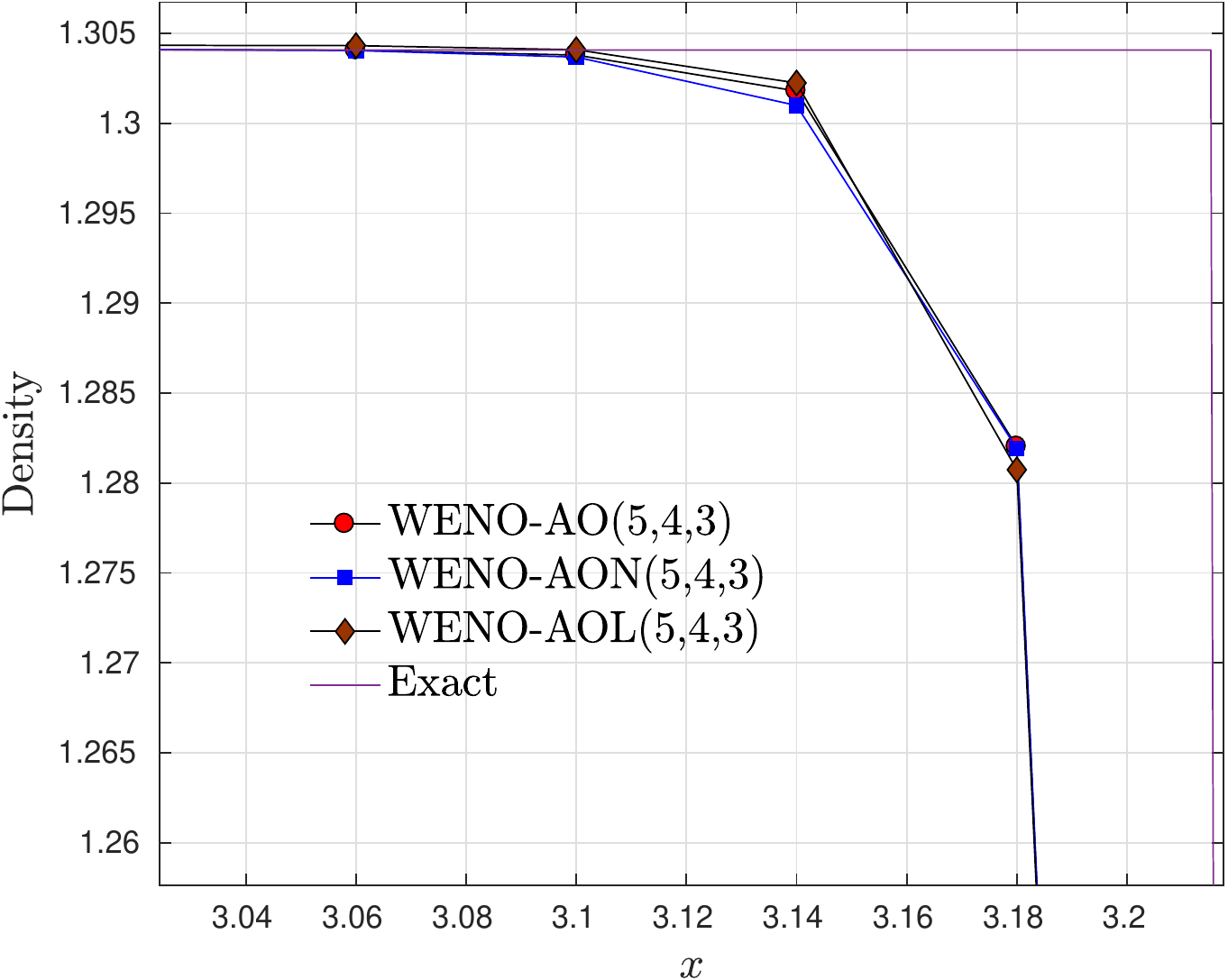}
 \end{tabular}
\end{center}
 \caption{ Comparison of the WENO-AO(5,4,3), WENO-AOL(5,4,3), and WENO-AON(5,4,3) schemes for the Example \ref{laxp}
 at time $T=1.3$ over a domain $[-4,4]$ with a uniform mesh of having 200 points.}
 \label{lax_comp_fig}
\end{figure}
\begin{example}\label{laxp}{\rm (Lax Problem)
Consider the Lax problem for which initial data is given by
   \begin{equation}\label{lax.system}
(\rho, u, p)(x,0)=\left\{\begin{array}{ll}
                (0.445,0.6980,3.528) ~~~~   x<0, \\
                (0.5,0.0,0.571) ~ ~~~~~~~~~~~ x>0.
                \end{array}\right.
\end{equation}}
\end{example}
The numerical solutions are computed at time $T=1.3$ over the domain $[-4,4]$ using 200 mesh points. We examine the numerical solutions obtained 
with different WENO schemes and compare it with the present WENO algorithms. In Figure \ref{lax_fig}, we have shown the comparison of density computed 
with WENO schemes
 along with the exact solution. It can be easily observed from Figure \ref{lax_fig} that WENO-AO(5,4,3) scheme resolves the shock and contact discontinuities 
 accurately without over- and under-shoot, whereas over- and under-shoot can be observed in case of WENO-ZQ scheme (see \cite[pp. 118, Figure 3.5]{zhu-qiu_16a}). The resolutions of WENO-AO(5,3) and
  WENO-AON(5,3) are comparable near the discontinuities. In Table \ref{table.lax}, we compared the $L^1$-errors of the WENO schemes increasing the
  number of mesh points for fixed time  $T=1.3$. The WENO-AO(5,4,3) scheme has least $L^1$-error among the considered WENO schemes.
  
In Figure \ref{lax_comp_fig}, we have compared the WENO-AO(5,4,3) scheme with WENO-AON(5,4,3) and WENO-AOL(5,4,3) scheme. In this case also WENO-AO(5,4,3) performs slightly better than WENO-AON(5,4,3) scheme and is comparable with WENO-AOL(5,4,3) scheme.

\begin{figure}[t!]
\begin{center}
 \begin{tabular}{cc}
\includegraphics[width=0.48\textwidth]{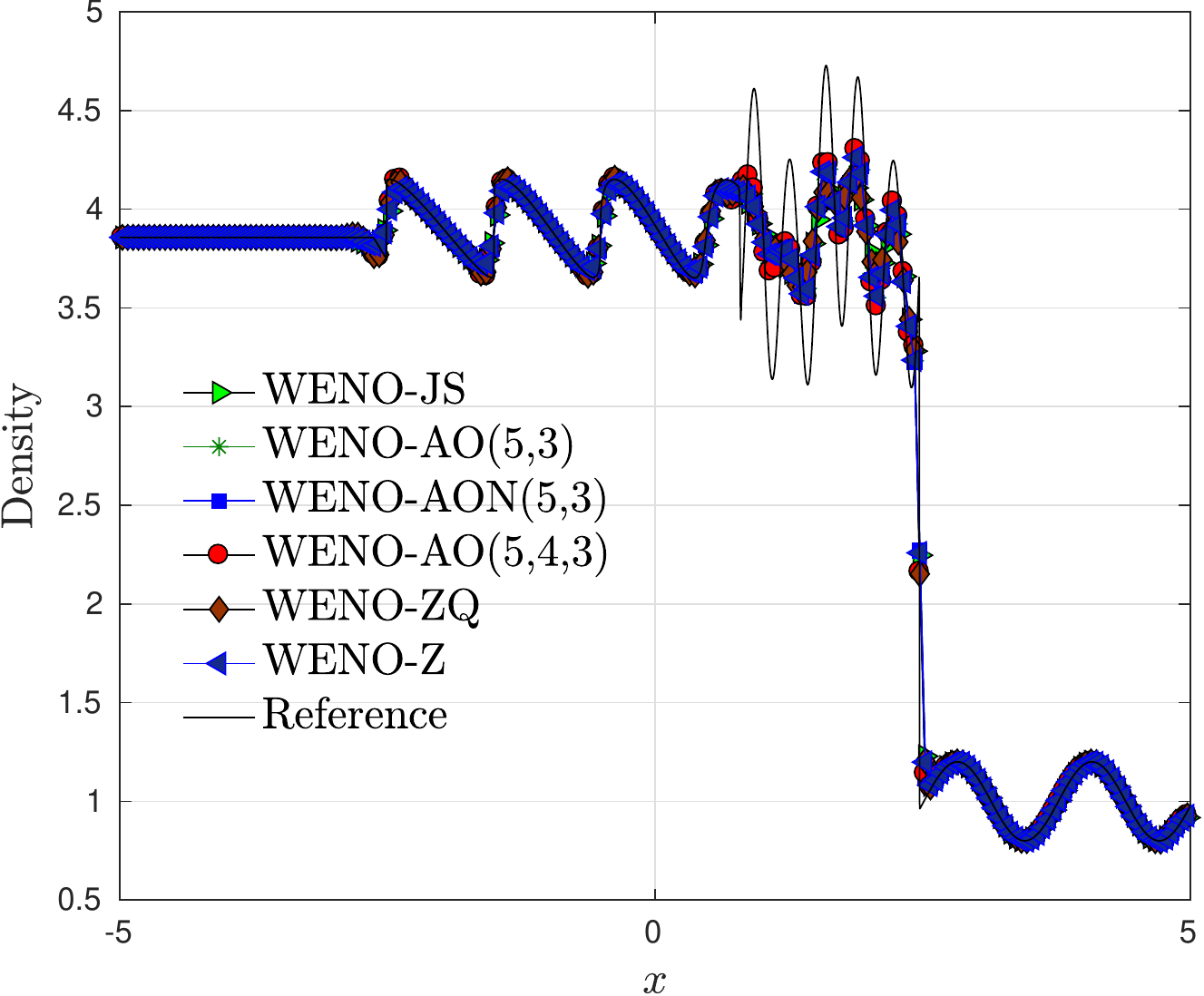}
\includegraphics[width=0.48\textwidth]{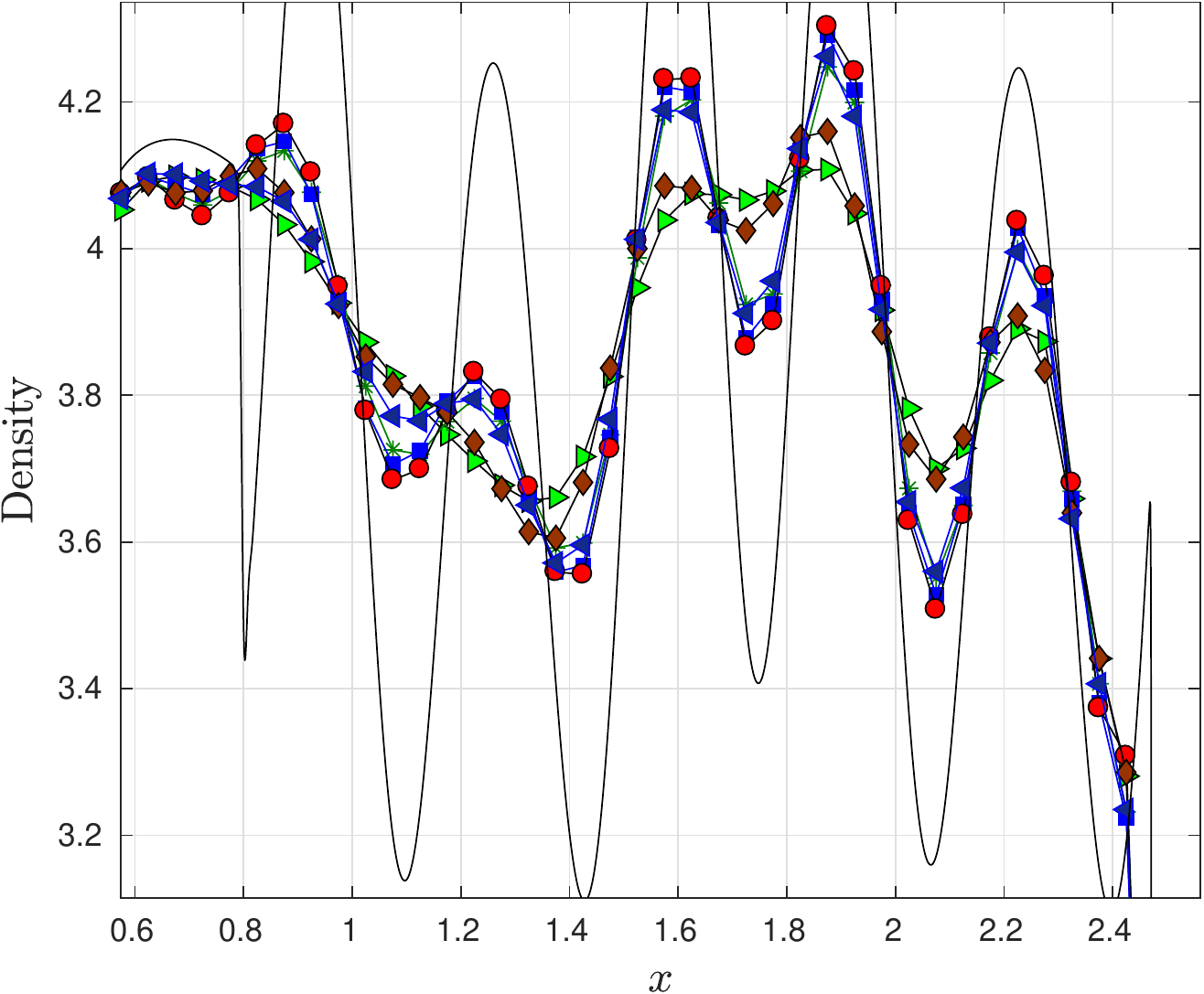}\\
\includegraphics[width=0.48\textwidth]{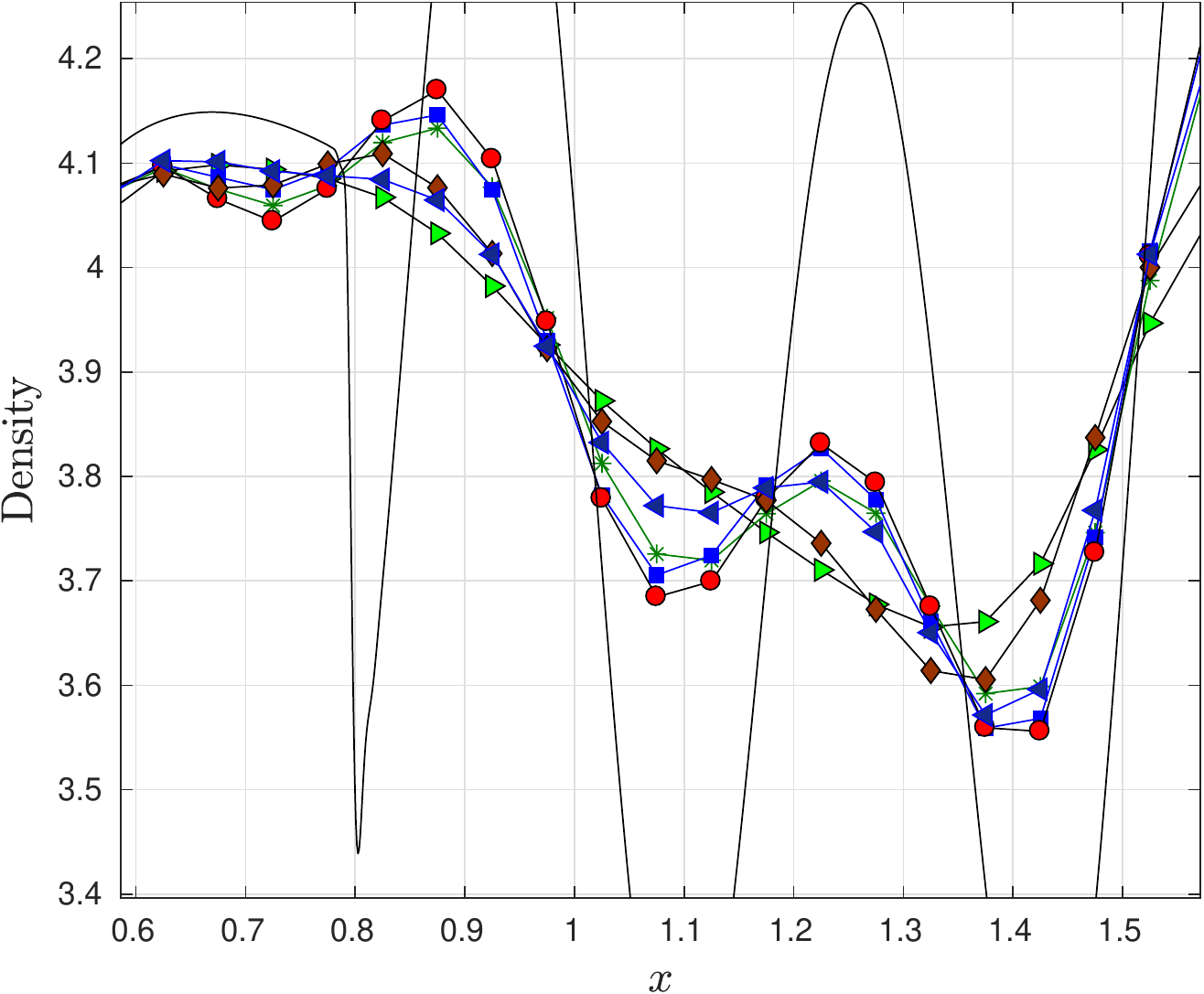}
\includegraphics[width=0.48\textwidth]{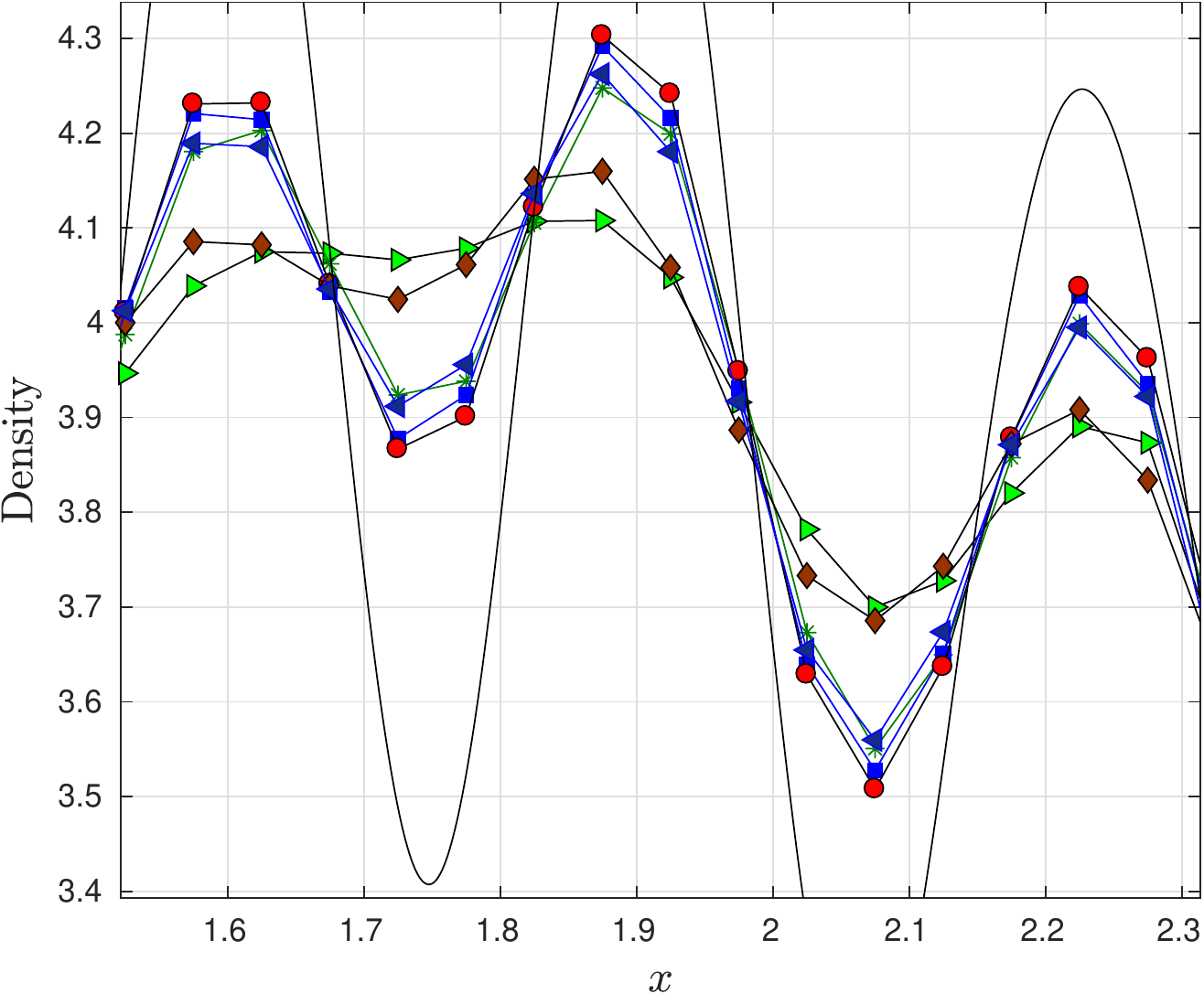}
 \end{tabular}
\end{center}
 \caption{Comparison of the WENO-JS, WENO-Z, WENO-ZQ, WENO-AO(5,3), WENO-AON(5,3), and WENO-AO(5,4,3) schemes for the Example \ref{shu-osh}
 at time $T=1.8$ over a uniform mesh of  resolution 200 points.}
 \label{shu_osh200}
\end{figure}

\begin{figure}[t!]
\begin{center}
 \begin{tabular}{cc}
\includegraphics[width=0.48\textwidth]{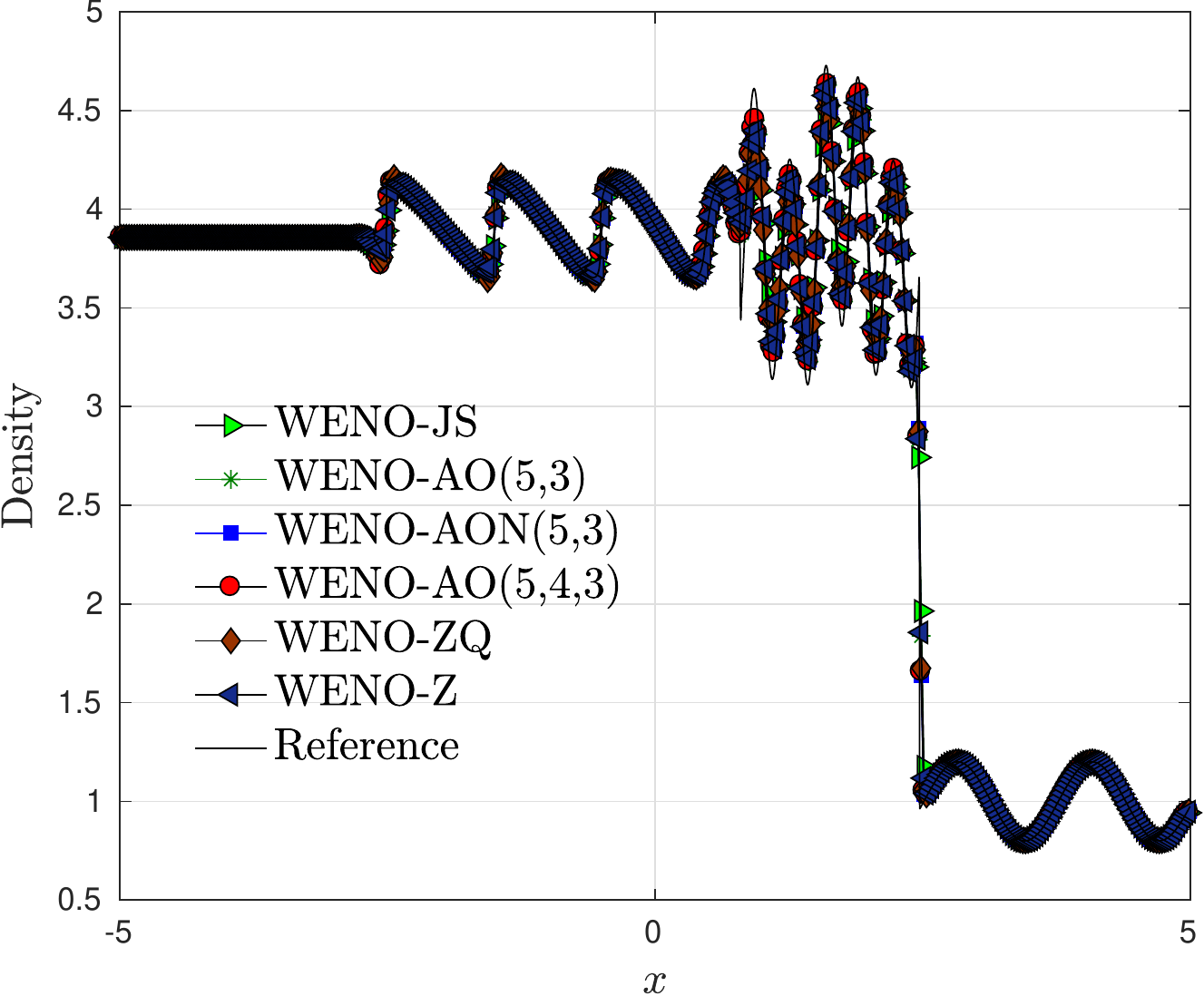}
\includegraphics[width=0.48\textwidth]{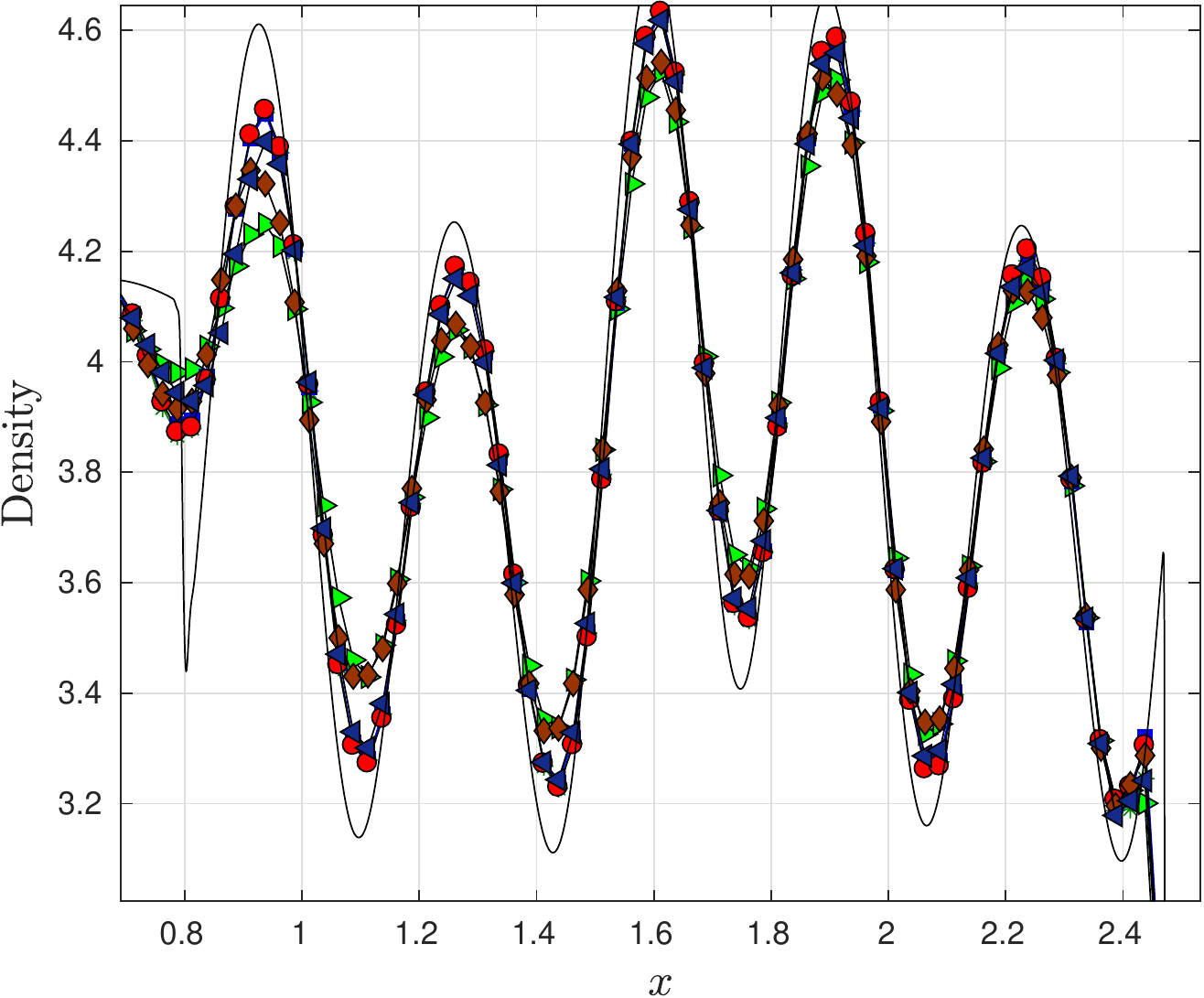}\\
\includegraphics[width=0.48\textwidth]{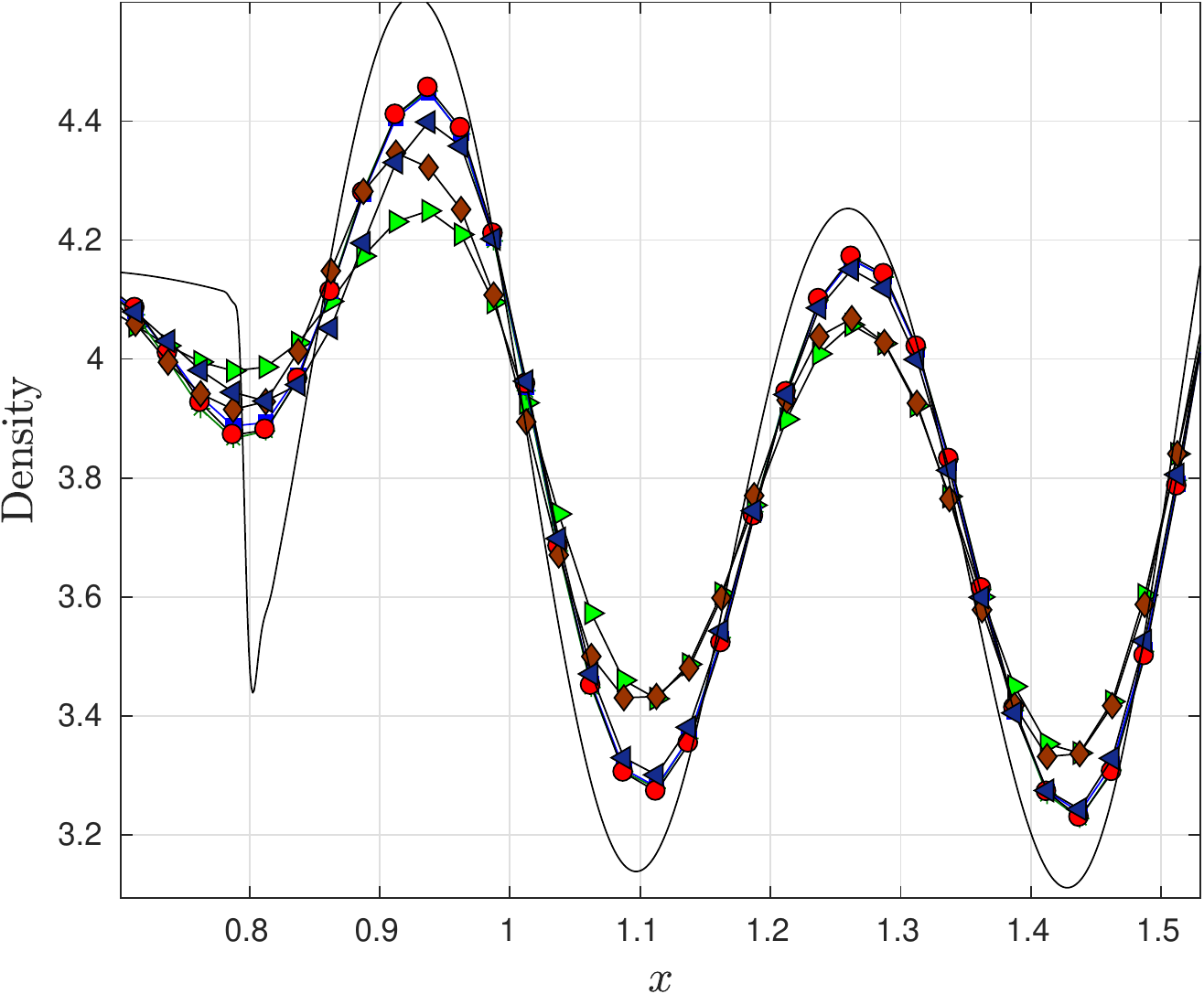}
\includegraphics[width=0.48\textwidth]{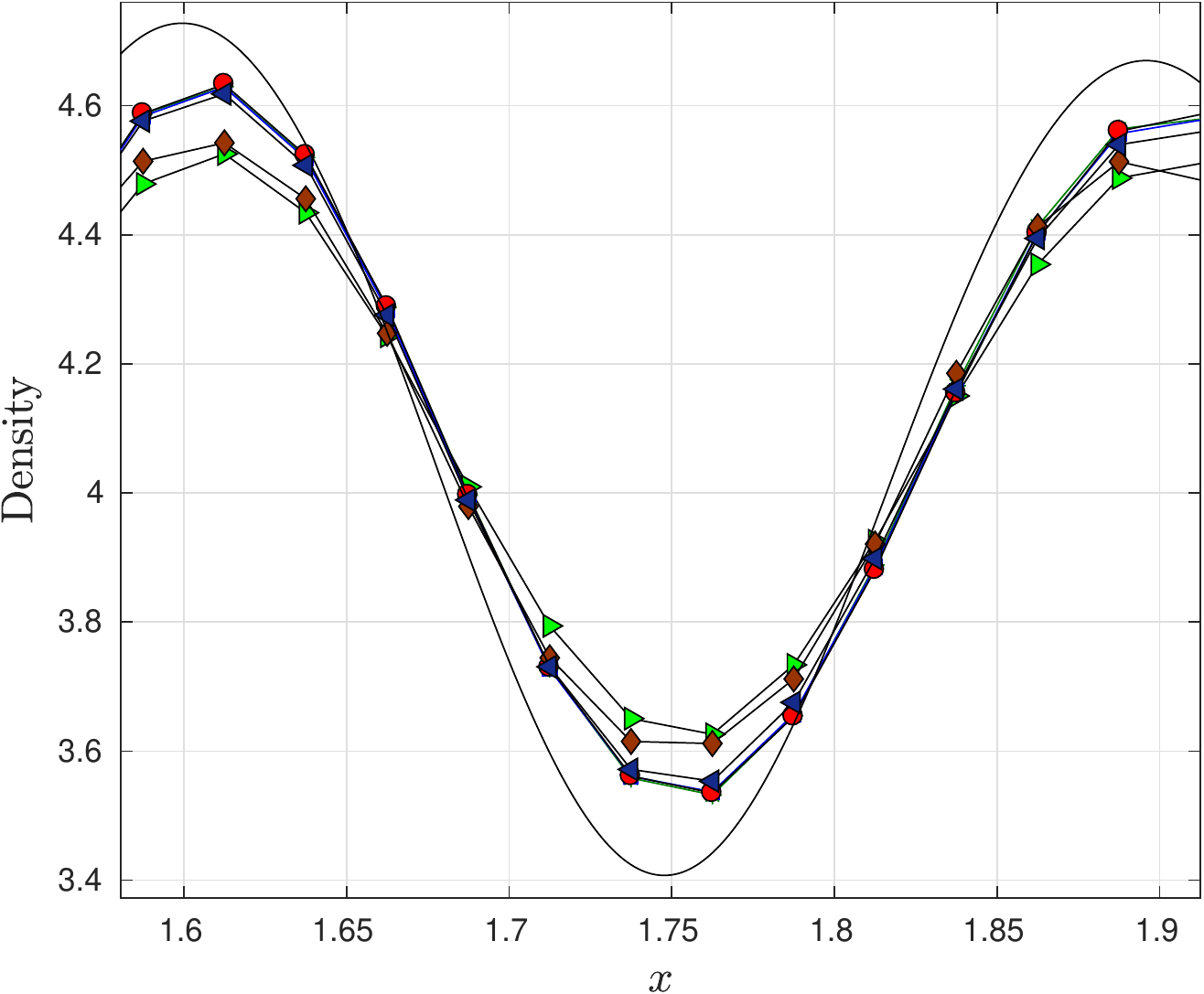}
 \end{tabular}
\end{center}
\caption{ Comparison of the WENO-JS, WENO-Z, WENO-ZQ, WENO-AO(5,3), WENO-AON(5,3), and WENO-AO(5,4,3) schemes for the Example \ref{shu-osh}
 at time $T=1.8$ over a uniform mesh of  resolution 400 points.}
 \label{shu_osh400}
\end{figure}

\begin{example}\label{shu-osh}{\rm (Shu-Osher Problem)
This test is used to simulate Shu and Osher problem  \cite{shu-osh_89a}, which involves the interaction of an entropy sine wave with a Mach 3
shock moving to the
 right. The simulations are performed over the domain $[-5,5]$ upto  time $T=1.8$. The initial conditions are
\begin{equation}\label{shu_osh}
(\rho, u, p)(x,0)=\left\{\begin{array}{ll}
                (3.857143,2.699369,10.33333) ~~~~   x<-4, \\
                (1+0.2\sin(5x),0.0,1.0) ~~~~~~~~~~~~  x>-4.
                \end{array}\right.
\end{equation}
with transmissive boundary conditions at both ends. Figure \ref{shu_osh200} (a) shows the density on a mesh  having 200 points for the WENO-Z, 
WENO-AO(5,3), WENO-AON(5,3), and WENO-AO(5,4,3) schemes, whereas in Figure \ref{shu_osh400}, we have shown the density computed with the considered schemes
on a mesh having
 400 points. The exact solution is a reference  solution computed using WENO-AO(5,3) scheme on a mesh with 10000 points. Figure
 \ref{shu_osh200} (b), (c), and (d) are the zoomed versions of the Figure \ref{shu_osh200} around the post-shock high frequency waves. It is observed that 
 the WENO-AO(5,4,3) scheme shows  significant improvement over the other
  considered schemes in resolving the high frequency wave in non-oscillatory fashion. From Figure \ref{shu_osh200} and \ref{shu_osh400}, we 
  observed that WENO-AON(5,3) scheme resolves shock more accurately than  
  WENO-AO(5,3) and 
   WENO-Z schemes on mesh with 200 points and is comparable with WENO-AO(5,3) scheme on a mesh with 400 points. In addition, the WENO-AO(5,4,3) shows 
   significantly lower smearing of the shock waves in comparison to other schemes.
  }
\end{example}

\begin{figure}[t!]
\begin{center}
 \begin{tabular}{cc}
\includegraphics[width=0.48\textwidth]{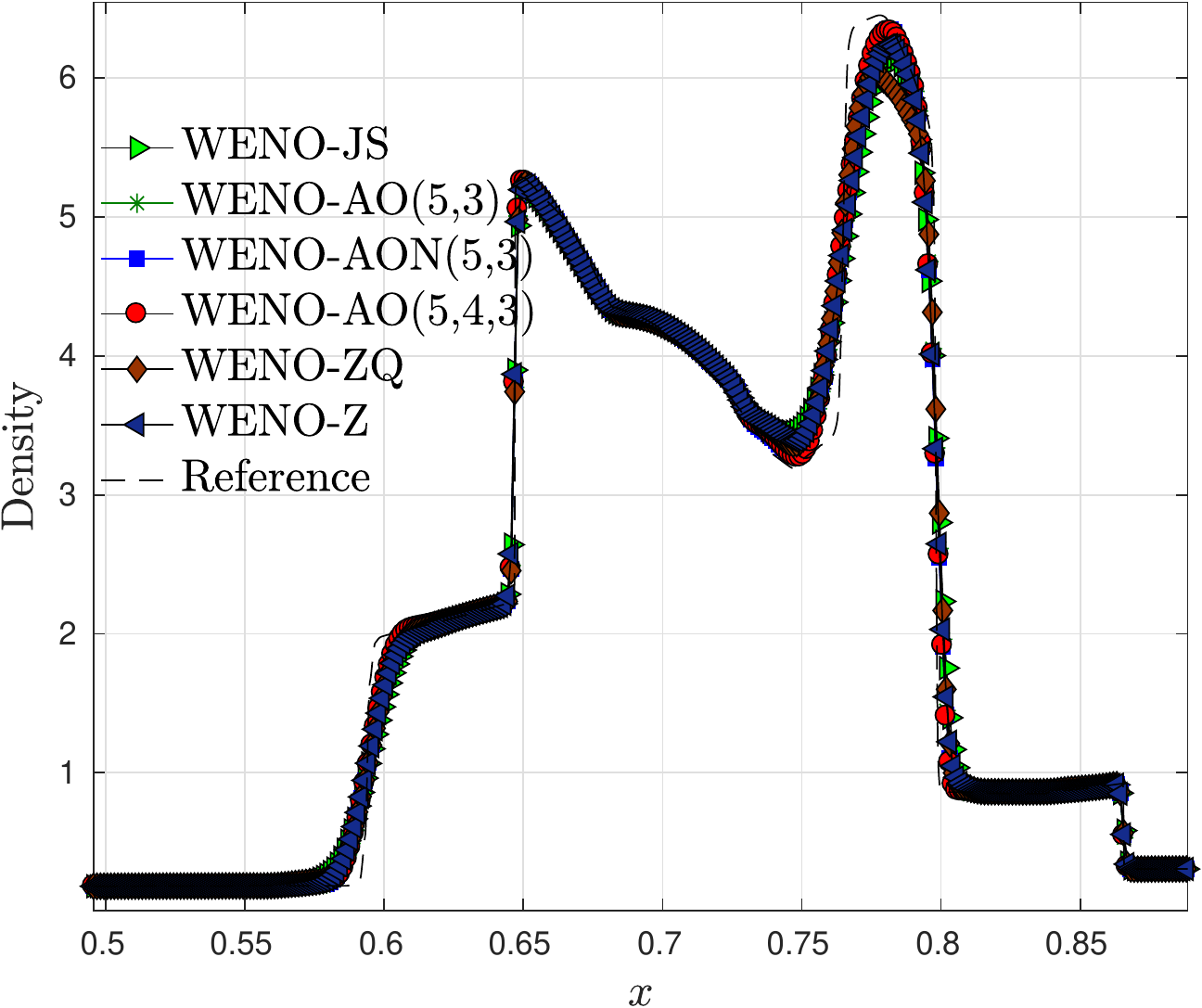}
\includegraphics[width=0.48\textwidth]{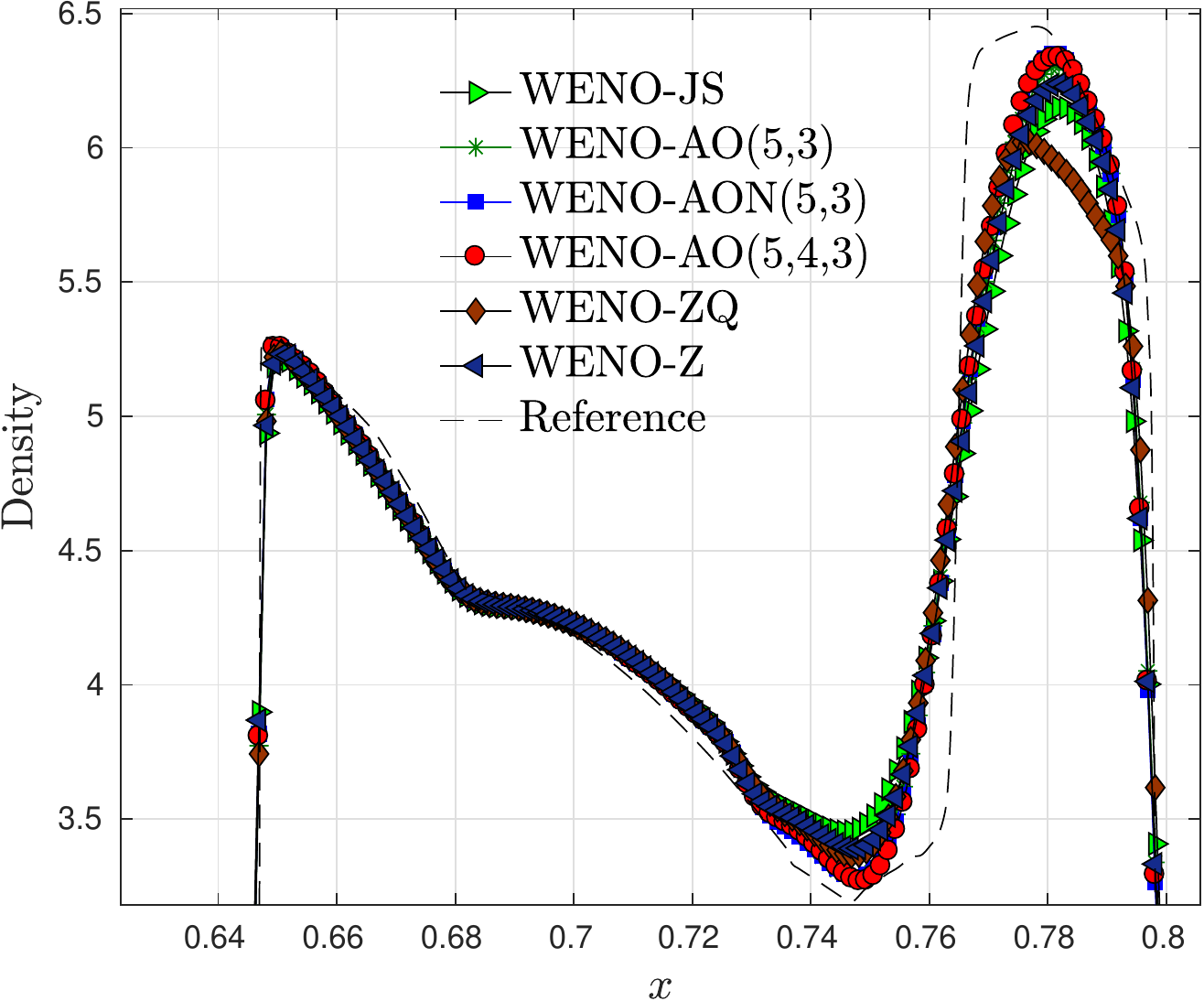}\\
\includegraphics[width=0.48\textwidth]{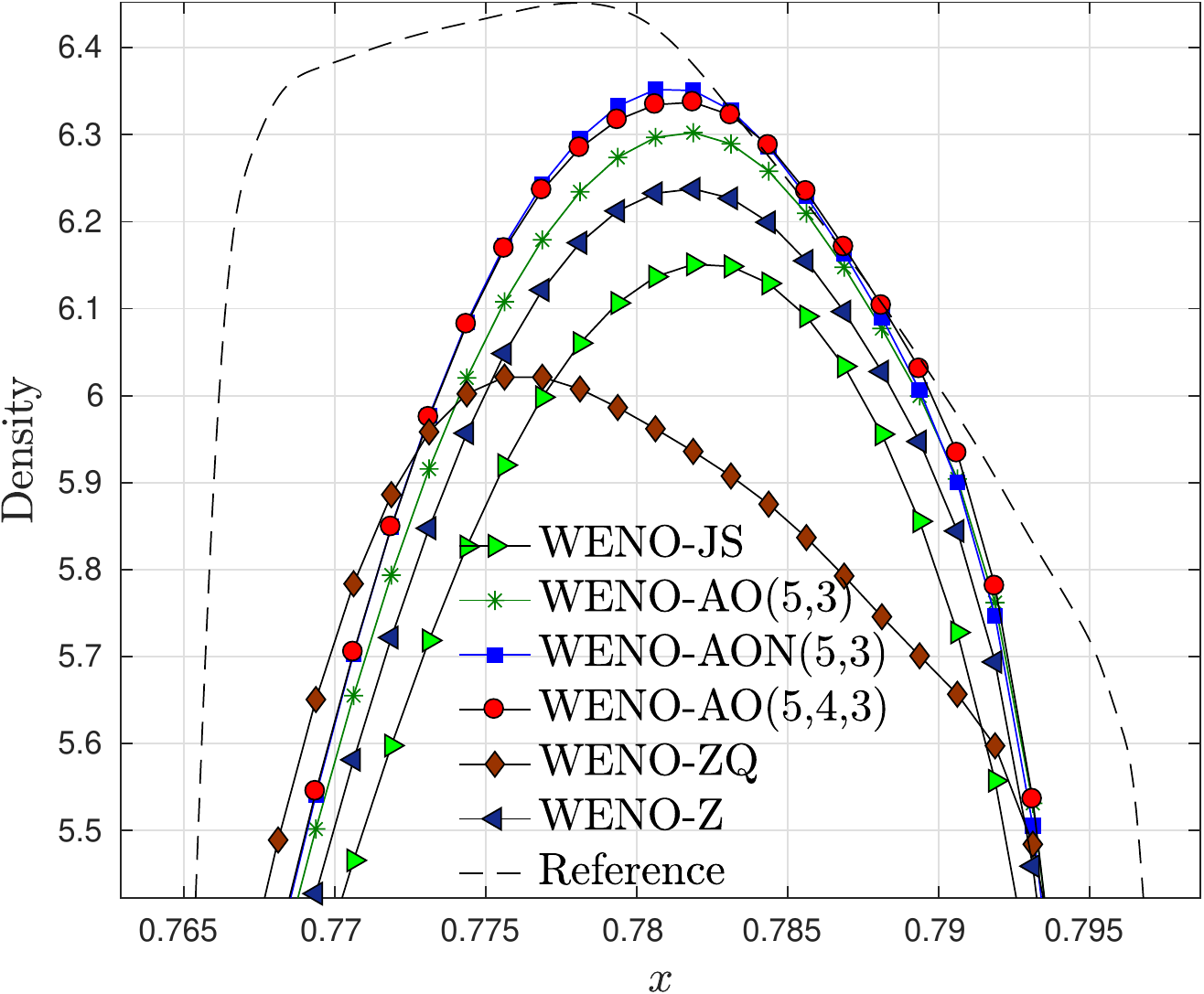}
\includegraphics[width=0.48\textwidth]{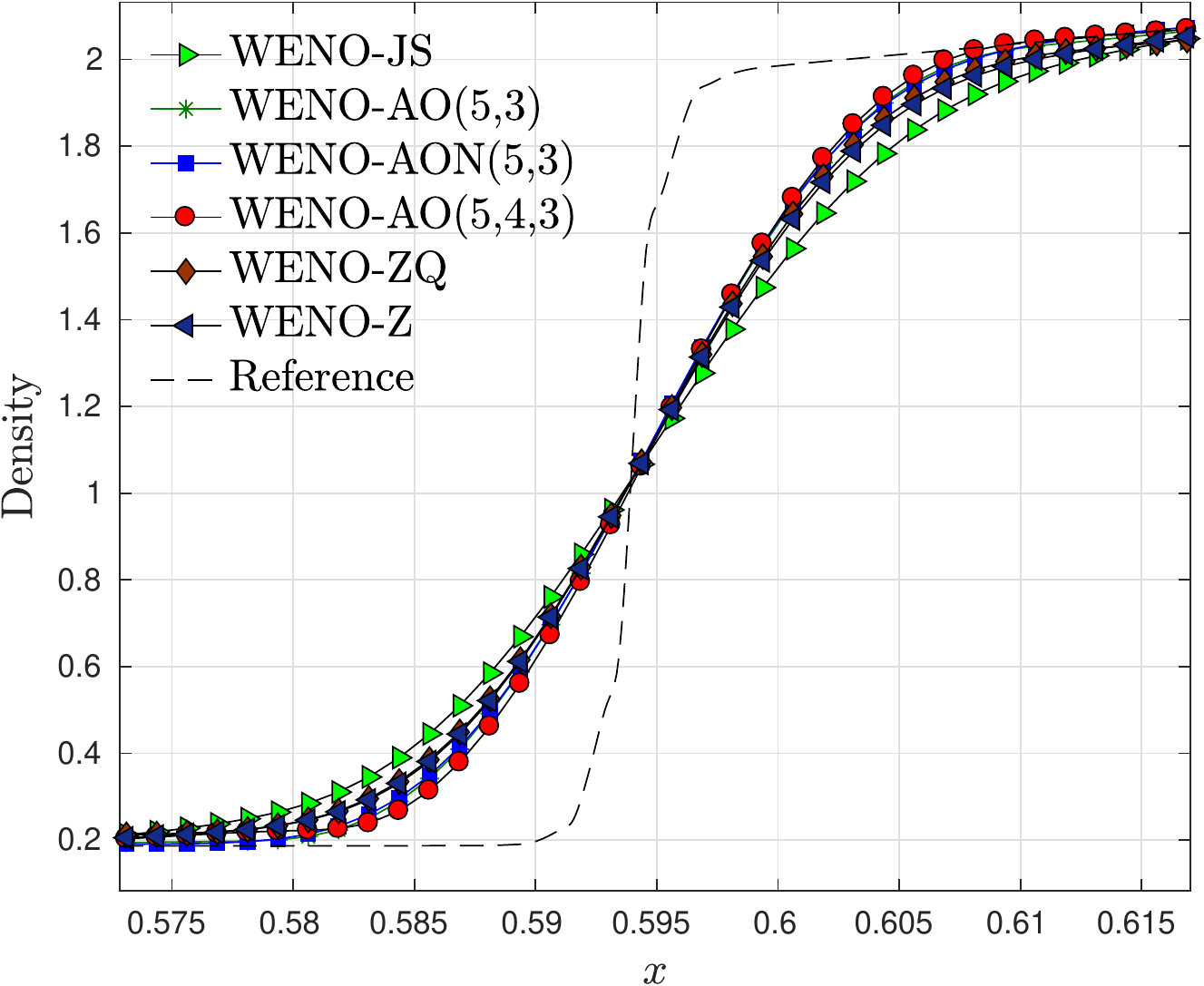}
 \end{tabular}
\end{center}
 \caption{Comparison of the WENO-JS, WENO-Z, WENO-ZQ, WENO-AO(5,3), WENO-AON(5,3), and WENO-AO(5,4,3) schemes for the Example \ref{blast}
 at time $T=0.038$ over a uniform mesh of  resolution 800 points.}
 \label{Blast.Figure}
\end{figure}

\begin{example}\label{blast}{\rm (Blast wave problem)
This one-dimensional test case involves the generation and interaction of two blast waves \cite{woo-col_84a}, with the 
following initial condition 
   \begin{equation}
(\rho, u, p)(x,0)=\left\{\begin{array}{ll}
                (1.0,0.0,1000) ~~~~   0.0<x<0.1, \\
                (1.0,0.0,0.01) ~ ~~~~0.1<x<0.9.\\
                (1.0,0.0,100) ~ ~~~~~ 0.9<x<1.0.
                \end{array}\right.
\end{equation}
and 
using reflective boundary conditions at both ends. The numerical solutions are computed at time $T=0.038$ over the domain $[0,1]$ with 800 mesh points
 and depicted in Figure \ref{Blast.Figure}. We have compared the numerical solutions
 with the reference solution, which we obtained using WENO5-JS scheme with 10000 grid points. It can be concluded from  careful examination 
 of Figure \ref{Blast.Figure} that resolution of the solution computed with WENO-AO(5,4,3) scheme is better than other WENO
  schemes.}
\end{example}

\begin{figure}
\begin{center}
 \begin{tabular}{cc}
 \includegraphics[width=0.48\textwidth]{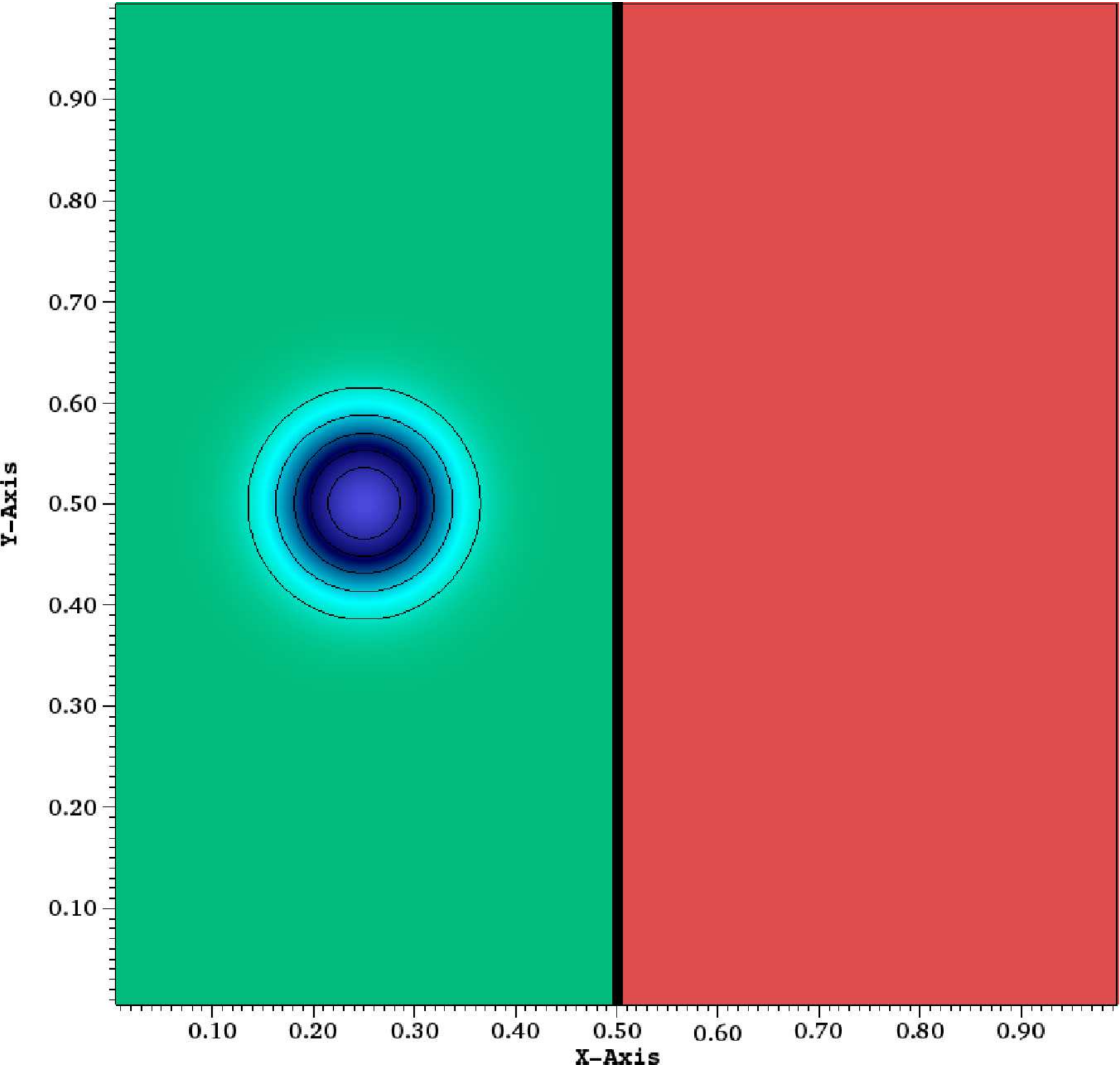} &
 \includegraphics[width=0.48\textwidth]{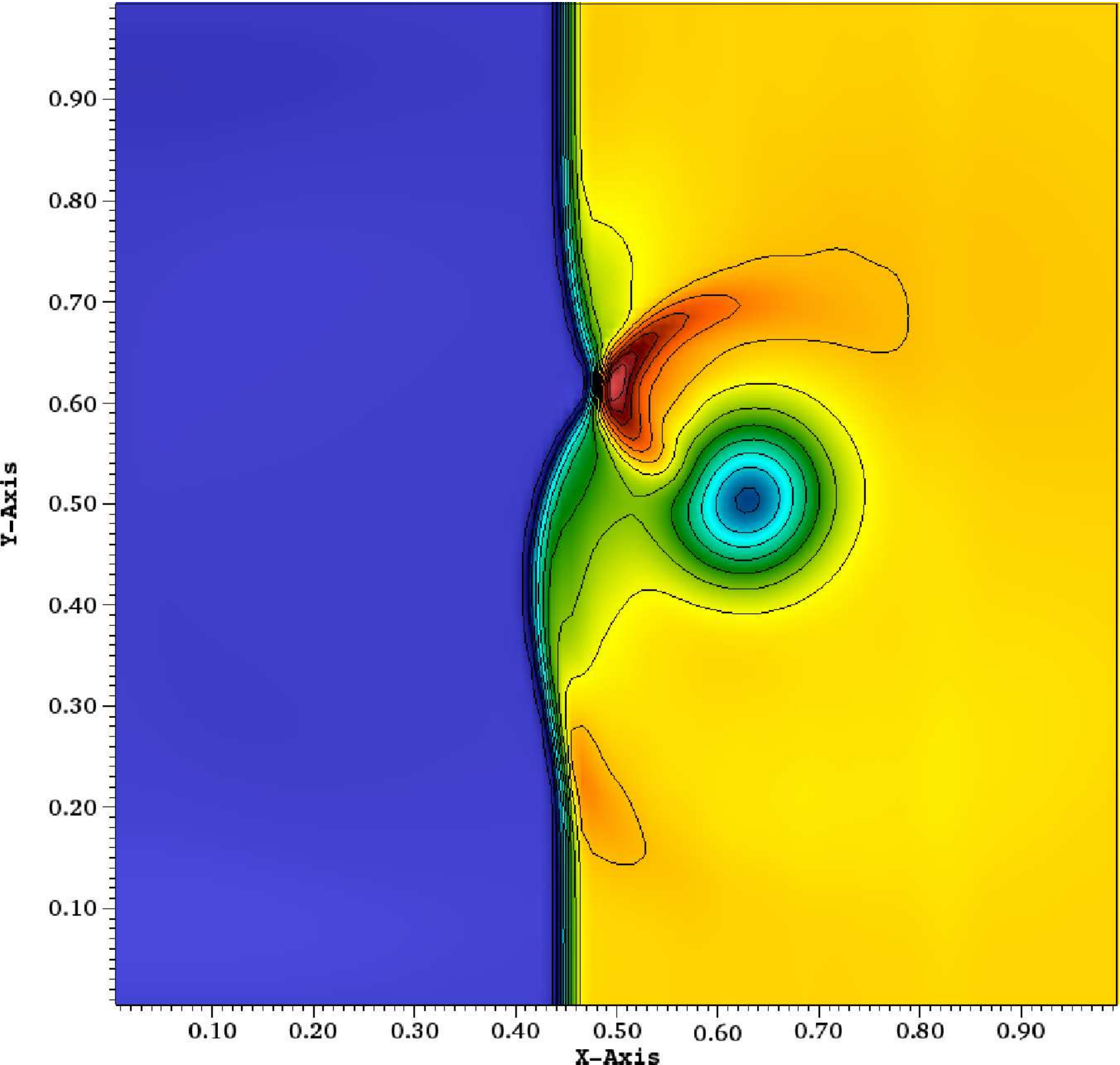}\\
 (a) Initial Density & (b) Final Density
 \end{tabular}
\end{center}
 \caption{Density contours of the Shock-vortex interaction problem  using 15 contours lines with range from 0.9 to 1.4, computed using WENO-AO(5,4,3) scheme at
 time $T=0.35$  with uniform mesh of resolution $100\times 100$.}
 \label{Figure.sv1}
\end{figure}

\begin{figure}[t!]
\begin{center}
 \begin{tabular}{cc}
\includegraphics[width=0.48\textwidth]{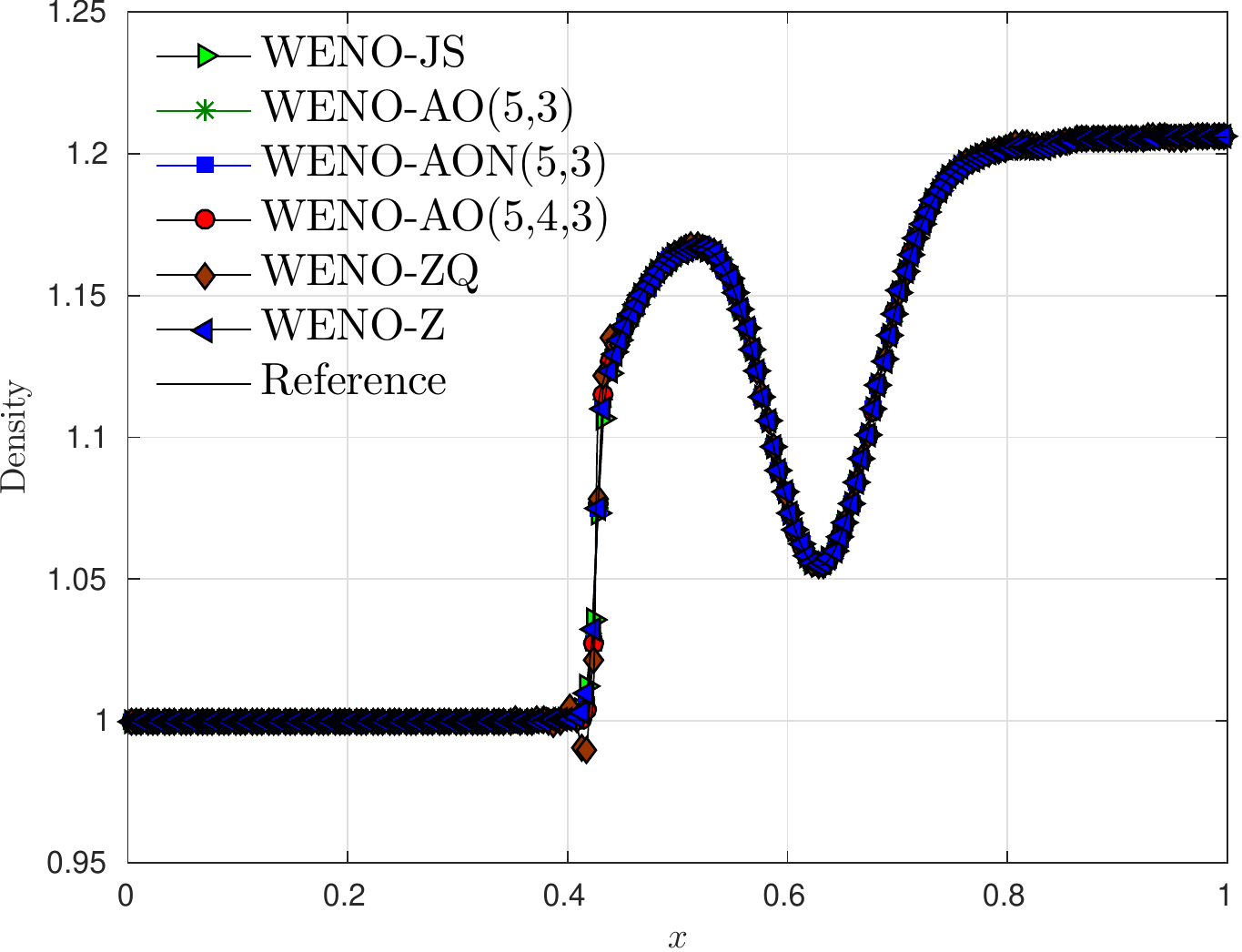}
\includegraphics[width=0.48\textwidth]{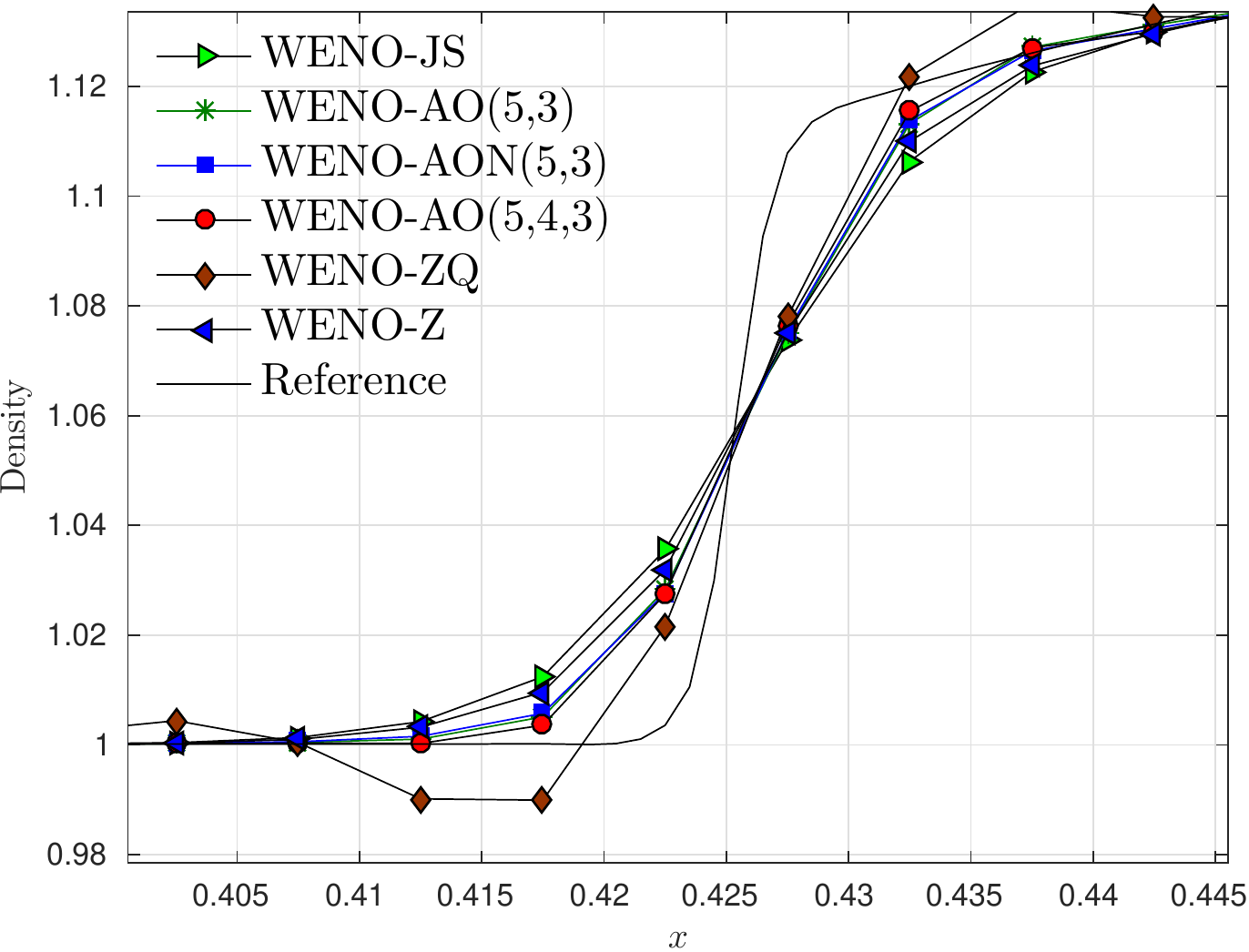}\\
\includegraphics[width=0.48\textwidth]{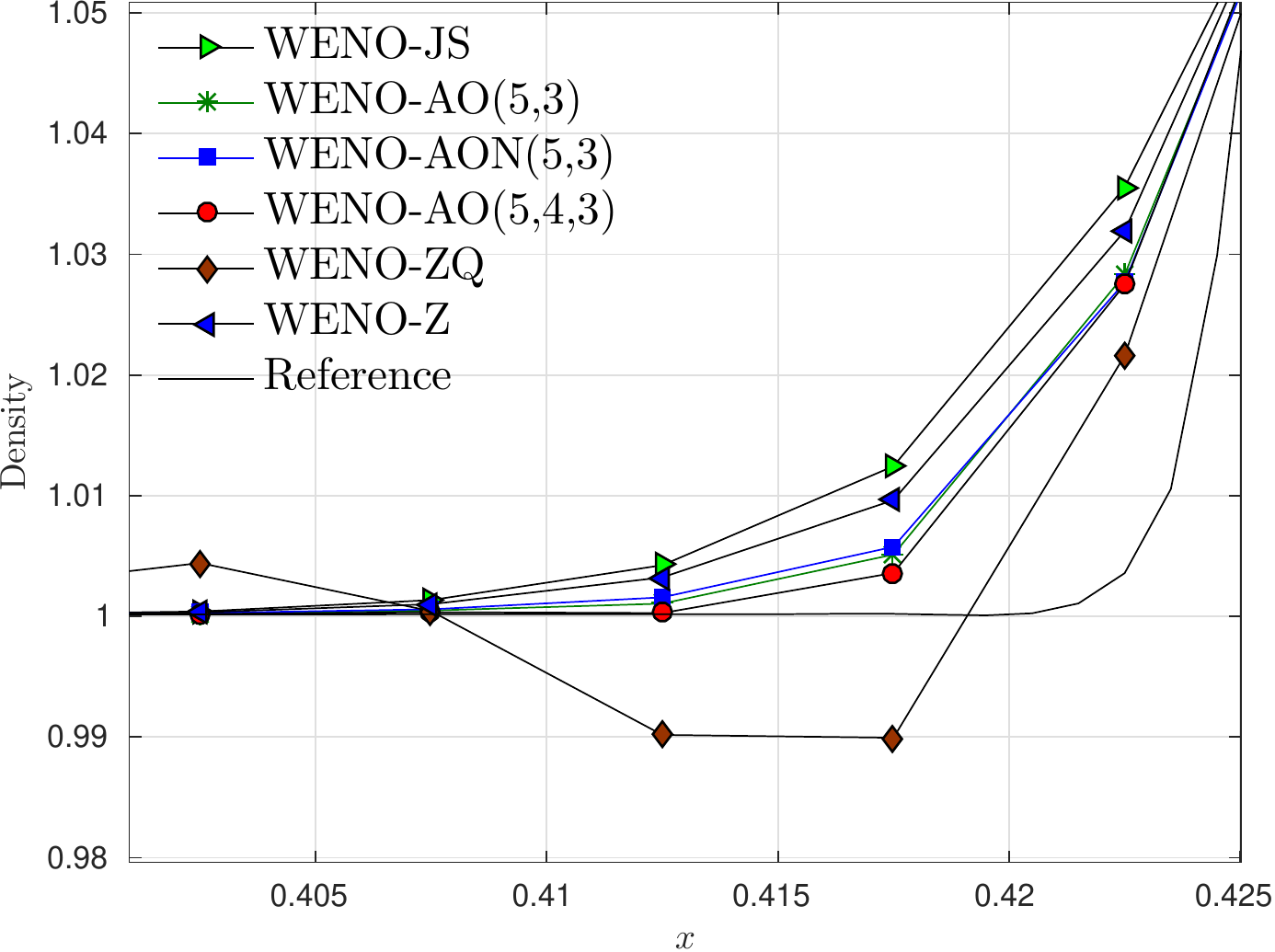}
\includegraphics[width=0.48\textwidth]{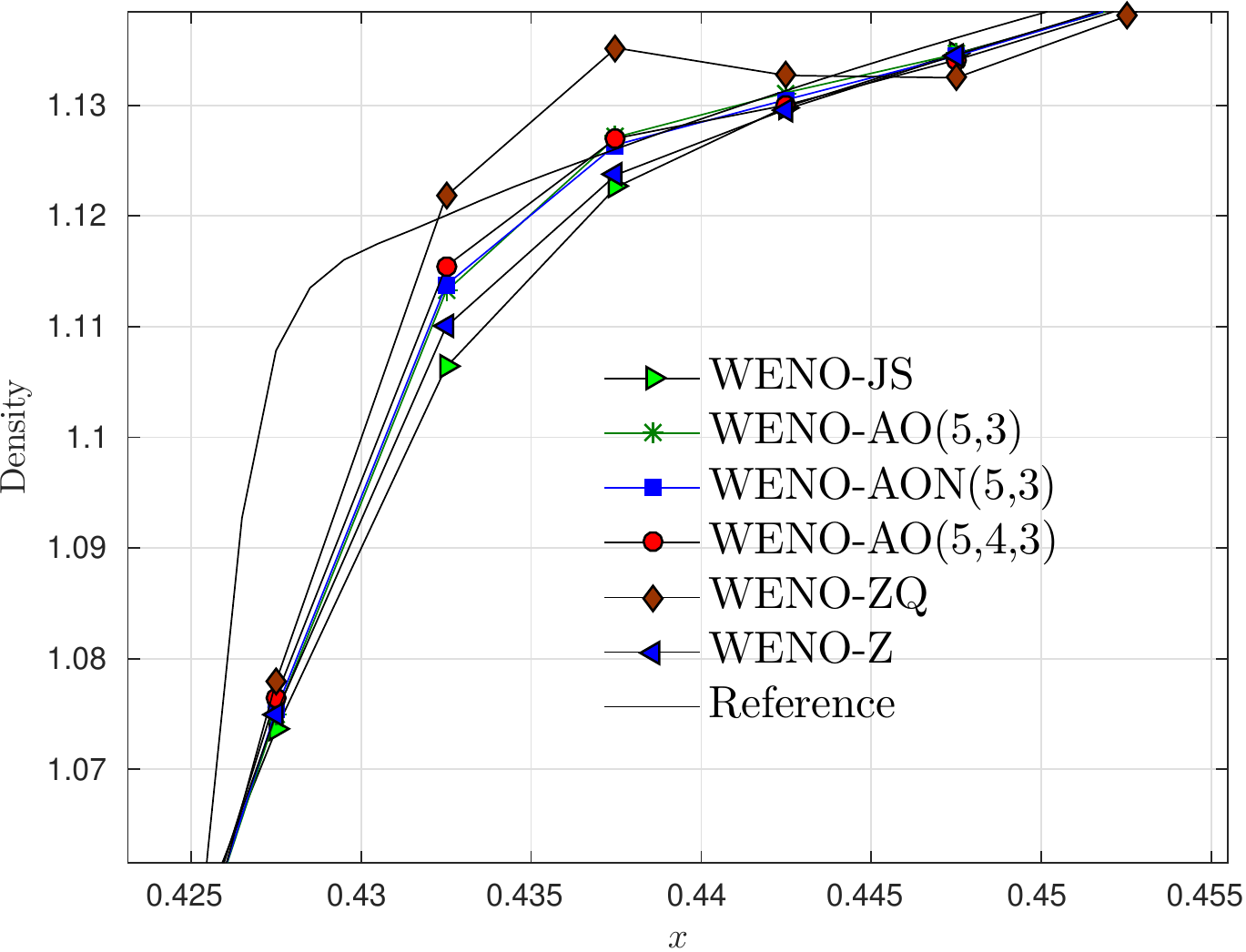}
 \end{tabular}
\end{center}
\caption{The cross sectional slices of density profile along the plane $y=0.5$ 
  computed with WENO-JS, WENO-Z, WENO-ZQ, WENO-AO(5,3), WENO-AON(5,3), and WENO-AO(5,4,3) schemes after the interaction (at time $T=0.35$) using uniform mesh of having
   resolution $200\times 200$.}
 \label{Figure.sv2}
\end{figure}

\begin{example}\label{ex:sv}{\rm (Shock-Vortex Interaction)
This 2D problem consists of the interaction of a left moving shock wave with a right moving vortex  \cite{pao-sal_81}, 
\cite{cha_99a}, \cite{ren-etal_03a}. 
Consider the system of 2D Euler equations \eqref{2deuler.crate} with 
  initial condition
  \begin{equation*}
\bold{V}_0=\left\{\begin{array}{ll}
                \bold{V}_L ~~~~   x<0.5, \\
                \bold{V}_R ~~ ~~ x\geq0.5.
                \end{array}\right.
\end{equation*}
where $\bold{V}_L$ and $\bold{V}_R$ denotes the left state  and right state of shock discontinuity, respectively.
The left state is taken as  $\bold{V}_L=(\rho_L, u_L,v_L,p_L)=(1, \sqrt{\gamma}, 0,1)$ and right state is given as
\begin{align*}
 p_R&=1.3, ~~~~\rho_R=\rho_L\Big(\frac{\gamma-1+(\gamma+1)p_R}{\gamma+1+(\gamma-1)p_R}\Big)\\
 u_R&=\sqrt(\gamma)+\sqrt{2}\Big(\frac{1-p_R}{\sqrt{\gamma-1+p_R(\gamma+1)}}\Big),~~v_R=0.0
\end{align*}
The left state $\bold{V}_L$ is superimposed onto a vortex given by following perturbations
\begin{align*}
 \delta u&= \epsilon \frac{(y-y_c)}{r_c}\exp(\alpha(1-r^2)),~~~\delta v =-\epsilon \frac{(x-x_c)}{r_c}\exp(\alpha(1-r^2))\\
 \delta \theta& = -\frac{\gamma-1}{4\alpha \gamma}\epsilon^2\exp(2\alpha(1-r^2)),~~ \delta s=0
\end{align*}
where $\delta \theta$ and $\delta s$ denotes the temprature and physical entropy, respectively, and $r^2=((x-x_c)^2+(y-y_c)^2)/r_c^2$. 

The domain is uniformly discretized with 200 grid points in both directions with transmissive boundary conditions, and the final time is taken as $T=0.35$.
 The parameters are chosen as $\epsilon_1=0.3$, $r_c=0.05$, $\alpha=0.204$ and center of vortex $(x_c,y_c)=(0.25, 0.5)$.
In Figure \ref{Figure.sv1} (a)-(b),  we have shown the initial and final positions of the shock and vortex in density profile, respectively, computed
 using WENO-AO(5,4,3) scheme.  In Figure \ref{Figure.sv2}, we have plotted the cross sectional slices of density profile along the plane $y=0.5$ 
  computed with WENO schemes after the interaction. The reference solution is computed using WENO-JS scheme over a uniform mesh  having resolution
   $1000\times 1000$. We can easily observe from the enlarged portion around shock, that  WENO-AO(5,4,3) resolves the shock in non-oscillatory manner and 
   performs better than other schemes. The WENO-AON(5,3) is comparable with WENO-AO(5,3) scheme and performs better than the WENO-JS and WENO-Z scheme.

 }
\end{example}

\begin{figure}[t!]
\begin{center}
 \begin{tabular}{cc}
\includegraphics[width=0.48\textwidth]{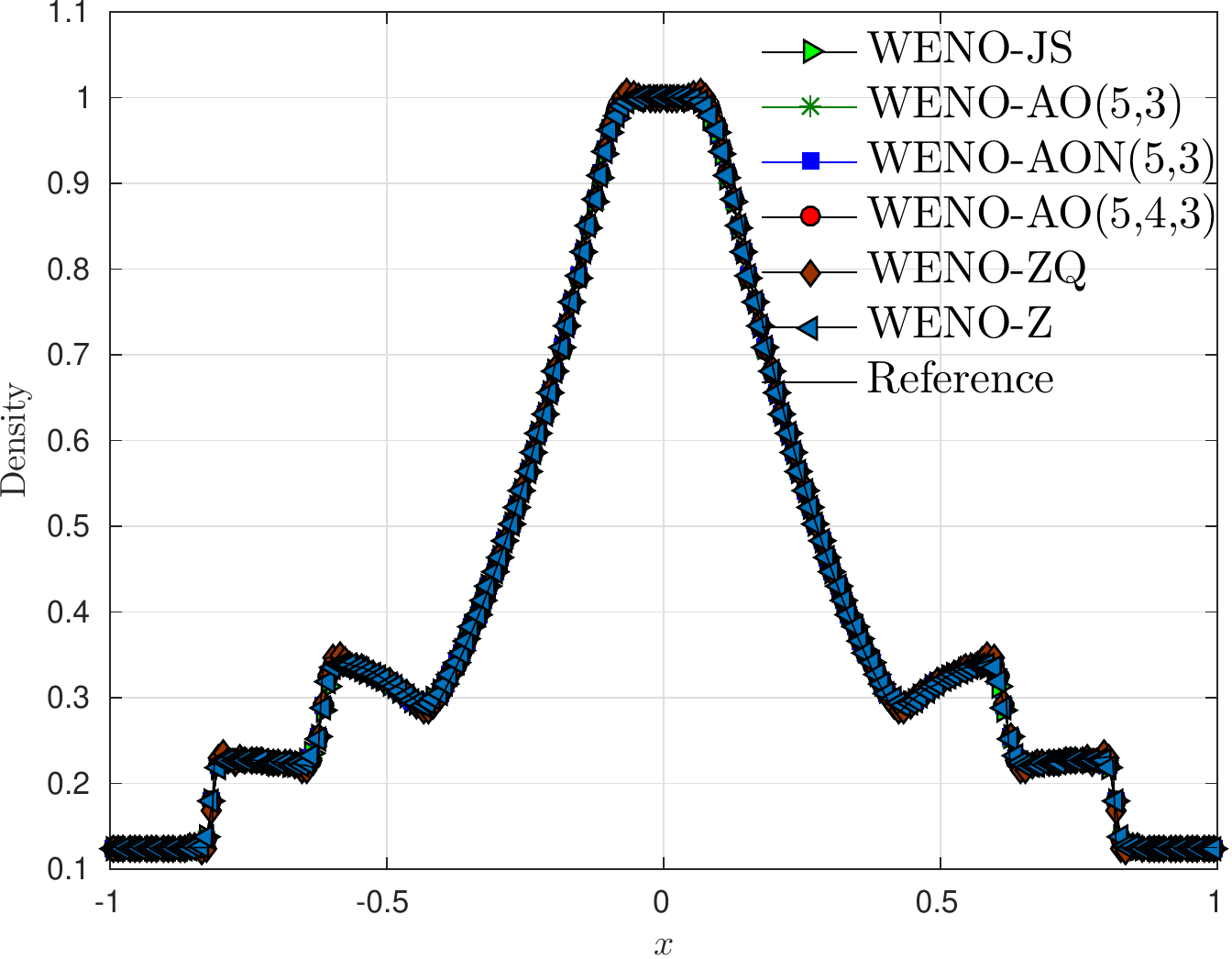}
\includegraphics[width=0.48\textwidth]{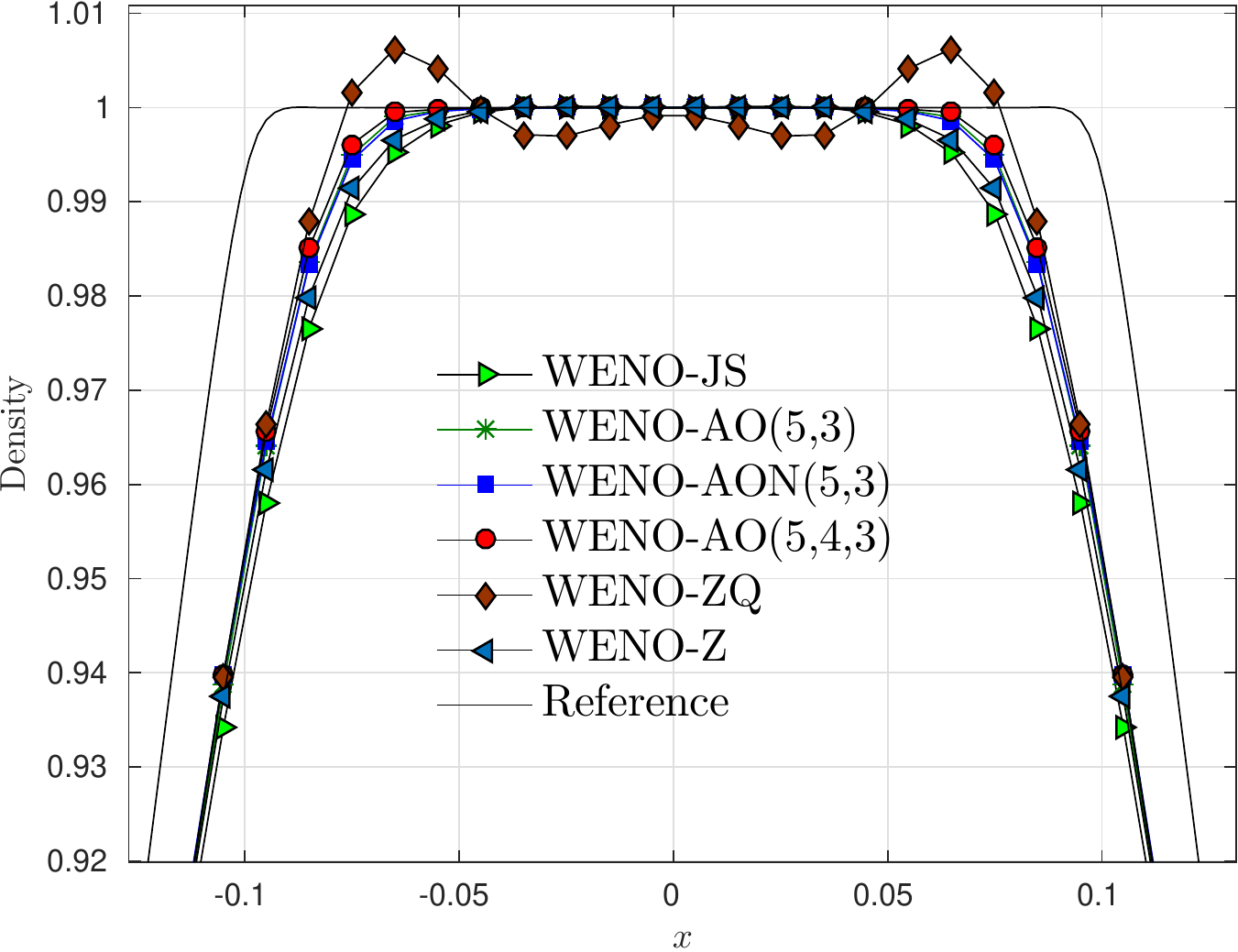}\\
\includegraphics[width=0.48\textwidth]{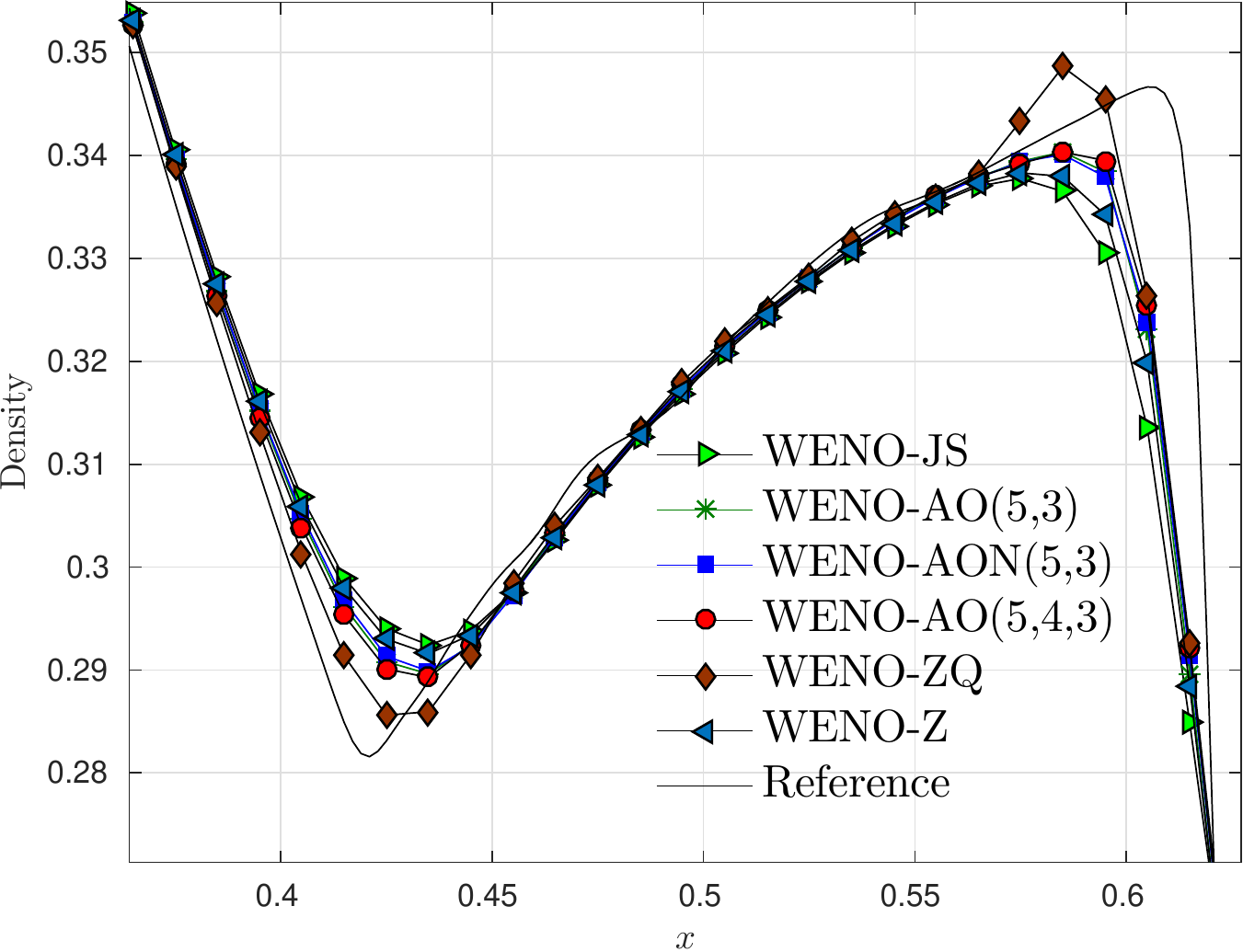}
\includegraphics[width=0.48\textwidth]{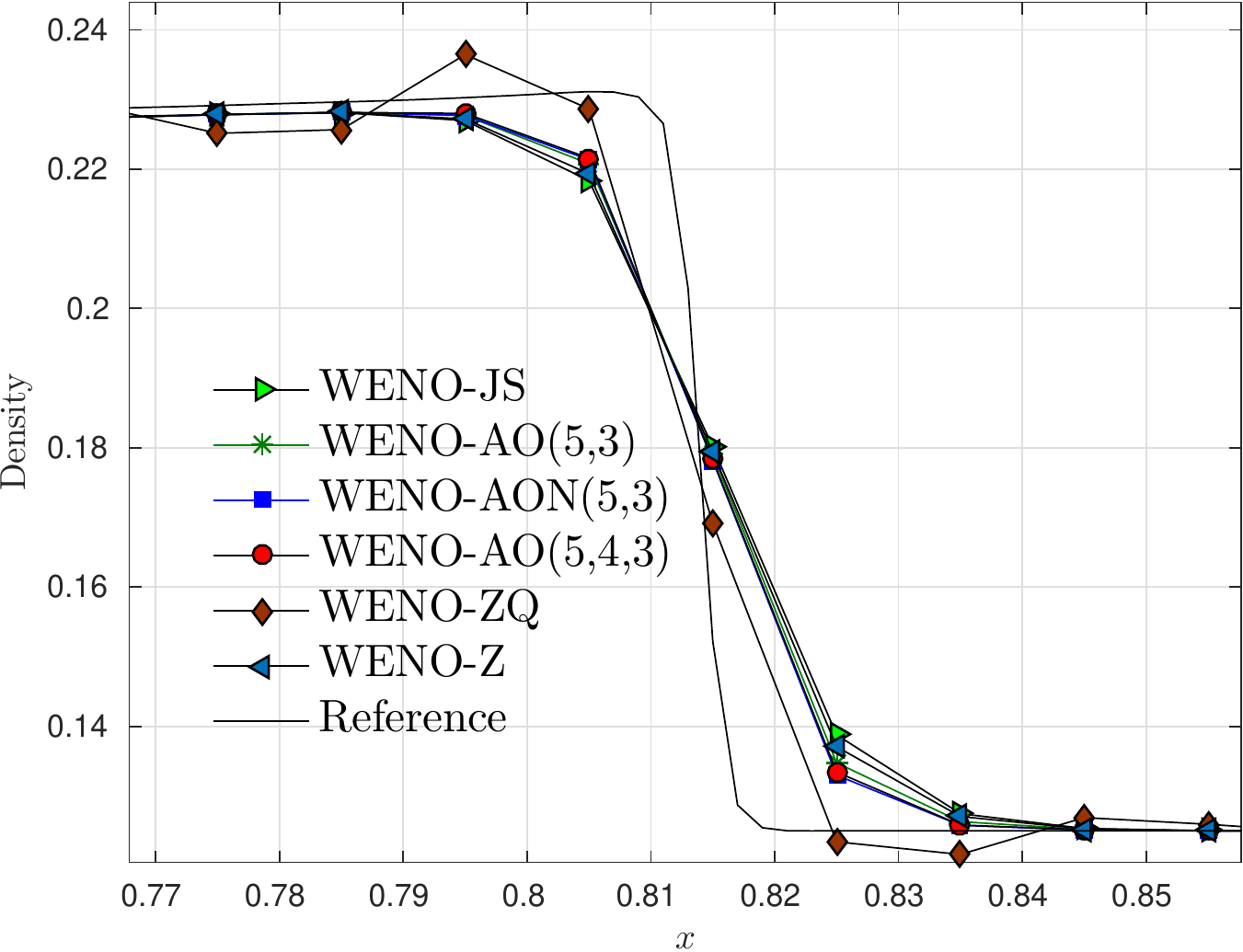}
 \end{tabular}
\end{center}
\caption{The cross sectional slices of density profile along the plane $y=0.0$ 
  computed with WENO-JS, WENO-Z, WENO-ZQ, WENO-AO(5,3), WENO-AON(5,3), and WENO-AO(5,4,3) schemes at time $T=0.25$ using uniform mesh of having
   resolution $200\times 200$.}
 \label{Figure.explosion}
\end{figure}

\begin{figure}
\begin{center}
\begin{tabular}{cc}
\includegraphics[width=0.48\textwidth]{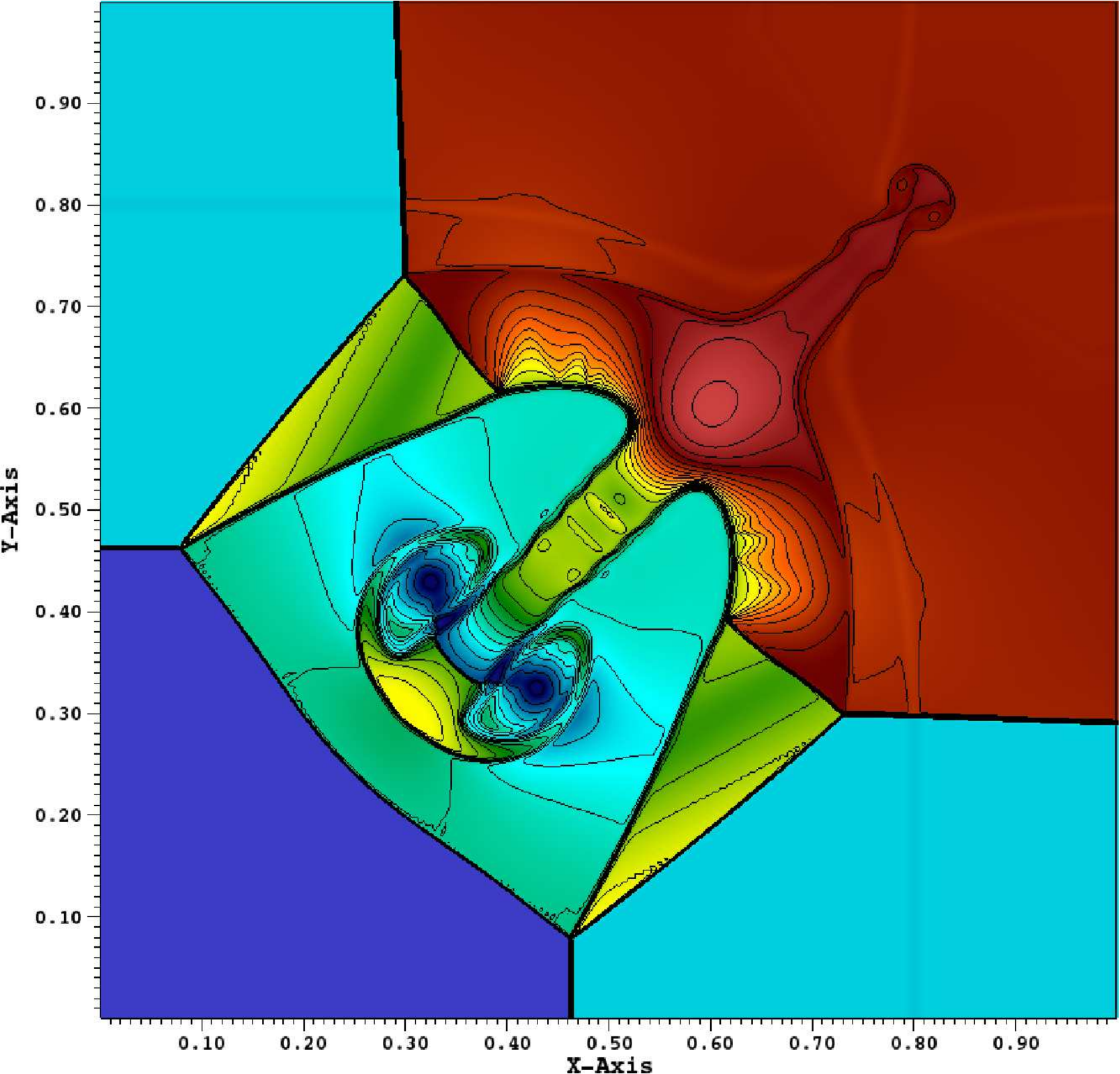} &
 \includegraphics[width=0.48\textwidth]{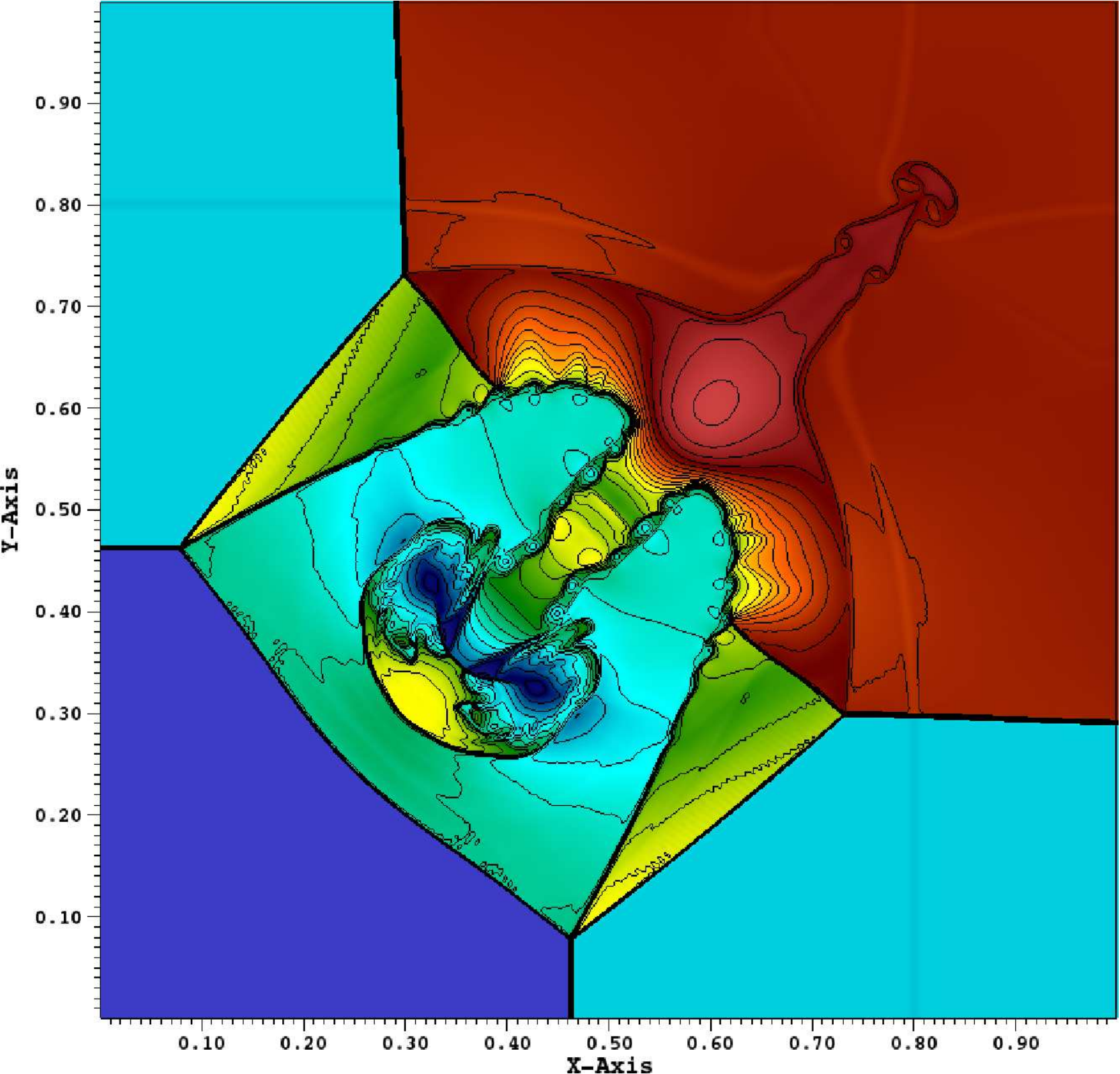} \\
 (a) WENO-JS & (b) WENO-Z\\
 \includegraphics[width=0.48\textwidth]{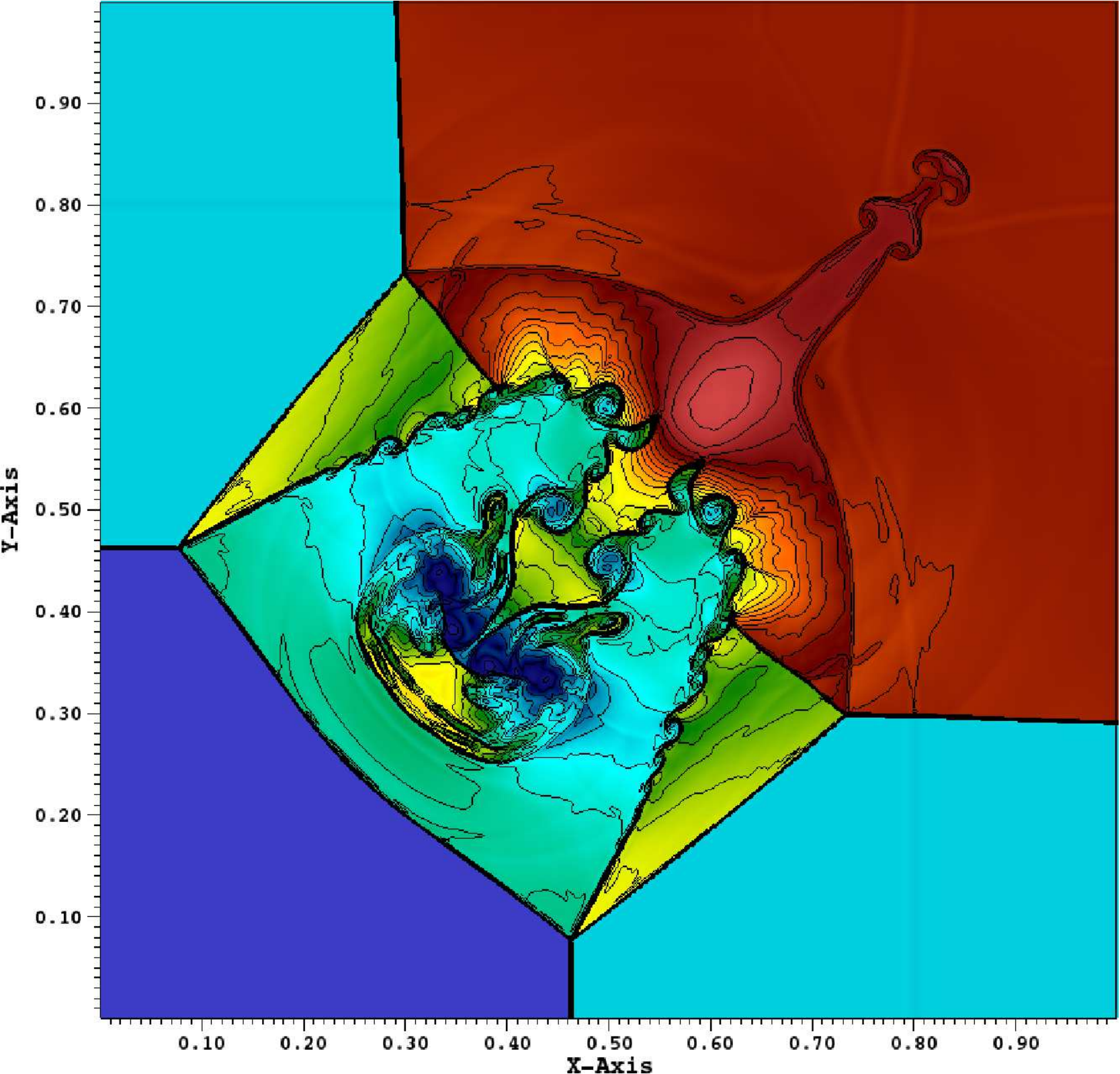} &
 \includegraphics[width=0.48\textwidth]{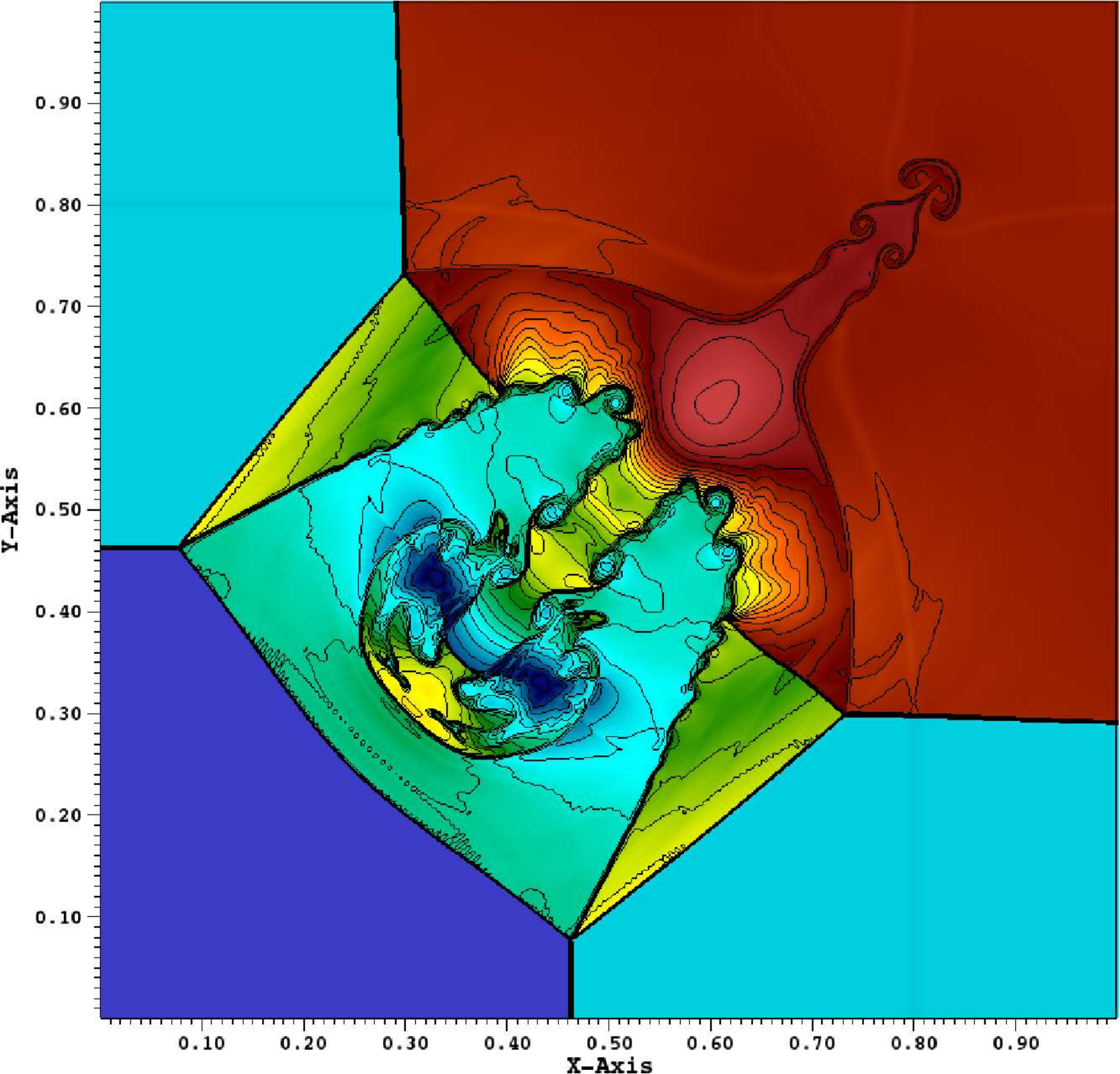} \\
 (c) WENO-ZQ & (d) WENO-AON(5,3)\\
 \includegraphics[width=0.48\textwidth]{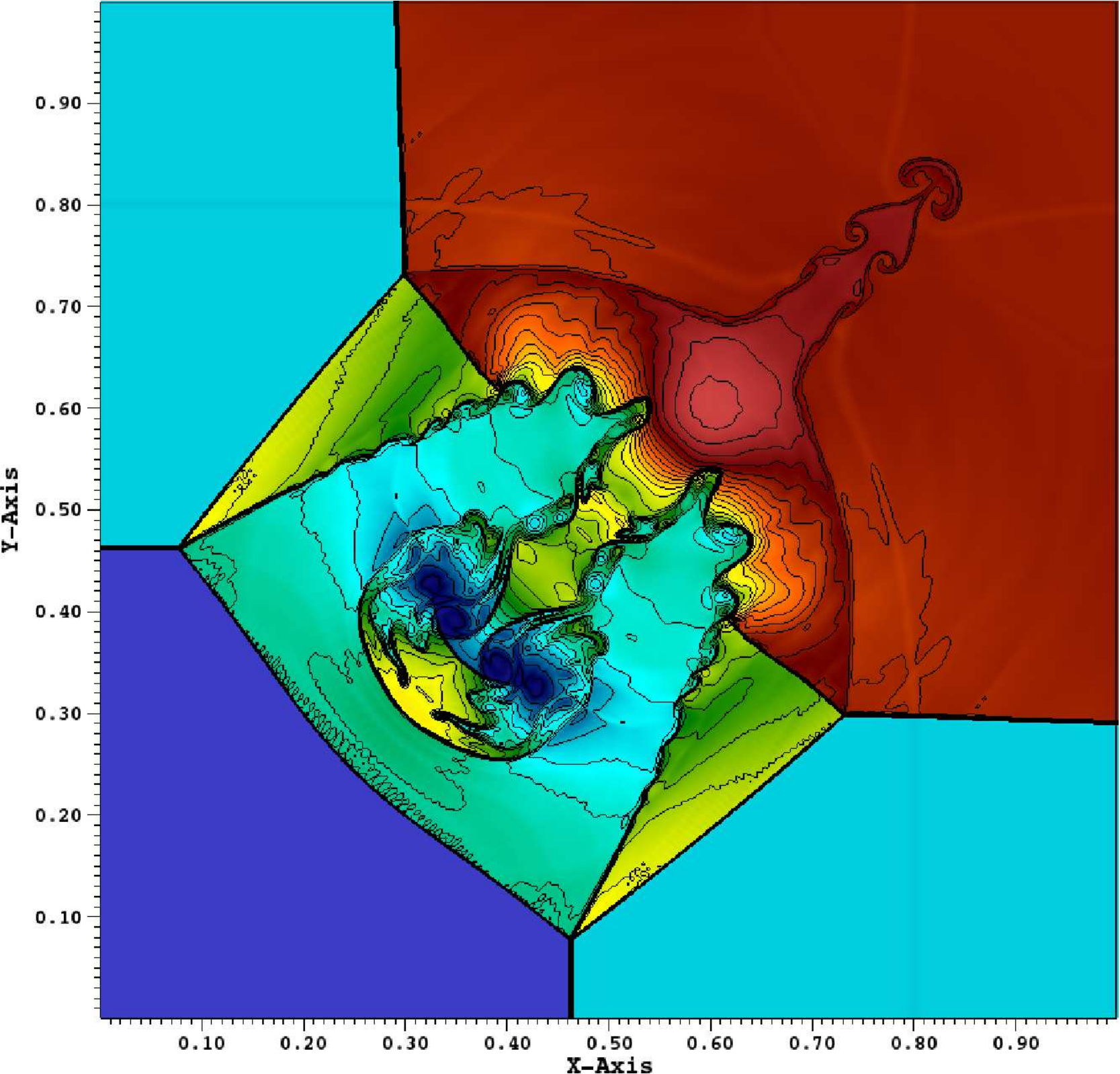} &
 \includegraphics[width=0.48\textwidth]{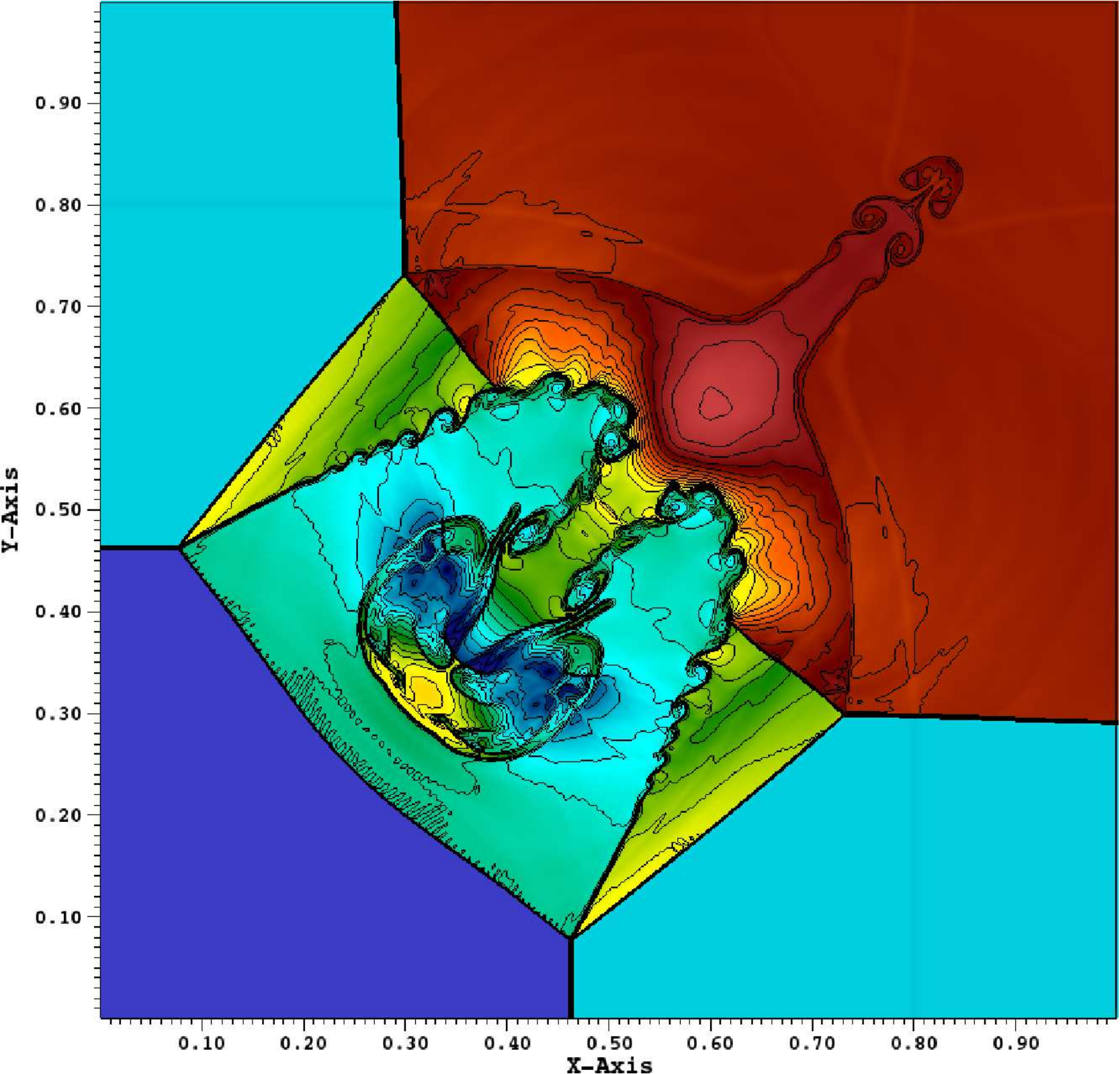} \\
 (e) WENO-AO(5,3) & (f) WENO-AO(5,4,3)
\end{tabular}
\end{center}
 \caption{Density contours for the Example \ref{ex:rp}  with 30 contour lines with range from 0.1 to 1.8, computed using WENO schemes at
 time $T=1.95$ over the domain $[0,1]\times[0,1]$  with mesh grid of size $800\times 800$.}
 \label{riemann.2d}
\end{figure}

\begin{figure}
\begin{center}
\begin{tabular}{cc}
\includegraphics[width=0.48\textwidth]{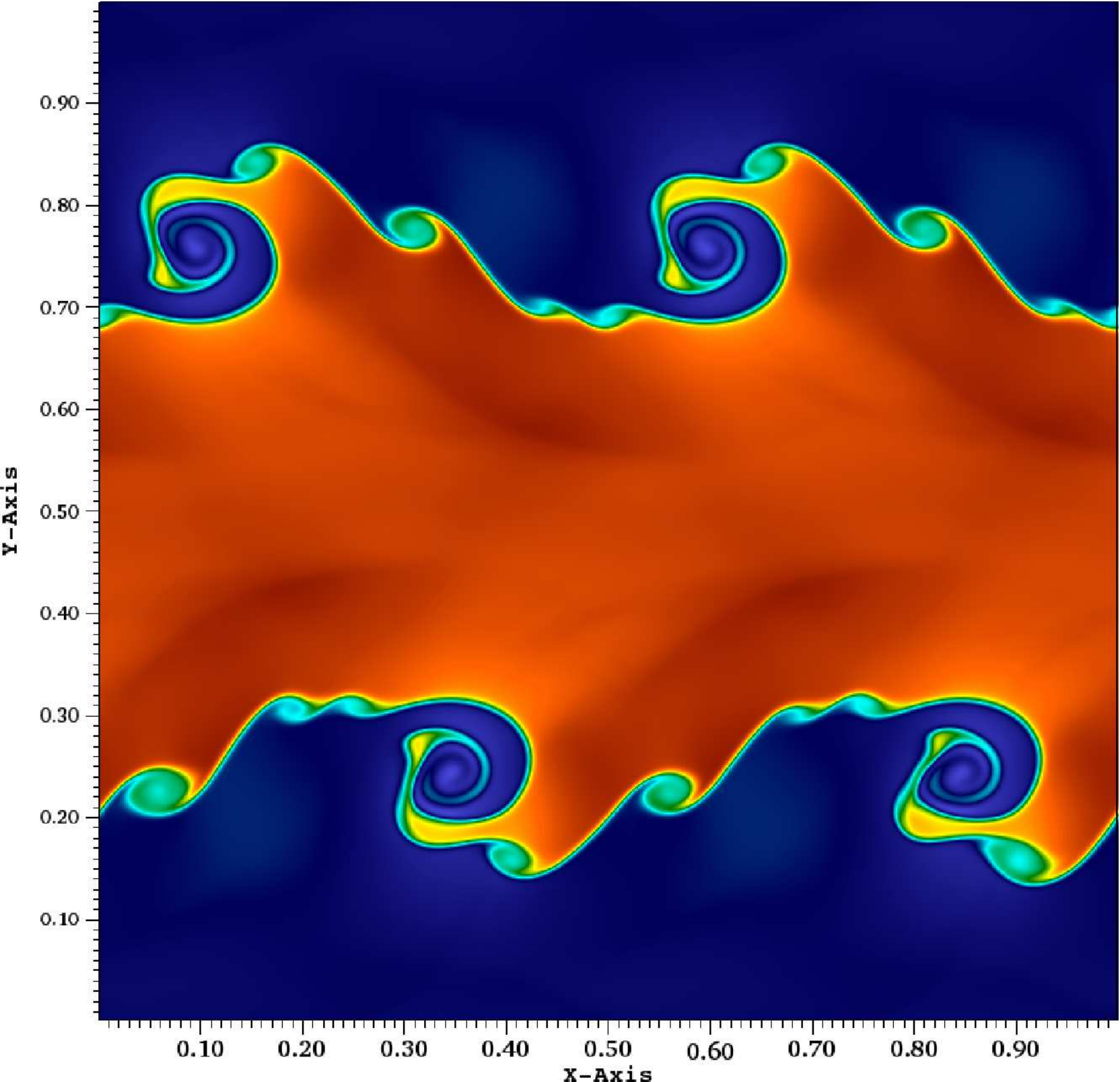} &
 \includegraphics[width=0.48\textwidth]{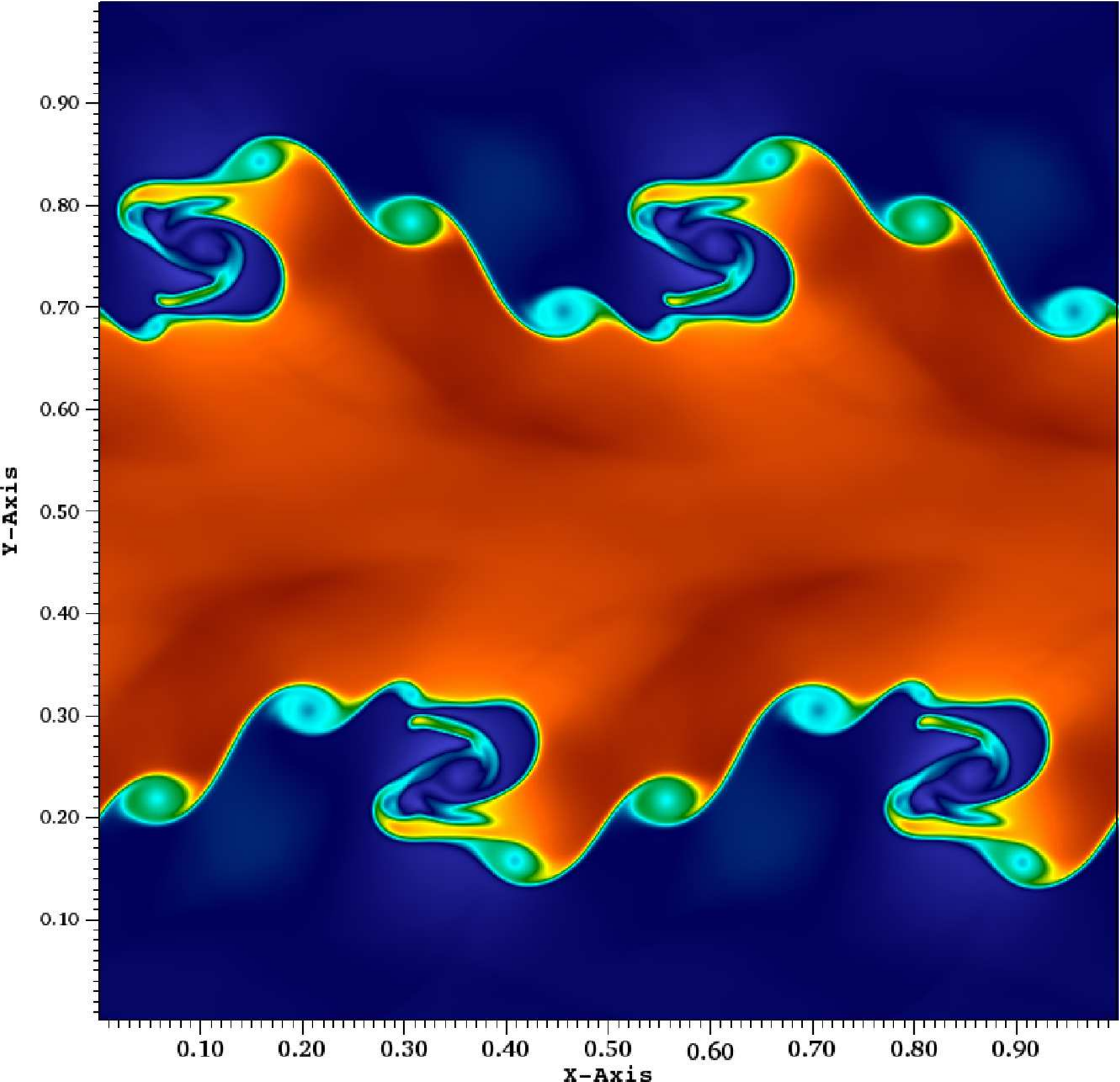} \\
 (a) WENO-JS & (b) WENO-Z\\
 \includegraphics[width=0.48\textwidth]{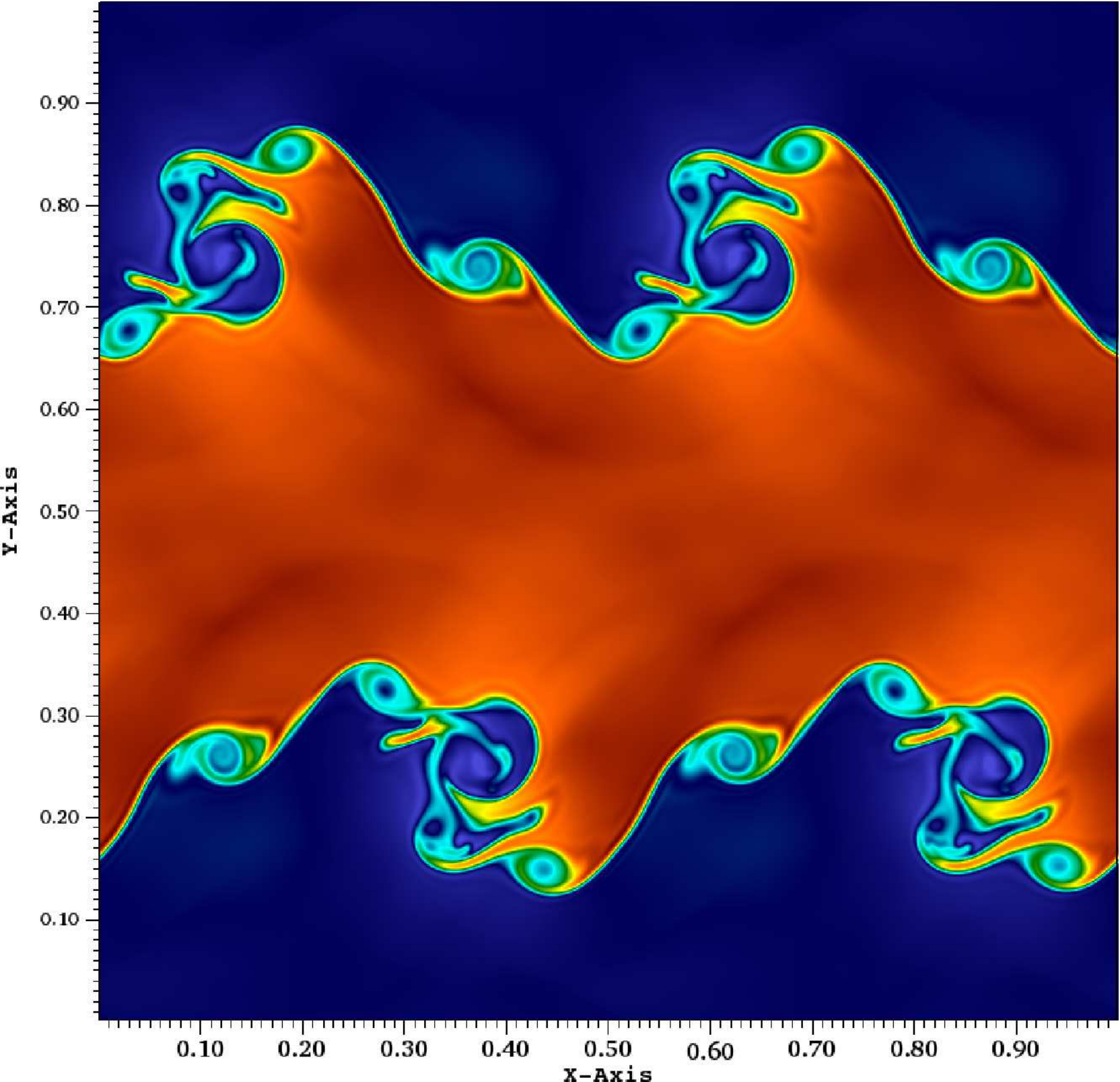} &
 \includegraphics[width=0.48\textwidth]{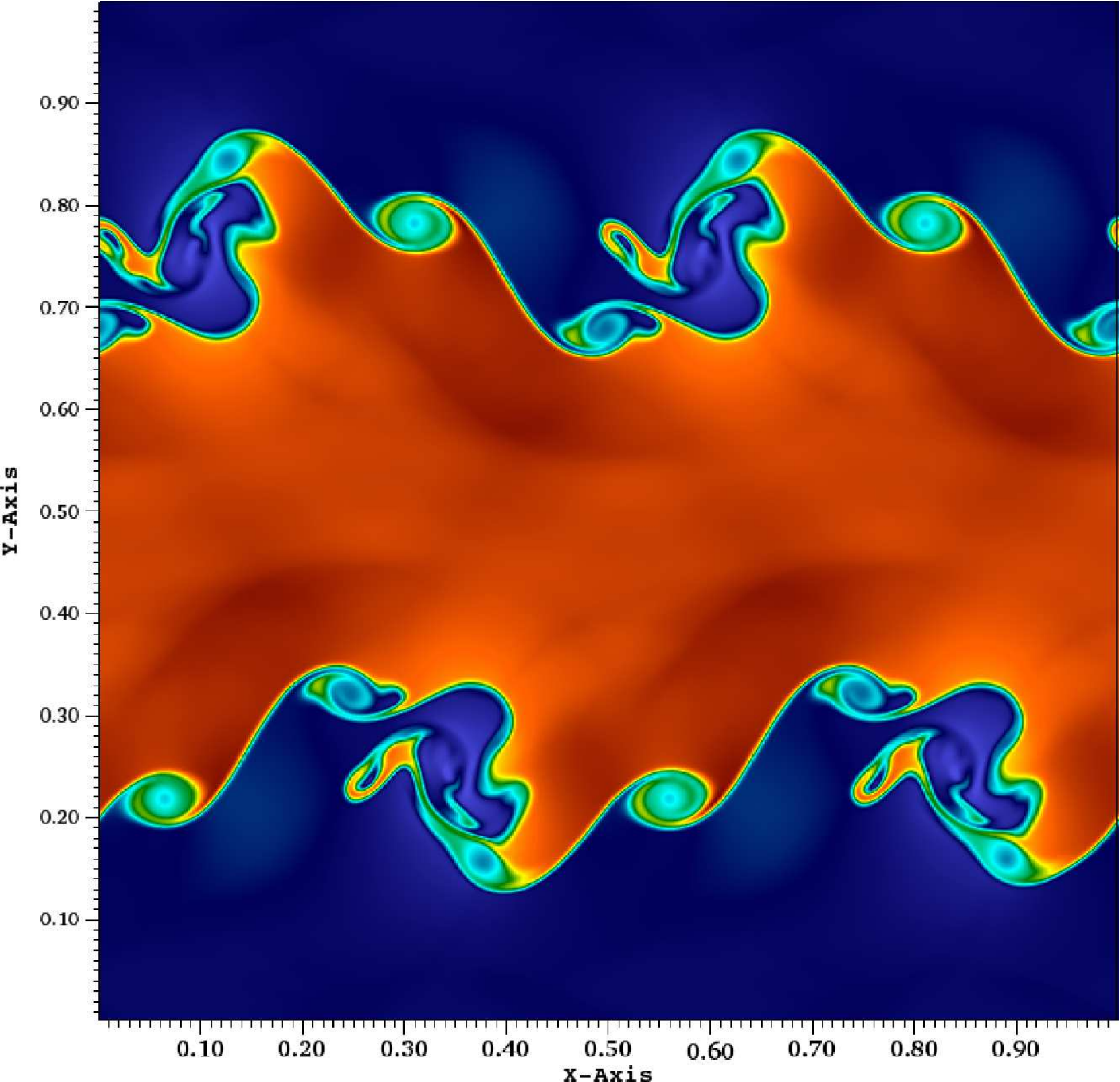} \\
 (c) WENO-ZQ & (d) WENO-AON(5,3)\\
 \includegraphics[width=0.48\textwidth]{kv_wenoao53n} &
 \includegraphics[width=0.48\textwidth]{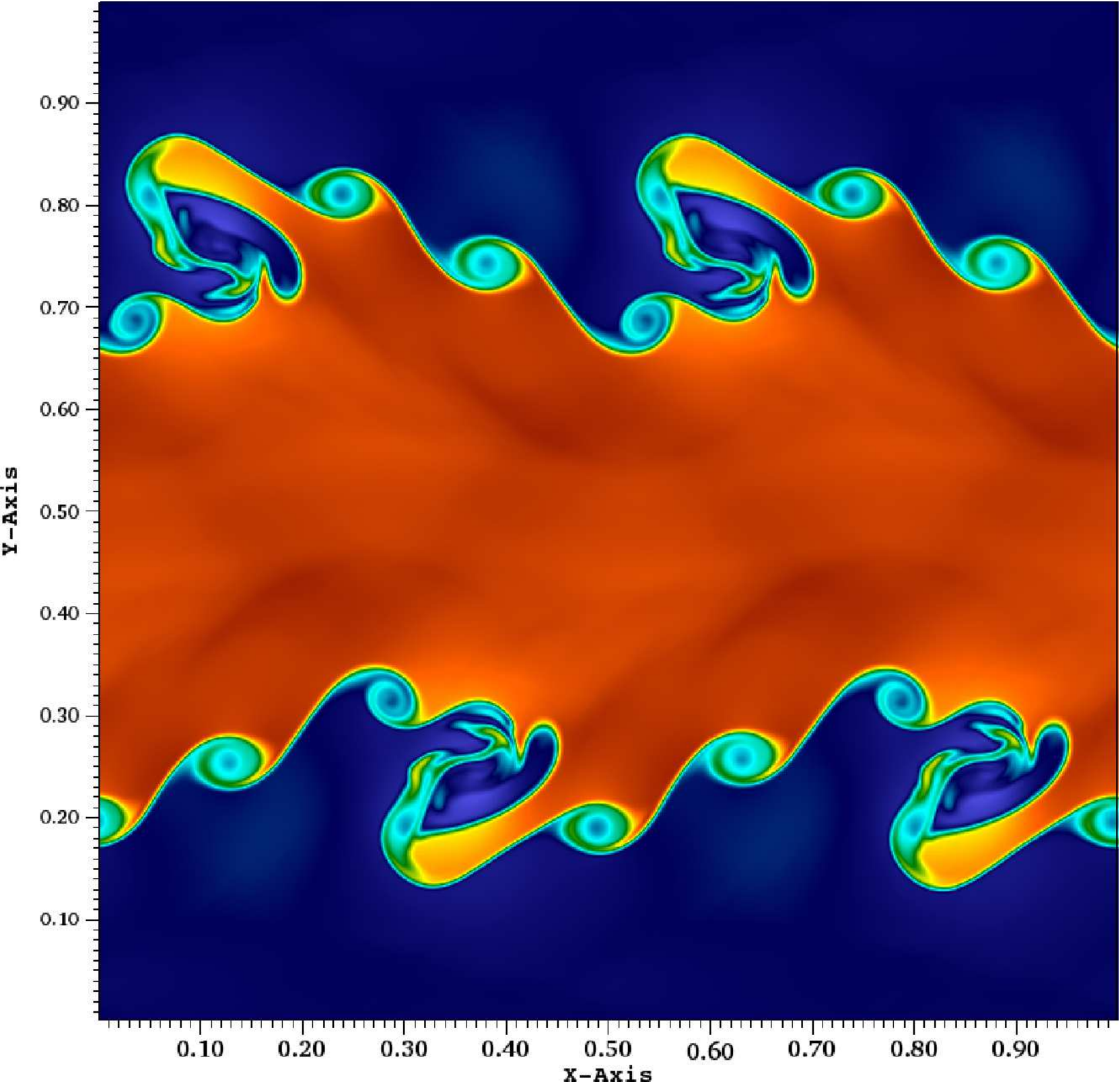} \\
 (e) WENO-AO(5,3) & (f) WENO-AO(5,4,3)
\end{tabular}
\end{center}
\caption{Density plots for the Example \ref{ex:kh} with range from 0.8 to 2.2, computed using WENO schemes at
 time $T=0.8$ over the domain $[0,1]\times[0,1]$  with mesh grid of size $512\times 512$.}
 \label{Figure.kh}
\end{figure}

\begin{figure}
\begin{center}
 \begin{tabular}{ccc}
  \includegraphics[width=0.3\textwidth]{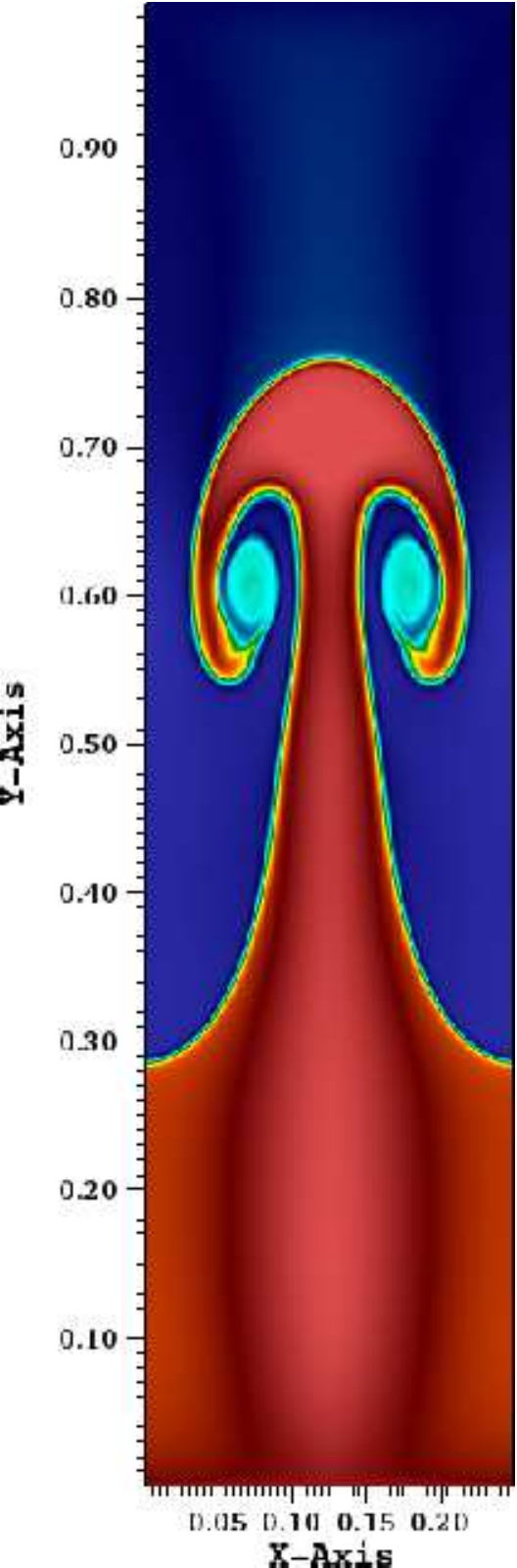}    &
  \includegraphics[width=0.3\textwidth]{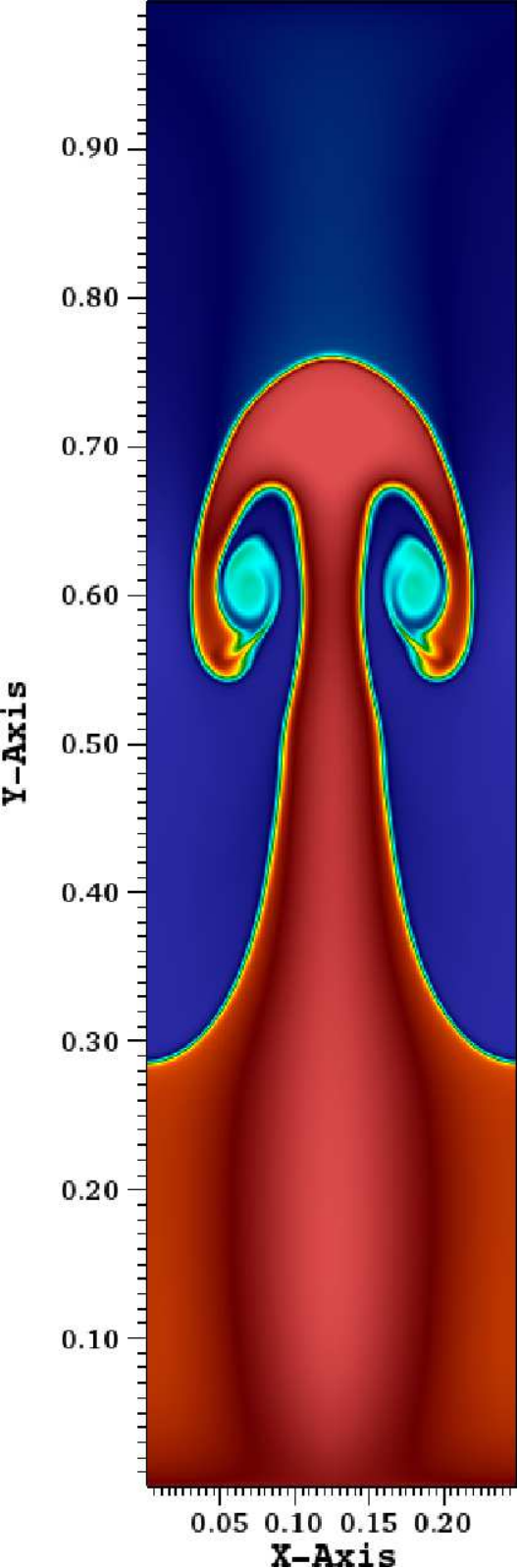}     &
  \includegraphics[width=0.3\textwidth]{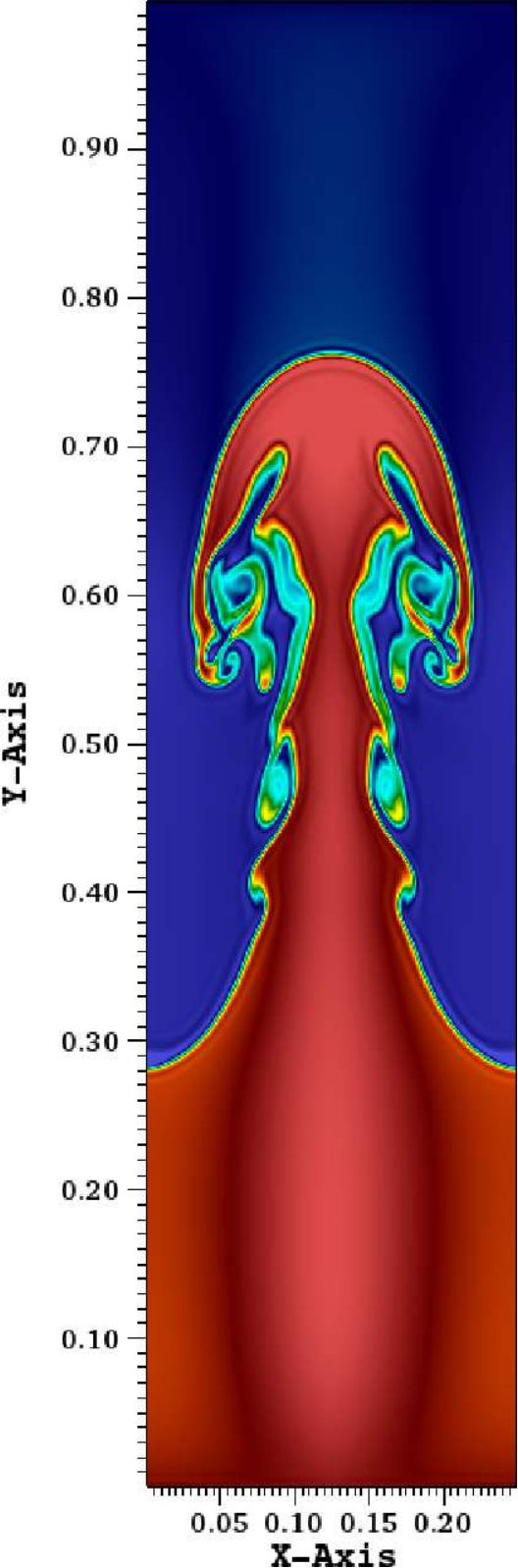}    \\
  (a) WENO-JS & (b) WENO-Z &  (c) WENO-ZQ   \\
   \includegraphics[width=0.3\textwidth]{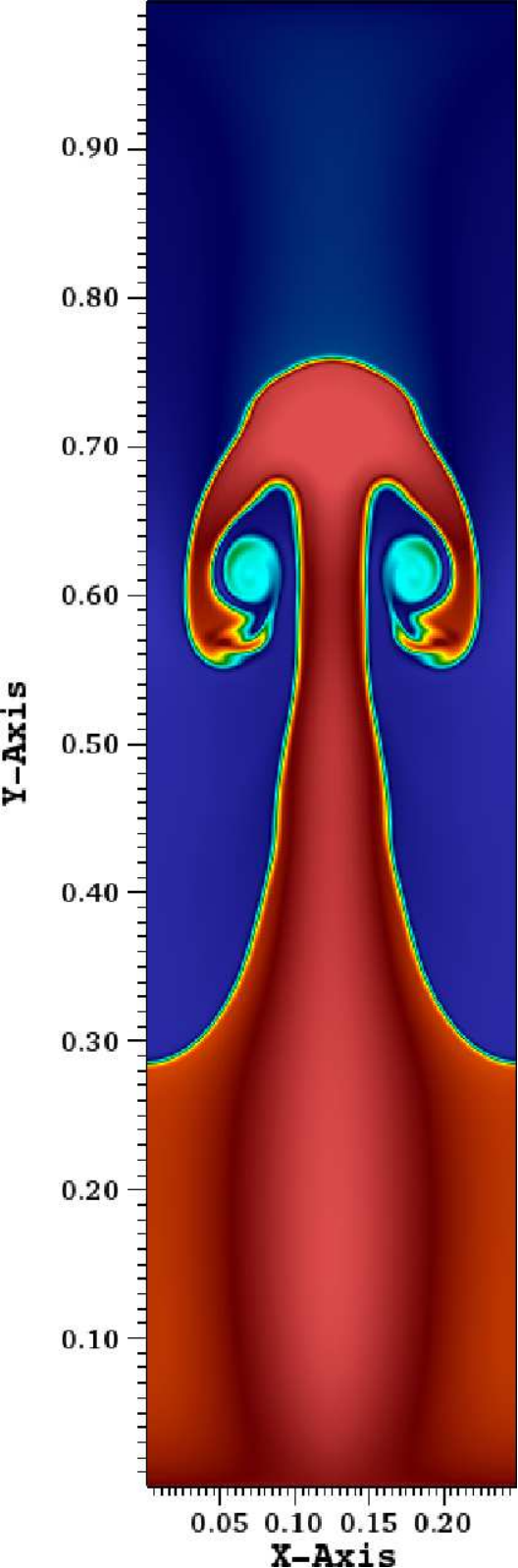}
  &  \includegraphics[width=0.3\textwidth]{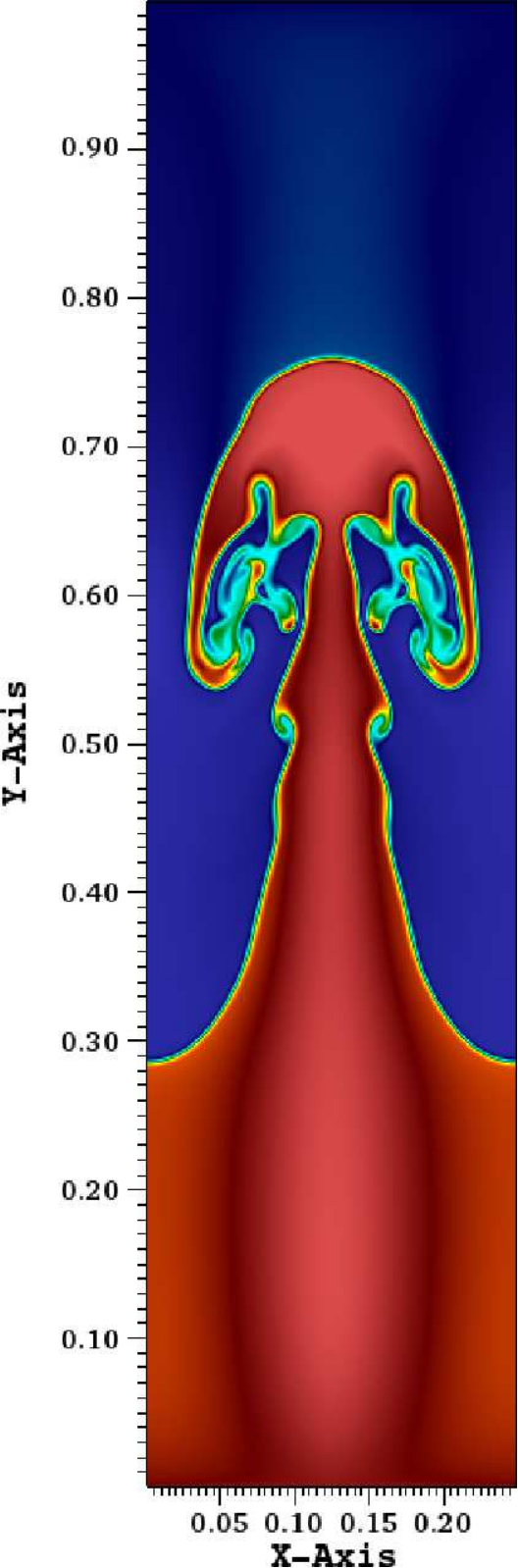}
  &  \includegraphics[width=0.3\textwidth]{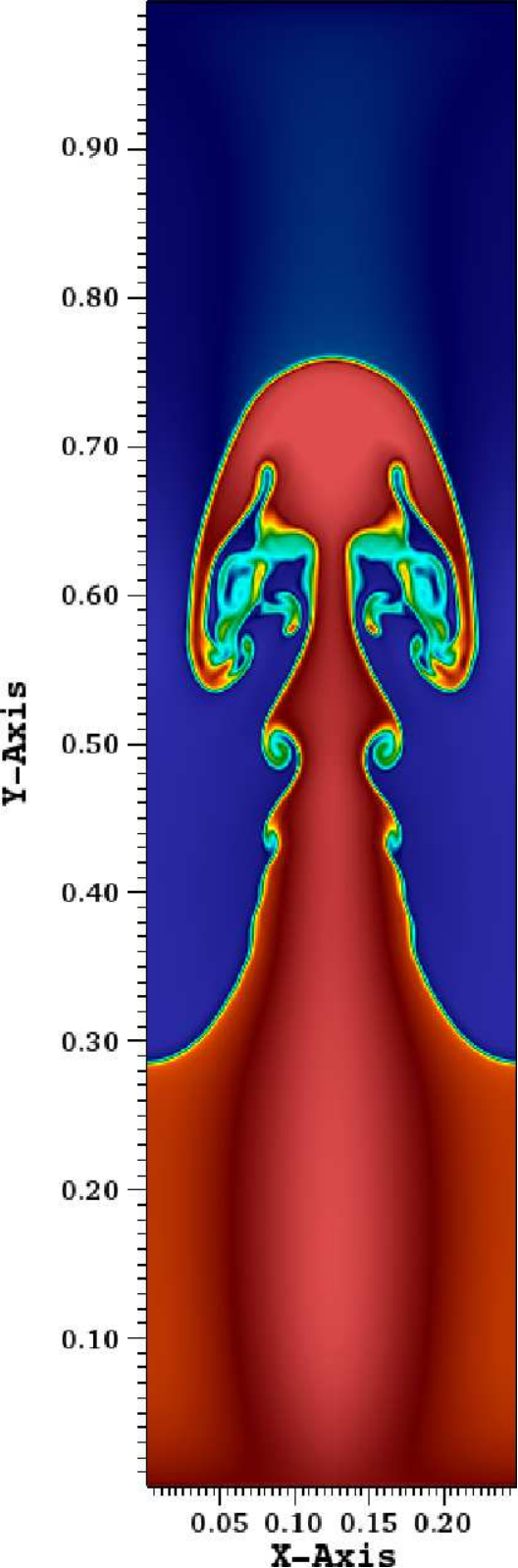}
   \\
  WENO-AON(5,3) & (e) WENO-AO(5,3) & (f) WENO-AO(5,4,3)
 \end{tabular}

\end{center}
 \caption{Density plots for the Example \ref{ex:rt}  with range from 0.8 to 2.4, computed using WENO schemes at
 time $T=1.95$ over the domain $[0,1/4]\times[0,1]$  with mesh grid of size $200\times 800$.}
 \label{Figure.rt}
\end{figure}

\begin{figure}
\begin{tabular}{c}
\includegraphics[width=1.0\textwidth]{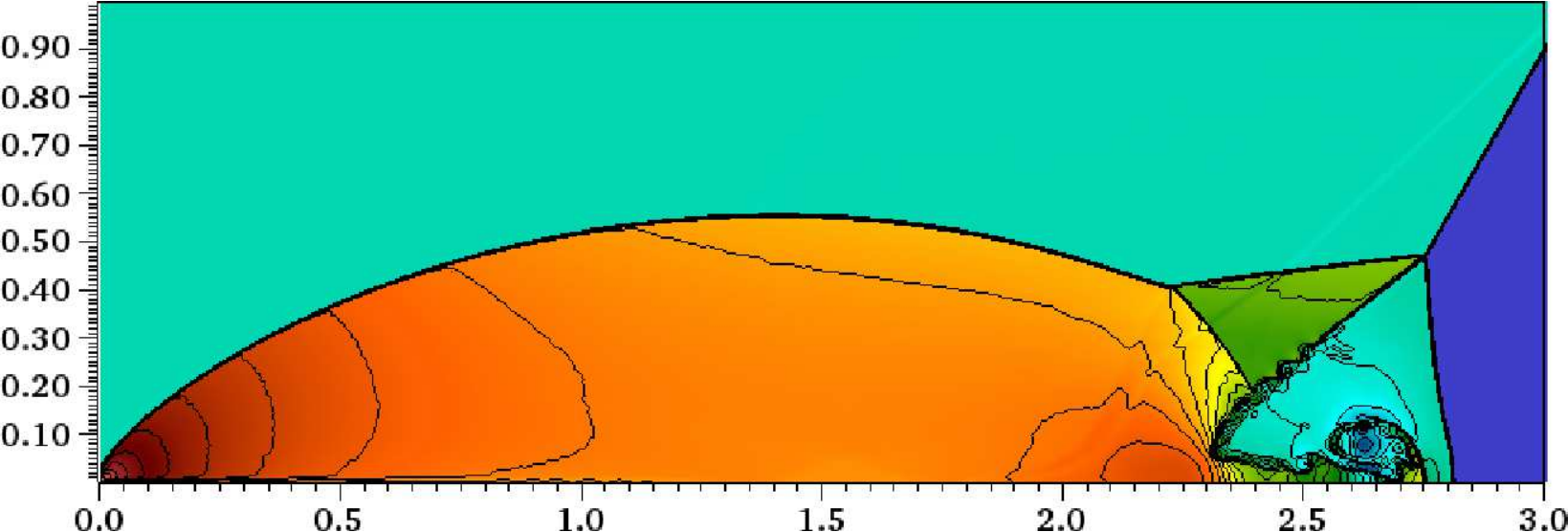}  \\
(a) WENO-AON(5,3)\\
\includegraphics[width=1.0\textwidth]{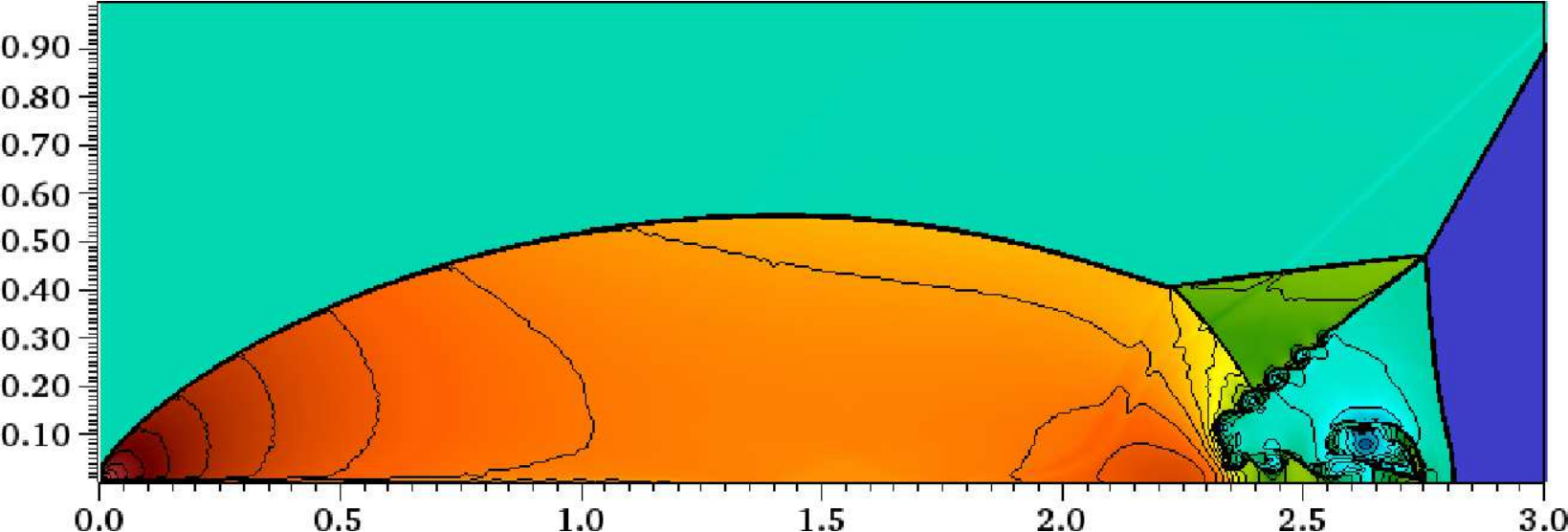}  \\
(b) WENO-AO(5,3)\\
\includegraphics[width=1.0\textwidth]{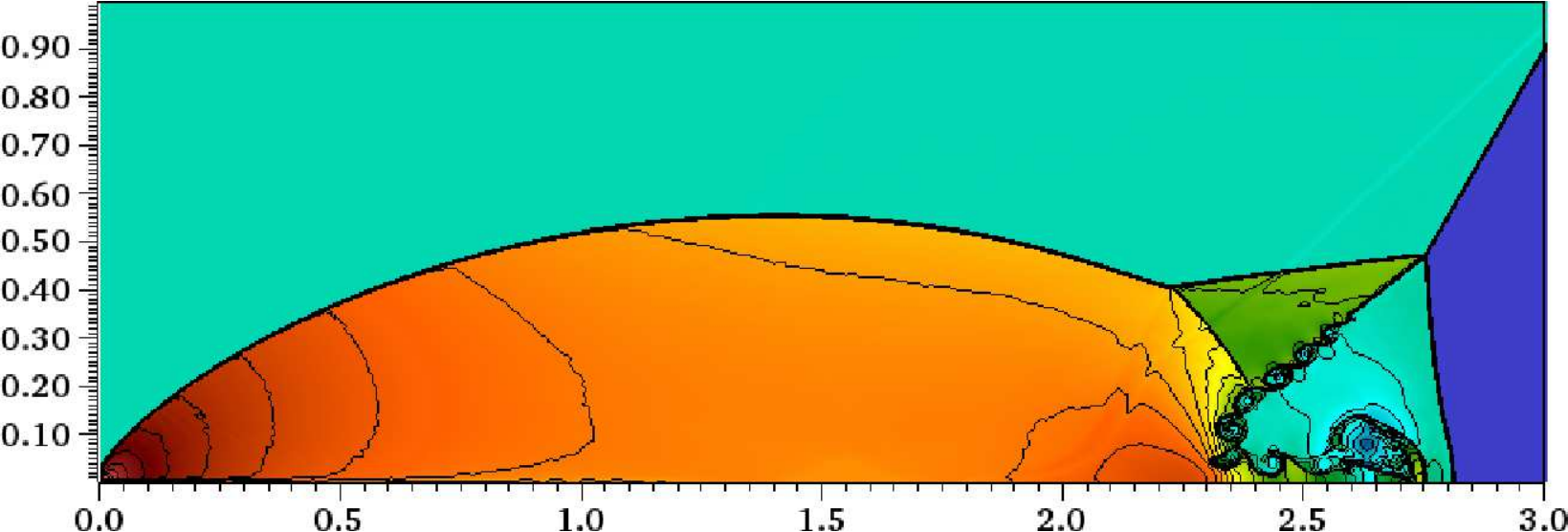} \\
(c) WENO-AO(5,4,3)\\
\end{tabular}
\caption{Density contours of the double-Mach shock reflection problem with 30 contour lines with range from 1 to 22.5, computed using WENO-AON(5,3),
WENO-AO(5,3), and WENO-AO(5,4,3) schemes at
 time $T=0.2$ over the domain $[0,4]\times[0,1]$ (shown over $[0,3]\times[0,1]$)  with mesh grid of size $1600\times 400$.}
 \label{Figure.Doublemach1}
\end{figure}

\begin{figure}
\begin{center}
\begin{tabular}{cc}
\includegraphics[width=0.48\textwidth]{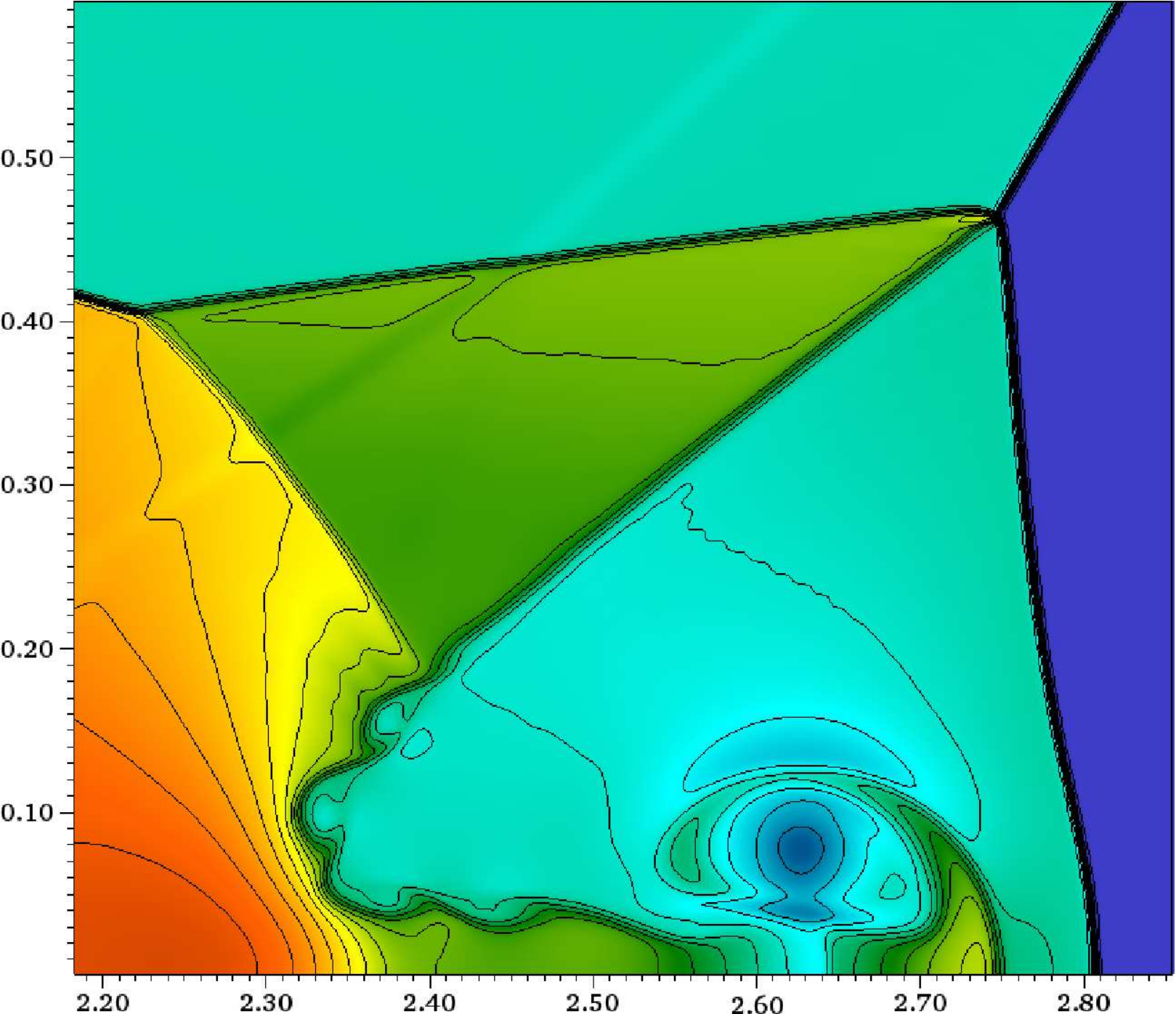} &
 \includegraphics[width=0.48\textwidth]{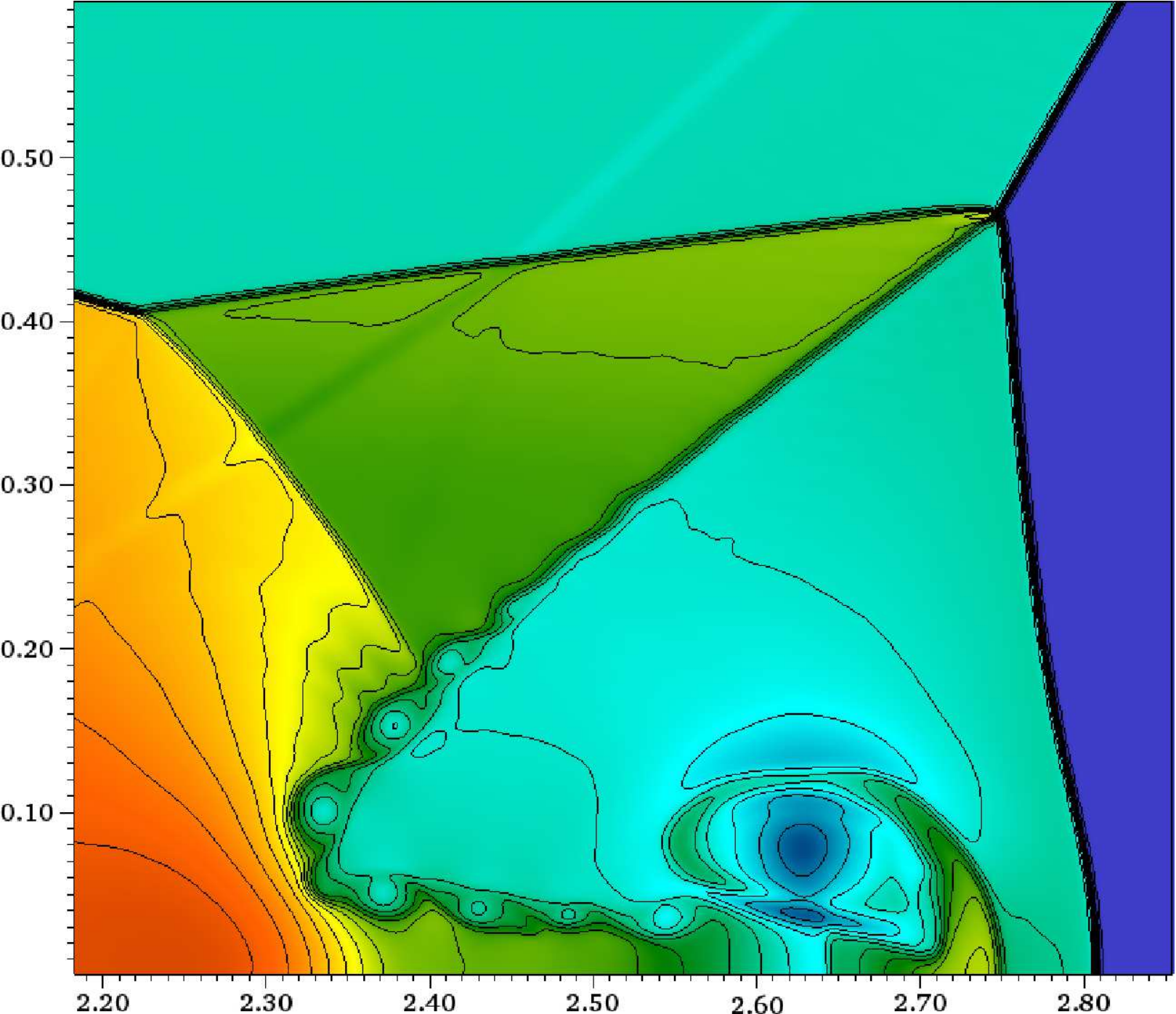} \\
 (a) WENO-JS & (b) WENO-Z\\
 \includegraphics[width=0.48\textwidth]{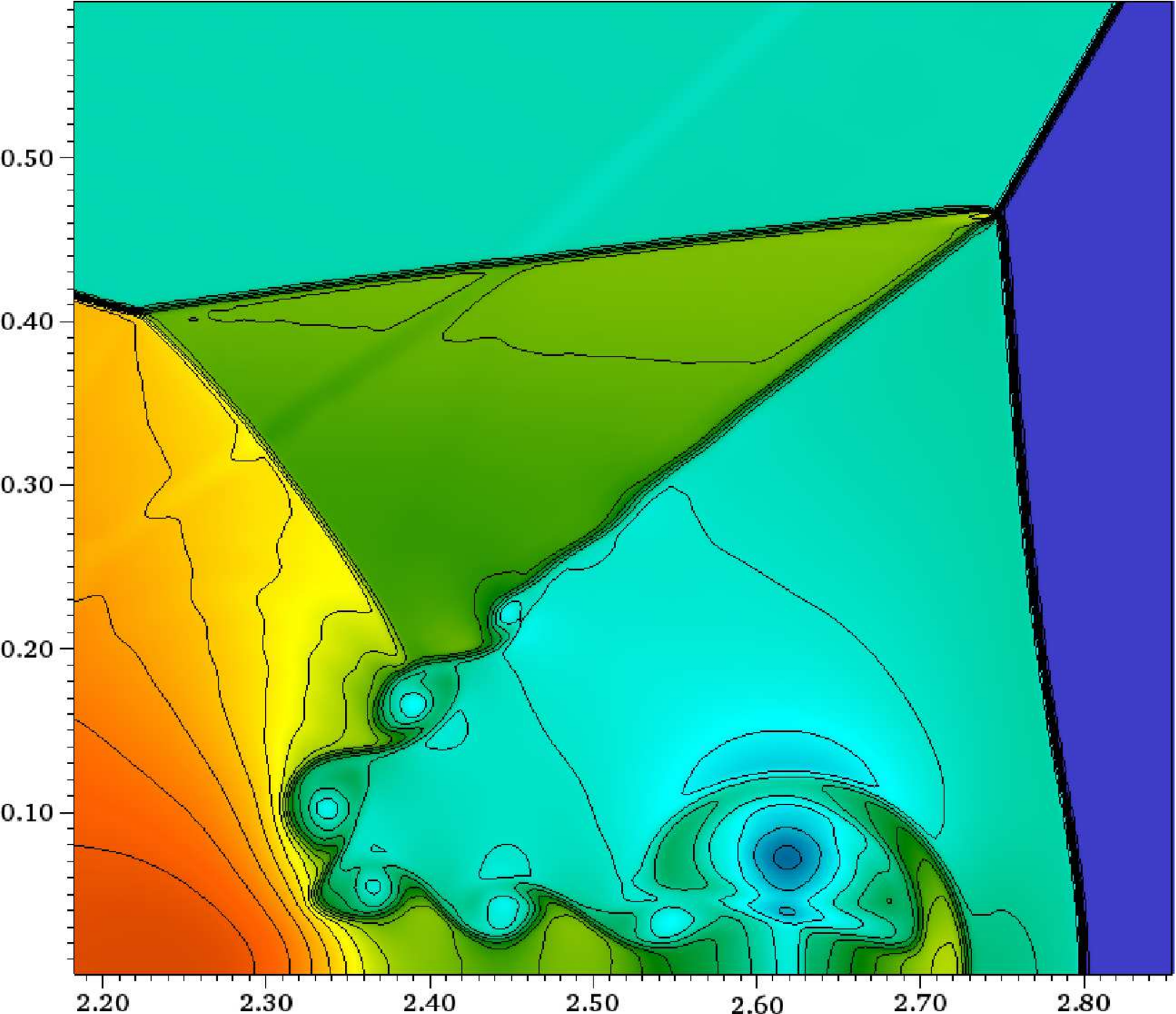} &
 \includegraphics[width=0.48\textwidth]{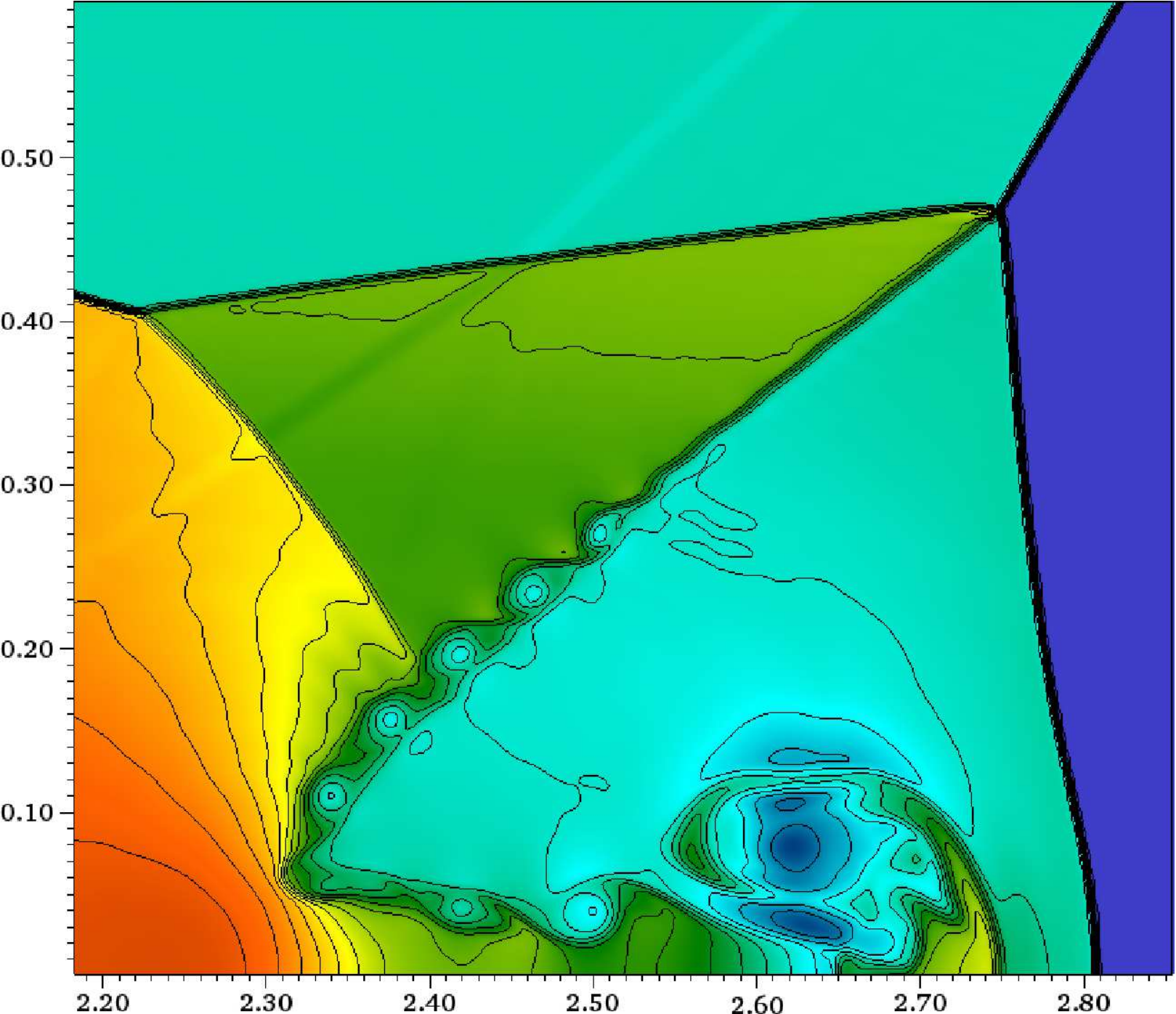} \\
 (c) WENO-ZQ & (d) WENO-AON(5,3)\\
 \includegraphics[width=0.48\textwidth]{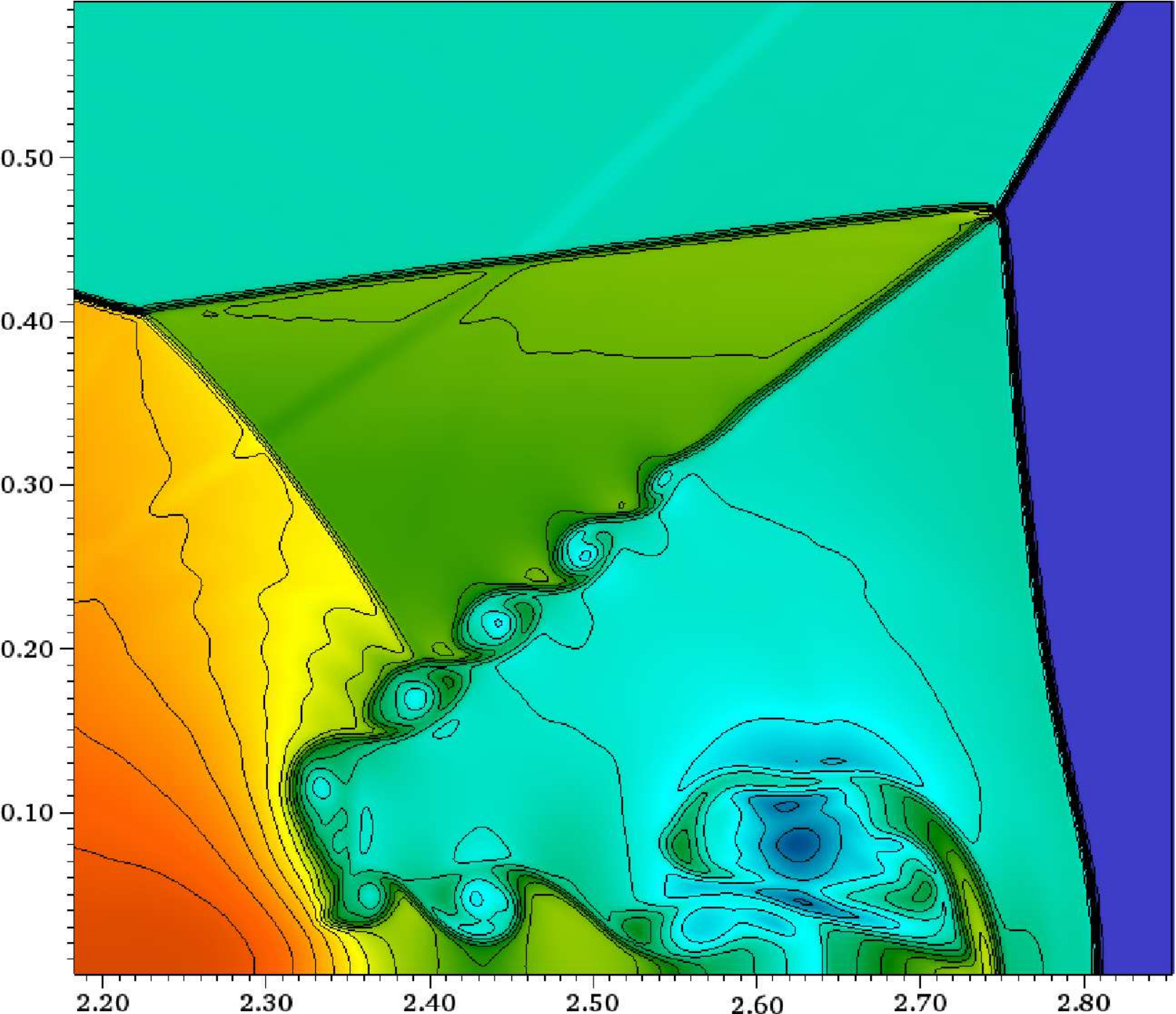} &
 \includegraphics[width=0.48\textwidth]{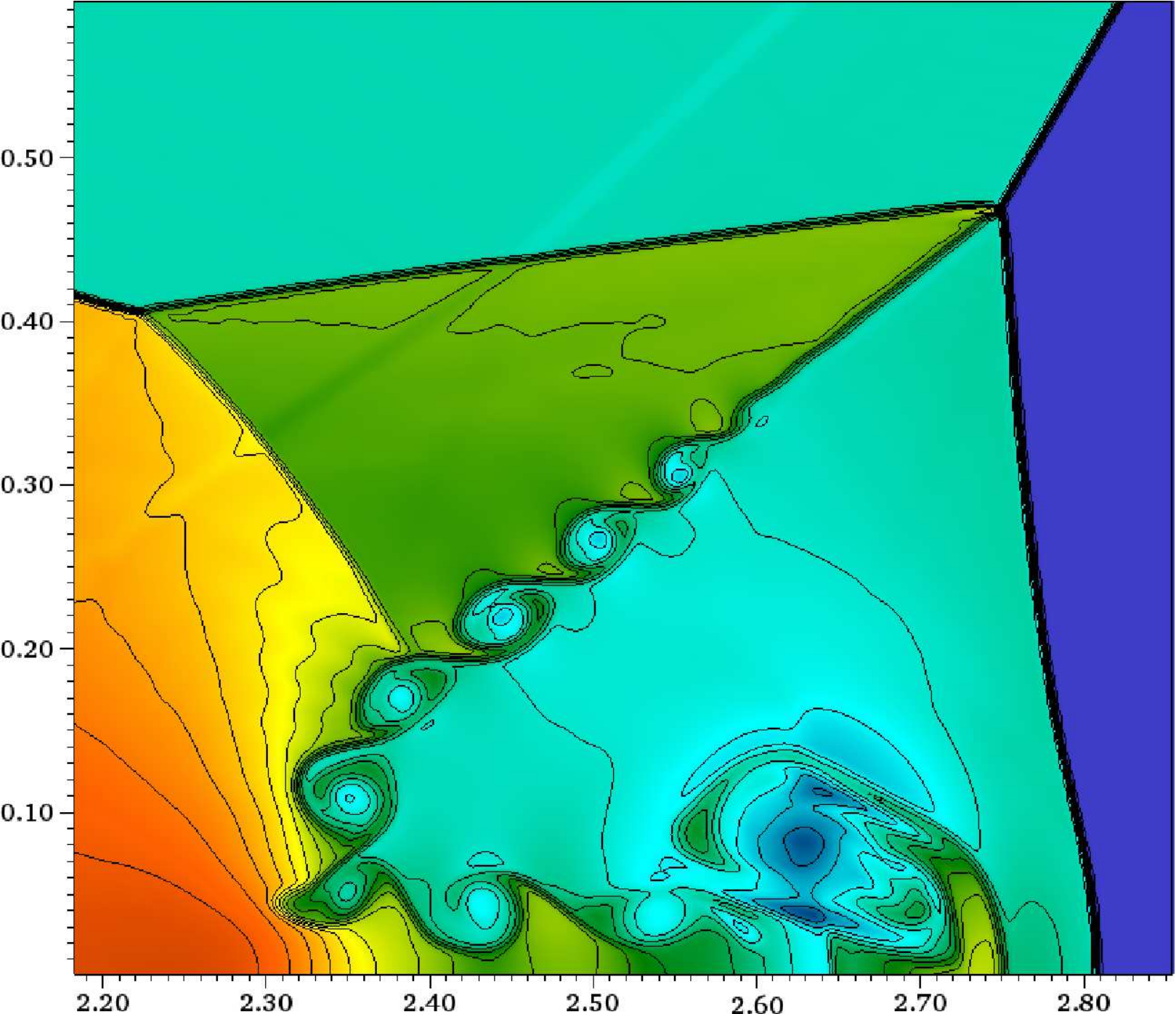} \\
 (e) WENO-AO(5,3) & (f) WENO-AO(5,4,3)
\end{tabular}
\end{center}
\caption{Density contours of the double-Mach shock reflection problem  with 30 contour lines with range from 1 to 22.5, computed using WENO schemes at
 time $T=0.2$ over the domain $[0,4]\times[0,1]$  with mesh grid of size $1600\times 400$.}
 \label{Figure.Doublemach2}
\end{figure}
\begin{example}\label{ex:exp} {\rm   (Explosion Problem) This test case   \cite{tor_90a} involves  two constant states of flow variables 
separated with a circle of radius $R=0.4$ centered at $(1,1)$ over a square domain $[0,2]\times [0,2]$. The initial conditions are given as
\begin{equation}\label{explosion}
(\rho, u,v, p)|_{t=0}=\left\{\begin{array}{ll}
                (1, 0, 0, 1), ~~~~~~~~~  \mbox{If} ~~\sqrt{x^2+y^2} <0.4, \\
                (0.125, 0, 0, 0.1), ~~~~~~~  \mbox{else}.  \\
                \end{array}\right.
\end{equation}
Numerical solutions are computed using WENO schemes at time $T=0.25$ over a uniform grid of resolution $200\times 200$. In Figure \ref{Figure.explosion},
we have shown the cross sectional slice of densities along the plane $y=0$ computed using WENO schemes. The numerical results are compared with the reference solution,
 which is computed using WENO-JS scheme over a uniform mesh of resolution $1000\times 1000$.  The WENO-AO(5,4,3) scheme performs better than the other schemes and
  it resolves the discontinuities without oscillation. The resolution of WENO-AO(5,3) and WENO-AO(5,4,3) schemes are almost comparable.}
\end{example}

\begin{example}\label{ex:rp}{\rm (Riemann Problem) This 2D Riemann problem is taken from \cite{sch-etal_93a}. The simulation is being done over unit 
square domain $[0,1]\times
[0,1]$, initially involves the constant states of flow variables over the each quadrant which is obtained by dividing unit square using 
lines $x=0.8$ and $y=0.8$. We consider
 the 2D Euler system of equations \eqref{2deuler.crate} with the following initial data
 \begin{equation}
(\rho, u,v, p)|_{t=0}=\left\{\begin{array}{ll}
                (1.5, 0, 0, 1.5), ~~~~~~~~~~~~~~~~~~~  \mbox{If} ~~x\geq 0.8, ~~y\geq 0.8, \\
                (0.5323, 1.206.0, 0, 0.3), ~~~~~~~  \mbox{If} ~~ x<0.8, ~~y\geq 0.8, \\
                (0.138, 1.206, 1.206, 0.029), ~~ \mbox{If} ~~ x <0.8, ~~y\geq 0.8, \\
                (0.5323, 0, 1.206, 0.3), ~~~~~~~~~  \mbox{If} ~~x\geq 0.8, ~~ y\leq 0.8.
                \end{array}\right.
\end{equation}
with Dirichlet boundary conditions. The numerical solutions are computed using WENO schemes at time $T=0.8$ with a grid of size $800\times 800$.
 A closer look of Figure \ref{riemann.2d} reveals that WENO schemes of adaptive order and WENO-ZQ yields better solution of complex structures in
 comparison to 
 WENO-JS and WENO-Z schemes. The  resolution of solution using WENO-AO(5,4,3) scheme is better than WENO-AO(5,3) scheme.
 }
\end{example}

\begin{example}\label{ex:kh}{\rm (Kelvin-Helmholtz instability)
 Consider the  system of Euler equations \eqref{2deuler.crate} with the initial conditions  over the periodic domain 
 $[0,1]\times [0,1]$ as follows  \cite{sam-kar_15a}
\begin{equation*}
  p=2.5
 \end{equation*}
 \begin{equation*}
\rho(x,y)=\left\{\begin{array}{ll}
                2, ~~~~  \mbox{If} ~~0.25 < y \leq 0.75, \\
                1, ~~~~ \mbox{else}, 
                \end{array}\right.
\end{equation*}
 \begin{equation*}
u(x,y)=\left\{\begin{array}{ll}
                0.5, ~~~~  \mbox{If} ~~0.25 < y \leq 0.75, \\
                -0.5, ~~ \mbox{else}, 
                \end{array}\right.
\end{equation*}
\begin{equation*}
 v(x,y)=w_0\sin(4\pi x)\left\{\exp\left[-\frac{(y-0.25)^2}{2\sigma^2}\right]+\exp\left[-\frac{(y-0.25)^2}{2\sigma^2}\right]\right\},
\end{equation*}
with $w_0=0.1$, $\sigma=0.05/\sqrt{2}$, and the adiabatic constant  $\gamma$ chosen to be $7/5$. The computational domain is discretized with 
512 cells in each direction, and final time is taken to be $T=0.8$. The numerical solutions are computed using the considered WENO schemes and depicted in 
 Figure \ref{Figure.kh}. We can observe that the  WENO schemes capture the complex structures and are able to capture small-scale vortices. The 
  WENO-AO(5,4,3) scheme resolves more vortices than WENO-AO(5,3) scheme, whereas WENO-AON(5,3) is comparable with other schemes.}
\end{example}

\begin{example}\label{ex:rt} {\rm (Rayleigh-Taylor instability)
  This kind of instability arises on an interface between  two fluids of different densities when an acceleration is directed from
 from heavier fluid to lighter fluid 
   \cite{shi-etal_03a}. Consider the Euler system of equations \eqref{2deuler.crate} with the following initial conditions 
   \begin{equation}
(\rho, u,v, p)(x,y,0)=\left\{\begin{array}{ll}
                (2,0,-0.025a\cos(8\pi x),2y+1), ~~~~~~   0\leq y \leq \frac{1}{2}, \\
                (1,0,-0.025a\cos(8\pi x),y+\frac{3}{2}), ~~~~~~~  \frac{1}{2}\leq y\leq 1.
                \end{array}\right.
\end{equation}
over the computational domain $[0, 1/4]\times [0,1]$. Reflective boundary conditions are imposed on the right and left boundaries. On the upper boundary,
we assign the values $(\rho, u,v, p)=(1,0,0,2.5)$ and for
the bottom boundary, we take 
$(\rho, u,v, p)=(2,0,0,1)$. The source term $Q=(0,0,\rho, \rho v)$ is added to the Euler equations \ref{2deuler.crate}. The computational domain is
discretized with a uniform mesh of resolution $200\times 800$ and  numerical solutions are computed at final time $T=1.95$.  The value of adiabatic constant
 $\gamma$ is taken to be 5/3. In Figure \ref{Figure.rt}, we 
 have shown the density contour plots computed with considered WENO schemes. The result of WENO-AON(5,3) is comparable with WENO-JS and WENO-Z schemes,
  whereas resolving power of  WENO-AO(5,4,3) scheme in resolving  the small-scale vortical structure is better than WENO-AO(5,3) scheme, since we can 
  observe more vortices in density plot of WENO-AO(5,4,3) scheme. The WENO-ZQ scheme is also capable of resolving more complex structures.
 }
\end{example}

\begin{example}\label{ex:dmr}{\rm (Double Mach Reflection)
Now we compare solutions obtained using WENO schemes in case of double-mach shock reflection problem \cite{woo-col_84a}. 
The double-Mach shock reflection
 problem is an important test case where a vertical shock moves horizontally into a wedge that is inclined by some angle. Numerical experiments 
 are performed
  over the domain $[0,4]\times [0,1]$ and solutions are computed at final time $T=0.2$, keeping the CFL number 0.3.  Here, we consider the Euler 
  system of equations \eqref{2deuler.crate} with the initial condition
  }
  \end{example}
 \begin{equation*}
(\rho, u,v, p)(x,y,0)=\left\{\begin{array}{ll}
                (8.0,8.25\cos(\frac{\pi}{6}),-8.25\sin(\frac{\pi}{6}), 116.5) ~~~~   x<x_0+\frac{y}{\sqrt{3}}, \\
                (1.4,0.0,0.0,1.0), ~~~~~~~~~~~~~~~~~~~~~~~~~~~~~  x\geq x_0+\frac{y}{\sqrt{3}},
                \end{array}\right.
\end{equation*}
where $x_0=\frac{1}{6}$. Inflow boundary conditions are imposed at $x=0$ using the post-shock value as above and outflow boundary conditions are implemented
using $\frac{\partial \bold{V}}{\partial x}=0$ at the  right boundary i.e $x=4$. At lower boundary, reflecting boundary conditions are applied to the
interval $[x_0,4]$ and $(\rho_y, u_y,v,p_y)=(0,0,0,0)$ for $[x_0,4]$. The position of shock wave at time $t$ on the upper boundary ($y=1$) is given by
 $s(t)=x_0+\frac{(1+20t)}{\sqrt{3}}$. The boundary conditions on the upper boundary can be implemented using shock position as follows
\begin{equation*}
(\rho, u,v, p)(x,1,t)=\left\{\begin{array}{ll}
                (8.0,8.25\cos(\frac{\pi}{6}),-8.25\sin(\frac{\pi}{6}), 116.5) ~~~~   0\leq x<s(t), \\
                (1.4,0.0,0.0,1.0), ~~~~~~~~~~~~~~~~~~~~~~~~~~~~~  s(t)\leq x\leq 4.
                \end{array}\right.
\end{equation*}
In Figure \ref{Figure.Doublemach1}, we have depicted the numerical densities obtained using WENO-AO(5,3), WENO-AON(5,3), and WENO-AO(5,4,3) schemes for
 the domain $[0,3]\times [0,1]$ over uniform mesh having resolution $1600\times 400$. Further, 
in Figure \ref{Figure.Doublemach2}, we have shown the region around the double Mach stems for all six considered WENO schemes. The resolving power of 
 WENO-AO(5,4,3) scheme is better than the other considered WENO schemes as we can see WENO-AO(5,4,3) scheme captures more number of small vortices 
 along the slip lines than the other schemes. The  vortices formed in WENO-AO(5,3) and WENO-AO(5,4,3) are bigger in size than that of other schemes.
  The resolving power of WENO-AON(5,3) are comparable with WENO-Z and WENO-ZQ schemes and better than WENO-JS scheme. 
 \begin{table}
\centering 
\small
\begin{tabular}{|c| c|c| c| c| c|  c| c| c|}
\hline 
Test                  & WENO-AON(5,3) & WENO-JS & WENO-Z   & WENO-ZQ & WENO-AO(5,4,3) \\
\hline
Example \ref{ex:sv}   &  0.9747      &0.8937   & 0.9187   &0.9664   & 1.0728\\ 
\hline
Example \ref{ex:exp}  &  0.9769      &0.8773   & 0.9063   &0.9602   & 1.0644\\ 
\hline 
Example \ref{ex:rp}   &  0.9714      &0.8733   & 0.9020   &0.9502   & 1.0796   \\
\hline
Example \ref{ex:kh}   &  0.9677      &0.8674   & 0.8911   &0.9639   & 1.0913  \\
\hline
Example \ref{ex:rt}   &  0.9587      &0.8691   & 0.8883   &0.9484   & 1.0549  \\
\hline
Example \ref{ex:dmr}  &  0.9507      &0.8838   & 0.8906   &0.9108   & 1.0819  \\
\hline
\end{tabular}
\caption{Relative computational cost for 2D-Examples with considered WENO schemes. }
\label{tab:timecomp}
\end{table}

\subsection{Computational cost comparison for 2D-Problems}
In this subsection, we compare the relative computational costs of numerical schemes with respect to WENO-AO(5,3) scheme for 2D system of Euler equations
with test cases \ref{ex:sv}-\ref{ex:dmr} taken in subsection \ref{sec:shock}. We have kept all other parameters such as mesh width, CFL etc
same as is given in test cases in 
subsection \ref{sec:shock}.  In Table \ref{tab:timecomp}, we depict the time taken to compute the solution by WENO schemes in relation to the base scheme 
WENO-AO(5,3). We can easily observe that 
WENO-JS and WENO-Z take almost $89-90\%$ time of computational cost of  WENO-AO(5,3) scheme, whereas WENO-AON(5,3) and WENO-ZQ  take $3-4\% $ less time
as compared to the base scheme.
We can also observe that  WENO-AO(5,4,3)
takes $6-9\%$ more time as compared to the base scheme. Since it includes an extra polynomial in the reconstruction than other schemes.

\section{Conclusion}\label{ccl}
In this article, two new algorithms of fifth order WENO scheme of adaptive order for hyperbolic conservation laws, named
 as WENO-AON(5,3) and
WENO-AO(5,4,3) are proposed. In the first algorithm WENO-AON(5,3), we have proposed a new simple smoothness indicator for the bigger stencil using linear 
combination of lower order 
smoothness indicators. 
With the aid of
numerical experiments,
 we have shown that WENO-AON(5,3) scheme is comparable with WENO-AO(5,3)\cite{bal-etal_16a} scheme in terms of accuracy and resolution of
 solution across 
 shocks and discontinuities with less computational cost.
 In the WENO-AO(5,4,3) scheme, we have used a convex combination of three quadratic polynomials, 
 one cubic polynomial and a fourth degree polynomial and
 we found that resolution of solutions obtained using WENO-AO(5,4,3) is better than that from WENO-JS, WENO-Z, WENO-ZQ, and WENO-AO(5,3) 
schemes.
  We have observed that WENO-AO(5,4,3) scheme takes $6-9\%$ more computational time  than WENO-AO(5,3) scheme but is more accurate, and the cost can be reduced using the smoothness indicator
   proposed in the first algorithm. 
 The present algorithms can be extended to WENO scheme of adaptive order more than five and will be a future topic of interest. 
The main advantage of these WENO schemes is that they don't require the existence
 of positive optimal linear weights. 
 Hence the cases, where existence of positive linear optimal weights are not assured, can be solved with these new WENO schemes of adaptive order in a 
 stable manner. 
% \section*{Acknowledgements} 

\bibliographystyle{IMANUM-BIB}
%\bibliography{IMANUM-refs}

\bibliography{comp_bib}

\end{document}